\documentclass[global,envcountsect,envcountsame]{svjour}
\journalname{AAECC}
\usepackage{times,url}
\usepackage{amsmath,amssymb,bbm}
\usepackage{ifthen,calc,float}
\usepackage{algorithm}
\usepackage{algorithm,algorithmic}
   \newcommand{\algorithmicreturn}{\algorithmicfont{return}\ }
   \newcommand{\algorithmicfont}[1]{\textbf{#1}}
\usepackage{fancybox}
    
\usepackage[arrow,graph]{xy}
\smartqed

\def\A{\mathcal{A}}
\def\B{\mathcal{B}}
\def\bq#1#2{\beta_{#1}^{(#2)}}
\def\bull{\hfill$\lhd$}
\def\C{\mathcal{C}}
\def\ccls#1#2{\mathrm{ccls}_{#2}{#1}}
\def\cnum#1{^{\quad\ovalbox{$\scriptstyle #1$}}}
\def\cnumt#1{{\scriptstyle\ovalbox{$#1$}}}
\def\cone#1#2{\C_{#1}(#2)}

\def\drl{\prec_{\mbox{\scriptsize degrevlex}}}
\def\ellv{\vec{\ell}}
\def\ev{\vec{e}}
\def\F{\mathcal{F}}
\def\fv{\vec{f}}
\def\G{\mathcal{G}}
\def\H{\mathcal{H}}
\def\hv{\vec{h}}
\def\I{\mathcal{I}}
\def\Isat{\I^{\mathrm{sat}}}
\def\idiv#1{\,|_{#1}\,}
\def\ispan#1#2{\lspan{#1}_{#2}}
\def\J{\mathcal{J}}

\def\kk{{\mathbbm k}}
\def\kv{\vec{k}}
\def\lc#1#2{\mathrm{lc}_{#1}#2}
\def\lcm#1#2{\mathrm{lcm}(#1,#2)}
\def\le#1#2{\mathrm{le}_{#1}#2} 

\def\lspan#1{\langle{#1}\rangle}
\def\lt#1#2{\mathrm{lt}_{#1}#2}
\def\M{\mathcal{M}}
\def\mult#1#2#3{{#1}_{#2}(#3)}

\def\NN{{\mathbbm N}}
\def\N{\mathcal{N}}
\def\Nn{\NN_0^n}
\def\nmult#1#2#3{\bar{#1}_{#2}(#3)}
\def\norm#1{|{#1}|}
\def\P{\mathcal{P}}
\def\pf{\mathfrak{p}}
\def\qf{\mathfrak{q}}
\def\QQ{{\mathbbm Q}}
\def\R{\mathcal{R}}
\def\RR{{\mathbbm R}}
\def\Rv{\vec{R}}
\def\S{\mathcal{S}}
\def\Sv{\vec{S}}
\def\sv{\vec{s}}
\def\syz#1#2{\mathrm{Syz}^{#1}(#2)}

\def\Spolym#1#2#3{\Sv_{#1}(#2,#3)}
\def\T{\mathcal{T}}
\def\TT{\mathbbm{T}}
\def\tv{\vec{t}}
\def\V{\mathcal{V}}
\def\W{\mathcal{W}}

\def\Y{\mathcal{Y}}
\def\Z{\mathcal{Z}}
\def\xv{\vec{x}}
\def\yv{\vec{y}}
\def\zv{\vec{z}}

\DeclareMathOperator{\ann}{Ann}
\DeclareMathOperator{\ass}{Ass}
\DeclareMathOperator{\chr}{char}
\DeclareMathOperator{\cls}{cls}
\DeclareMathOperator{\depth}{depth}
\DeclareMathOperator{\sdepth}{sdepth}
\DeclareMathOperator{\pd}{proj\,dim}
\DeclareMathOperator{\reg}{reg}
\DeclareMathOperator{\sat}{sat}
\DeclareMathOperator{\supp}{supp}

\begin{document}
\title{A Combinatorial Approach to Involution and $\delta$-Regularity II:
  Structure Analysis of Polynomial Modules with Pommaret Bases}
\titlerunning{Involution and $\delta$-Regularity II}
\author{Werner M.  Seiler}
\institute{AG ``Computational Mathematics'',
  Universit\"at Kassel, 34132 Kassel, Germany\\   
  \url{www.mathematik.uni-kassel.de/~seiler}\\ 
  \email{seiler@mathematik.uni-kassel.de}} 
\date{Received: date / Revised version: date}
\maketitle
\begin{abstract}
  Much of the existing literature on involutive bases concentrates on their
  efficient algorithmic construction.  By contrast, we are here more concerned
  with their structural properties.  Pommaret bases are particularly useful in
  this respect.  We show how they may be applied for determining the Krull and
  the projective dimension, respectively, and the depth of a polynomial
  module.  We use these results for simple proofs of Hironaka's criterion for
  Cohen-Macaulay modules and of the graded form of the Auslander-Buchsbaum
  formula, respectively.
  
  Special emphasis is put on the syzygy theory of Pommaret bases and its use
  for the construction of a free resolution.  In the monomial case, the
  arising complex always possesses the structure of a differential algebra and
  it is possible to derive an explicit formula for the differential.
  Furthermore, in this case one can give a simple characterisation of those
  modules for which our resolution is minimal.  These observations generalise
  results by Eliahou and Kervaire.
  
  Using our resolution, we show that the degree of the Pommaret basis with
  respect to the degree reverse lexicographic term order is just the
  Castelnuovo-Mumford regularity.  This approach leads to new proofs for a
  number of characterisations of this regularity proposed in the literature.
  This includes in particular the criteria of Bayer/Stillman and
  Eisenbud/Goto, respectively.  We also relate Pommaret bases to the recent
  work of Bermejo/Gimenez and Trung on computing the Castelnuovo-Mumford
  regularity via saturations.
  
  It is well-known that Pommaret bases do not always exist but only in
  so-called $\delta$-regular coordinates.  Fortunately, generic coordinates
  are $\delta$-regular.  We show that several classical results in commutative
  algebra, holding only generically, are true for these special coordinates.
  In particular, they are related to regular sequences, independent sets of
  variables, saturations and Noether normalisations.  Many properties of the
  generic initial ideal hold also for the leading ideal of the Pommaret basis
  with respect to the degree reverse lexicographic term order.  We further
  present a new simple criterion for detecting $\delta$-singularity leading to
  a deterministic approach for finding $\delta$-regular coordinates that is
  more efficient than all solutions proposed in the literature so far.
\end{abstract}

\section{Introduction}

Rees~\cite{rees:polymod} introduced a combinatorial decomposition of finitely
generated polynomial modules and related it for graded modules to the Hilbert
series.  Later, more general decompositions of $\kk$-algebras were studied by
Stanley and several other authors (see e.\,g.\ 
\cite{bg:comb,bcs:stanley,rps:hilbert,rps:diophant}), especially in the
context of Cohen-Macaulay complexes but also for other applications like
invariant theory or the theory of normal forms of vector fields with nilpotent
linear part.  Sturmfels and White~\cite{sw:comb} presented algorithms to
compute various combinatorial decompositions.

Apparently all these authors have been unaware that similar decompositions are
implicitly contained in the Janet-Riquier theory of differential equations.
In fact, they represent the fundamental idea underlying this theory.
Involutive bases combine this idea with concepts from Gr\"obner bases.  As we
have seen in Part I, one may consider involutive bases as those Gr\"obner
bases which automatically induce a combinatorial decomposition of the ideal
they generate.

The main goal of this second part is to show that Pommaret bases possess a
number of special properties not shared by other involutive bases which makes
them particularly useful for the structure analysis of polynomial modules.  A
number of important invariants can be directly read off a Pommaret basis.  One
reason for this is that Pommaret bases induce the special type of
decomposition introduced by Rees~\cite{rees:polymod} and which now carries his
name.

Pommaret bases (for ideals in power series rings) implicitly appeared already
in the classical work of Hironaka \cite{hh:resol} on the resolution of
singularities.  Later, Amasaki \cite{ma:applweier,ma:gener} followed up this
idea and explicitly introduced them under the name Weierstra\ss\ bases because
of their connection to the Weierstra\ss\ Preparation Theorem.  In his study of
their properties, Amasaki obtained to some extent similar results than we
present in this article, however in a different way.

This second part is organised as follows.  Section~\ref{sec:pombas} discusses
the problem of $\delta$-regularity and presents a very simple and effective
criterion for detecting $\delta$-singular coordinate systems based on a
comparison of the Janet and Pommaret multiplicative variables.  The proof of
this criterion furthermore implies a method for the automatic construction of
$\delta$-regular coordinates without destroying too much sparsity.  As a first
application, we determine the depth of a polynomial ideal $\I$ and a maximal
$\I$-regular sequence.

The following section studies combinatorial decompositions of general
polynomial modules using involutive bases.  A trivial application, already
given by Janet \cite{ja:lec} and Stanley \cite{rps:hilbert}, is the
determination of the Hilbert series and thus of the Krull dimension.  For
Pommaret bases an alternative characterisation of the dimension can be given
which is related to Gr\"obner's approach via maximal independent sets of
variables \cite{wg:mag,kw:dim}.  Extending our previous results on the depth
from submodules to arbitrary polynomial modules, we obtain as a simple
corollary Hironaka's criterion for Cohen-Macaulay modules.

Section~\ref{sec:noether} discusses the relation between $\delta$-regularity
and Noether normalisation.  It turns out that searching $\delta$-regular
coordinates for an ideal $\I$ is equivalent to putting simultaneously $\I$ and
all primary components of $\lt{\prec}{\I}$ into Noether position.  As a
by-product we provide a number of equivalent characterisations for monomial
ideals possessing a Pommaret basis and show how an irredundant primary
decomposition of such ideals can be easily obtained.  These results are
heavily based on recent work by Bermejo and Gimenez \cite{bg:scmr}.

Section~\ref{sec:syz} develops the syzygy theory of involutive bases.  We show
that the involutive standard representations of the non-multiplicative
multiples of the generators induce a Gr\"obner basis (for an appropriately
chosen term order) of the first syzygy module.  Essentially, this involutive
form of Schreyer's theorem follows from the ideas behind Buchberger's second
criterion for redundant $S$-polynomials.  For Janet and Pommaret bases the
situation is even better, as the arising Gr\"obner basis is then again a Janet
and Pommaret basis, respectively.

In the next three sections we construct by iteration of this result free
resolutions of minimal length.  We first outline the construction for
arbitrary polynomial modules with a Pommaret basis.  Then we specialise to
monomial modules where one can always endow the resolution with the structure
of a differential algebra and provide an explicit formula for the
differential.  In the special case of a stable monomial module the resolution
is even minimal.  Most of these results are inspired by and generalisations of
the work of Eliahou and Kervaire \cite{ek:res}.

In Section \ref{sec:cmr} we show that the degree of a Pommaret basis with
respect to the degree reverse lexicographic order equals the
Castelnuovo-Mumford regularity of the ideal.  Together with our results on the
construction of $\delta$-regular coordinates, this leads to a simple effective
method for the computation of this important invariant.  As corollaries we
recover characterisations of the Castelnuovo-Mumford regularity previously
proposed by Bayer/Stillman \cite{bs:mreg} and Eisenbud/Goto \cite{eg:freeres}.
In the following section we discuss the relation between regularity and
saturation from the point of view of Pommaret bases.  Here we make contact
with recent works of Trung \cite{ngt:cmr} and Bermejo/Gimenez \cite{bg:scmr}.

Finally, we apply the previously developed syzygy theory to the construction
of involutive bases in iterated polynomial algebras of solvable type.  A
rather technical appendix clarifies the relation between Pommaret bases and
the Sturmfels-White approach \cite{sw:comb} to the construction of Rees
decompositions.

\section{Pommaret Bases and $\delta$-Regularity}\label{sec:pombas}

We saw in Part I (Example~2.12) that not every monoid ideal in $\Nn$ possesses
a finite Pommaret basis: the Pommaret division is not Noetherian.  Obviously,
this is also implies that there are polynomial ideals
$\I\subseteq\P=\kk[x_1,\dots,x_n]$ without a finite Pommaret basis for a given
term order.  However, we will show now that at the level of \emph{polynomial}
ideals this may be considered as solely a problem of the chosen variables
$\xv$.  For this purpose, we take in the sequel the following point of view:
term orders are defined for exponent vectors, i.\,e.\ on the monoid $\Nn$;
performing a linear change of variables $\tilde{\xv}=A\xv$ leads to new
exponent vectors in each polynomial which are then sorted according to the
same term order as before.

\begin{definition}\label{def:dregid}
  The variables\/ $\xv$ are\/ \emph{$\delta$-regular} for the ideal\/
  $\I\subseteq\P$ and the term order\/ $\prec$, if\/ $\I$ possesses a finite
  Pommaret basis for\/ $\prec$.
\end{definition}

Given our definition of an involutive basis, it is obvious that
$\delta$-regularity concerns the existence of a Pommaret basis for the monoid
ideal $\le{\prec}{\I}$.  A coordinate transformation generally yields a new
leading ideal which may possess a Pommaret basis.  In fact, we will show in
this section that for every polynomial ideal $\I\subseteq\P$ variables $\xv$
exist such that $\I$ has a finite Pommaret basis provided that the chosen term
order $\prec$ is class respecting.\footnote{Recall from the appendix of Part I
  that any class respecting term order coincides on terms of the same degree
  with the reverse lexicographic order.}

Besides showing the mere existence of $\delta$-regular variables, we will
develop in this section an effective approach to recognising $\delta$-singular
coordinates and transforming them into $\delta$-regular ones.  It is inspired
by the work of Gerdt \cite{vpg:janpom} on the relation between Pommaret and
Janet bases and the key ideas have already been used in the context of the
combined algebraic-geometric completion to involution of linear differential
equations \cite{wms:geocompl}.  We begin by proving two useful technical
lemmata.  The number $\max_{h\in\H}{\deg{h}}$ is called the \emph{degree} of
the finite set $\H\subset\P$ and denoted by $\deg{\H}$.

\begin{lemma}\label{lem:truncbas}
  Let the set\/ $\H$ be a homogeneous Pommaret basis of the homogeneous
  ideal\/ $\I\subseteq\P$.  Then for any degree\/ $q\geq\deg{\H}$ a Pommaret
  basis of the truncated ideal\/ $\I_{\geq q}=\bigoplus_{p\geq q}\I_p$ is
  given by
  \begin{equation}\label{eq:truncbas}
    \H_q=\bigl\{x^\mu h\mid h\in\H,\ |\mu|+\deg{h}=q,\
                            \forall j>\cls{h}:\mu_j=0\bigr\}\;.
  \end{equation}
  Conversely, if\/ $\I_{\geq q}$ possesses a finite Pommaret basis, then so
  does\/ $\I$.
\end{lemma}

\begin{proof}
  According to the conditions in (\ref{eq:truncbas}), each polynomial $h\in\H$
  is multiplied by terms $x^\mu$ containing only variables which are
  multiplicative for it.  Thus trivially $\cls{(x^\mu h)}=\cls{\mu}$.
  Furthermore, $\H_q$ is involutively head autoreduced, as $\H$ is.  Now let
  $f\in\I_{\geq q}$ be an arbitrary homogeneous polynomial.  As $\H$ is a
  Pommaret basis of $\I$, it has a standard representation
  $f=\sum_{h\in\H}P_hh$ with polynomials $P_h\in\kk[x_1,\dots,x_{\cls{h}}]$.
  Hence $f$ can be written as a linear combination of polynomials $x^\nu h$
  where $|\nu|=\deg{f}-\deg{h}\geq q-\deg{h}$ and where $x^\nu$ contains only
  multiplicative variables.  We decompose $\nu=\mu+\rho$ with
  $|\mu|=q-\deg{h}$ and $\rho_j=0$ for all $j>\cls{\mu}$.  Thus $x^\nu
  h=x^\rho(x^\mu h)$ with $x^\mu h\in\H_q$ and $x^\rho$ contains only
  variables multiplicative for it.  But this trivially implies the existence
  of a standard representation $f=\sum_{h'\in\H_q}P_{h'}h'$ with
  $P_{h'}\in\kk[x_1,\dots,x_{\cls{h'}}]$ and thus $\H_q$ is a Pommaret basis
  of $\I_{\geq q}$.
  
  The converse is also very simple.  Let $\H_q$ be a finite Pommaret basis of
  the truncated ideal $\I_{\geq q}$ and $\H_p$ head autoreduced $\kk$-linear
  bases of the components $\I_p$ for $0\leq p<q$.  If we set
  $\H=\bigcup_{p=0}^q\H_q$, then $\le{\prec}{\H}$ is obviously a weak Pommaret
  basis of the full monoid ideal $\le{\prec}{\I}$ and by Proposition~5.7 of
  Part I an involutive head autoreduction yields a strong basis.\qed
\end{proof}

\begin{lemma}\label{lem:stable}
  With the same notations as in Lemma~\ref{lem:truncbas}, let\/
  $\N=\le{\prec}{\H_q}$.  If\/ $\nu\in\N$ with\/ $\cls\nu=k$,
  then\/\footnote{Recall from Part I that $\ell_i$ denotes for any number
    $\ell\in\NN$ the multi index where all entries except the $i$th one vanish
    and the $i$th one is given by $\ell$.}  $\nu-1_k+1_j\in\N$ for all\/
  $k<j\leq n$.  Conversely, let\/ $\N\subseteq(\Nn)_q$ be a set of multi
  indices of degree\/ $q$.  If for each\/ $\nu\in\N$ with\/ $\cls\nu=k$ and
  each\/ $k<j\leq n$ the multi index\/ $\nu-1_k+1_j$ is also contained in\/
  $\N$, then the set\/ $\N$ is involutive for the Pommaret division.
\end{lemma}

\begin{proof}
  $j$ is non-multiplicative for $\nu$.  As $\N$ is an involutive basis of
  $\le{\prec}{\I_{\geq q}}$, it must contain a multi index $\mu$ with
  $\mu\idiv{P}\nu+1_j$.  Obviously, $\cls{(\nu+1_j)}=k$ and thus $\cls\mu\geq
  k$.  Because of $|\mu|=|\nu|$, the only possibility is
  $\mu=\nu+1_j-1_k$.  The converse is trivial, as each non-multiplicative
  multiple of $\nu\in\N$ is of the form $\nu+1_j$ with $j>k=\cls{\nu}$ and
  hence has $\nu-1_k+1_j$ as an involutive divisor.\qed
\end{proof}

As in practice one is defining an ideal $\I\subseteq\P$ by some finite
generating set $\F\subset\I$, we introduce a concept of $\delta$-regularity
for such sets.  Assume that $\F$ is involutively head autoreduced with respect
to an involutive division $L$ and a term order $\prec$; then we call the total
number of multiplicative variables of the elements of $\F$ its
\emph{involutive size} and denote it by
\begin{equation}\label{eq:invsize}
  \norm{\F}_{L,\prec}=\sum_{f\in\F}\norm{\mult{X}{L,\prec,\F}{f}}\;.
\end{equation}

Let $\tilde{\xv}=A\xv$ be a linear change of coordinates with a regular matrix
$A\in\kk^{n\times n}$.  It transforms each polynomial $f\in\P$ into a
polynomial $\tilde f\in\tilde{\P}=\kk[\tilde x_1,\dots,\tilde x_n]$ of the
same degree.  Thus $\F$ is transformed into a set
$\tilde{\F}\subset\tilde{\P}$ which generally is no longer involutively head
autoreduced.  Performing an involutive head autoreduction yields a set
$\tilde{\F}^\triangle$.  The leading exponents of $\tilde{\F}^\triangle$ may
be very different from those of $\F$ and thus $\norm{\F}_{L,\prec}$ may differ
from $\norm{\tilde{\F}^\triangle}_{L,\prec}$.

\begin{definition}\label{def:dregalg}
  Let the finite set\/ $\F\subset\P$ be involutively head autoreduced with
  respect to the Pommaret division and a term order\/ $\prec$.  The
  coordinates\/ $\xv$ are called\/ \emph{$\delta$-regular} for\/ $\F$ and\/
  $\prec$, if after any linear change of coordinates\/ $\tilde{\xv}=A\xv$ the
  inequality $\norm{\F}_{P,\prec}\geq\norm{\tilde{\F}^\triangle}_{P,\prec}$
  holds.
\end{definition}

\begin{example}\label{ex:dsing}
  Let us reconsider Example~2.12 of Part I where $\F=\{xy\}\subset\kk[x,y]$.
  Independent of how we order the variables, the class of $xy$ is $1$.  After
  the change of coordinates $x=\tilde x+\tilde y$ and $y=\tilde y$, we obtain
  the set $\tilde{\F}=\{\tilde y^2+\tilde x\tilde y\}\subset\kk[\tilde
  x,\tilde y]$.  If we use the degree reverse lexicographic order, the leading
  term is $\tilde y^2$ which is of class $2$.  Thus it is possible to enlarge
  the involutive size of $\F$ and the original coordinates are not
  $\delta$-regular.  The new coordinates obviously are, as a higher class than
  $2$ is not possible.\bull
\end{example}

Note that generally $\delta$-regularity of the variables $\xv$ for a finite
set $\F$ according to Definition \ref{def:dregalg} and for the ideal
$\I=\lspan{\F}$ according to Definition \ref{def:dregid} are independent
properties.  For a concrete instance where the two definitions differ see
Example \ref{ex:dregalg} below where one easily checks that the used
coordinates are $\delta$-regular for the whole submodule but not for the given
generating set, as any transformation of the form $x=\bar x+a\bar y$ with
$a\neq0$ will increase the involutive size.  The main point is that
$\delta$-regularity for the ideal $\I$ is concerned with the monoid ideal
$\le{\prec}{\I}$ whereas $\delta$-regularity for the set $\F$ depends on the
ideal $\lspan{\le{\prec}{\F}}\subseteq\le{\prec}{\I}$.

\begin{proposition}\label{prop:dreg}
  Let\/ $\H$ be a Pommaret basis of the ideal\/ $\I\subseteq\P$ for the term
  order\/~$\prec$.  Then the coordinates\/ $\xv$ are\/ $\delta$-regular for\/
  $\H$ and\/ $\prec$.
\end{proposition}

\begin{proof}
  By definition, $\le{\prec}{\H}$ is a Pommaret basis of $\le{\prec}{\I}$.
  Thus it suffices to consider homogeneous ideals and, by
  Lemma~\ref{lem:truncbas}, we may further reduce to the case that all
  generators are of the same degree $q$.  The involutive size of $\H_q$ is
  then the vector space dimension of $\I_{q+1}$.  This value is independent of
  the chosen coordinates and no involutively head autoreduced set can lead to
  a larger involutive size.  Thus our coordinates are $\delta$-regular for
  $\H$, too.\qed
\end{proof}

Thus we conclude that $\delta$-regularity is necessary for the existence of
Pommaret bases and from an algorithmic point of view $\delta$-singular
coordinates are ``bad''.  Note that for the completion Algorithm~3 of Part I
$\delta$-regularity for the current basis $\H$ is as important as
$\delta$-regularity for the ideal $\I=\lspan{\H}$.  Even if the variables
$\xv$ are $\delta$-regular for $\I$, it may still happen that the algorithm
does not terminate, as it tries to construct a non-existing Pommaret basis for
$\lspan{\le{\prec}{\H}}$.

\begin{example}\label{ex:dregalg}
  One of the simplest instance where this termination problem occurs is not
  for an ideal but for a submodule of the free $\kk[x,y]$-module with basis
  $\{\ev_1,\ev_2\}$.  Consider the set $\F=\{y^2\ev_1,xy\ev_1+\ev_2,x\ev_2\}$
  and any term order for which $xy\ev_1\succ\ev_2$.  The used coordinates are
  not $\delta$-regular for $\F$, as any transformation of the form $x=\bar
  x+a\bar y$ with $a\neq0$ will increase the involutive size.  Nevertheless,
  the used coordinates are $\delta$-regular for the submodule $\lspan{\F}$.
  Indeed, adding the generator $y\ev_2$ (the $S$-``polynomial'' of the first
  two generators) makes $\F$ to a reduced Gr\"obner basis which is
  simultaneously a minimal Pommaret basis.
  
  Note that the termination of the involutive completion algorithm depends
  here on the precise form of the used term order.  If we have $xy^k\ev_2\prec
  xy^2\ev_1$ for all $k\in\NN$, then the algorithm will not terminate, as in
  the $k$th iteration it will add the generator $xy^k\ev_2$.  Otherwise it
  will treat at some stage the non-multiplicative product
  $y\cdot(xy\ev_1+\ev_2)$ and thus find the decisive generator $y\ev_2$.  This
  is in particular the case for any degree compatible order.\bull
\end{example}

Fortunately, most coordinates are $\delta$-regular for a given set $\F$ of
polynomials.  Choosing an arbitrary reference coordinate system, we may
identify every system of coordinates with the regular matrix $A\in\kk^{n\times
  n}$ defining the linear transformation from our reference system to it.

\begin{proposition}\label{prop:dsing}
  The coordinate systems that are\/ $\delta$-singular for a given finite
  involutively head autoreduced set\/ $\F\subset\P$ form a Zariski closed set
  in\/ $\kk^{n\times n}$.
\end{proposition}

\begin{proof}
  We perform first a linear coordinate transformation with an undetermined
  matrix $A$, i.\,e.\ we treat its entries as parameters.  This obviously
  leads to a $\delta$-regular coordinate system, as each polynomial in
  $\tilde{\F}^\triangle$ will get its maximally possible class.
  $\delta$-singular coordinates are defined by the vanishing of certain
  (leading) coefficients.  These coefficients are polynomials in the entries
  of $A$.  Thus the set of all $\delta$-singular coordinate systems can be
  described as the zero set of an ideal.\qed
\end{proof}

Thus the $\delta$-singular coordinate systems form a variety of dimension less
than $n^2$ in an $n^2$-dimensional space and if we randomly choose a
coordinate system, it is $\delta$-regular with probability $1$.  Nevertheless,
$\delta$-singular coordinate systems exist for most sets $\F$ (exceptions are
sets like $\F=\{x^2+y^2\}$ over $\kk=\RR$ where one easily shows that any
coordinate system is $\delta$-regular), and in applications such coordinates
often show up.  We will see below in the proof of Theorem~\ref{thm:dreg} that
it suffices for the construction of $\delta$-regular coordinates to consider
only transformations $\xv\mapsto A\xv$ where the matrix $A$ satisfies
$A_{ii}=1$ and $A_{ij}=0$ for $j>i$.

Our goal is to develop a simple criterion that coordinates are
$\delta$-singular for a given set $\F$ and a class respecting term order.  The
basic idea is to compare the multiplicative variables assigned by the Pommaret
and the Janet division, respectively.  The definitions of these two divisions
are very different.  Somewhat surprisingly, they yield very similar results,
as shown by Gerdt and Blinkov \cite[Prop.~3.10]{gb:invbas}.

\begin{proposition}\label{prop:janpom}
  Let the finite set\/ $\N\subset\Nn$ be involutively autoreduced with respect
  to the Pommaret division.  Then\/
  $\mult{N}{P}{\nu}\subseteq\mult{N}{J,\N}{\nu}$ for all\/ $\nu\in\N$.
\end{proposition}

For later use, we mention the following simple corollaries which further study
the relationship between the Janet and the Pommaret division.  Recall that any
set $\N\subset\Nn$ is involutively autoreduced with respect to the Janet
division.  We first note that by an involutive autoreduction of $\N$ with
respect to the \emph{Pommaret} division its \emph{Janet} span can become only
larger but not smaller.

\begin{corollary}\label{cor:janpomred}
  Let\/ $\N\subset\Nn$ be an arbitrary finite set of multi indices and set\/
  $\N_P=\N\setminus\bigl\{\nu\in\N\mid\exists\mu\in\N:\mu\idiv{P}\nu\bigr\}$,
  i.\,e.\ we eliminate all multi indices possessing a Pommaret divisor in\/
  $\N$.  Then\/ $\ispan{\N}{J}\subseteq\ispan{\N_P}{J}$.
\end{corollary}

\begin{proof}
  If $\mu^{(1)}\idiv{P}\mu^{(2)}$ and $\mu^{(2)}\idiv{P}\nu$, then trivially
  $\mu^{(1)}\idiv{P}\nu$.  Thus for each eliminated multi index
  $\nu\in\N\setminus\N_P$ another multi index $\mu\in\N_P$ exists with
  $\mu\idiv{P}\nu$.  Let $\cls\mu=k$.  By the proposition above
  $\{1,\dots,k\}\subseteq\mult{N}{J,\N_P}{\mu}$.  Assume that an index $j>k$
  exists with $j\in\mult{N}{J,\N}{\nu}$.  By definition of the Pommaret
  division, $\mu_i=\nu_i$ for all $i>k$.  Thus $\mu\in(\nu_{j+1},\dots,\nu_n)$
  and $j\in\mult{N}{J,\N}{\mu}$.  As by the second condition on an involutive
  division $\mult{N}{J,\N}{\mu}\subseteq\mult{N}{J,\N_P}{\mu}$ for all
  $\mu\in\N_P$, we conclude that $j\in\mult{N}{J,\N_P}{\mu}$ and
  $\cone{J,\N}{\nu}\subset\cone{J,\N_P}{\mu}$.  But this immediately implies
  $\ispan{\N}{J}\subseteq\ispan{\N_P}{J}$.  \qed
\end{proof}

This implies that any Pommaret basis is simultaneously a Janet basis (a
similar result is contained in \cite[Thm.~17]{vpg:janpom}).  Thus if $\H$ is a
Pommaret basis, then $\mult{X}{P,\prec}{h}=\mult{X}{J,\H,\prec}{h}$ for all
polynomials $h\in\H$.

\begin{corollary}\label{cor:janpom}
  Let the finite set\/ $\H\subset\P$ be involutive with respect to the
  Pommaret division (and some term order).  Then\/ $\H$ is also involutive for
  the Janet division.
\end{corollary}

\begin{proof}
  By the proposition above, it is obvious that the set $\H$ is at least weakly
  involutive with respect to the Janet division.  For the Janet division any
  weakly involutive set is strongly involutive, if no two elements have the
  same leading terms.  But as $\H$ is a Pommaret basis, this cannot
  happen.\qed
\end{proof}

We show now that for $\delta$-regular coordinate systems and a class
respecting term order the inclusions in Proposition~\ref{prop:janpom} must be
equalities.  In other words, if a variable $x_\ell$ exists which is
multiplicative for an element of $\F$ with respect to the Janet division but
non-multiplicative with respect to the Pommaret division, then the used
coordinates are $\delta$-singular.  Our proof is constructive in the sense
that it shows us how to find coordinates leading to a larger involutive span.

\begin{theorem}\label{thm:dreg}
  Let the finite set\/ $\F\subset\P$ be involutively head autoreduced for the
  Pommaret division and a class respecting term order\/ $\prec$.  Furthermore
  assume that the field\/ $\kk$ contains more than\/ $|\F|\deg{\F}$ elements.
  If\/ $\norm{\F}_{J,\prec}>\norm{\F}_{P,\prec}$, then the coordinates\/ $\xv$
  are\/ $\delta$-singular for\/ $\F$.
\end{theorem}

\begin{proof}
  By the proposition above, we have $\mult{X}{P}{f}\subseteq\mult{X}{J,\F}{f}$
  for all $f\in\F$.  Assume that for a polynomial $h\in\F$ the strict
  inclusion $\mult{X}{P}{h}\subset\mult{X}{J,\F}{h}$ holds.  Thus at least one
  variable $x_\ell\in\mult{X}{J,\F}{h}$ with $\ell>k=\cls{h}$ exists.  We
  perform the linear change of variables $x_i=\tilde x_i$ for $i\neq k$ and
  $x_k=\tilde x_k+a\tilde x_\ell$ with a yet arbitrary parameter
  $a\in\kk\setminus\{0\}$.  This induces the following transformation of the
  terms $x^\mu\in\TT$:
  \begin{equation}\label{eq:termtrafo}
    x^\mu=\sum_{j=0}^{\mu_k}\binom{\mu_k}{j}a^j\tilde x^{\mu-j_k+j_\ell}\;.
  \end{equation}
  Let $\le{\prec}{h}=\mu$.  Thus $\mu=[0,\dots,0,\mu_k,\dots,\mu_n]$ with
  $\mu_k>0$.  Consider the multi index $\nu=\mu-(\mu_k)_k+(\mu_k)_\ell$;
  obviously, $\cls{\nu}>k$.  Applying our transformation to the polynomial $h$
  leads to a polynomial $\tilde h$ containing the term $\tilde x^\nu$.  Note
  that $\nu$ cannot be an element of $\le{\prec}{\F}$.  Indeed, if it was, it
  would be an element of the same set $(\mu_{\ell+1},\dots,\mu_n)$ as $\mu$.
  But this contradicts our assumption that $\ell$ is multiplicative for the
  multi index $\mu$ with respect to the Janet division, as by construction
  $\nu_\ell>\mu_\ell$.
  
  Transforming all polynomials $f\in\F$ yields the set $\tilde{\F}$ on which
  we perform an involutive head autoreduction in order to obtain the set
  $\tilde{\F}^\triangle$.  Under our assumption on the size of the ground field
  $\kk$, we can choose the parameter $a$ such that after the transformation
  each polynomial $\tilde f\in\tilde{\F}$ has at least the same class as the
  corresponding polynomial $f\in\F$, as our term order respects classes.  This
  is a simple consequence of (\ref{eq:termtrafo}): cancellations of terms may
  occur only, if the parameter $a$ is a zero of some polynomial (possibly one
  for each member of $\F$) with a degree not higher than $\deg{\F}$.  By the
  definition of the Pommaret division, if
  $\le{\prec}{f_2}\idiv{P}\le{\prec}{f_1}$, then
  $\cls{\le{\prec}{f_2}}\geq\cls{\le{\prec}{f_1}}$.  Hence even after the
  involutive head autoreduction the involutive size of $\tilde{\F}^\triangle$
  cannot be smaller than that of $\F$.
  
  Consider again the polynomial $h$.  The leading term of the transformed
  polynomial $\tilde h$ must be greater than or equal to $\tilde x^\nu$.  Thus
  its class is greater than $k$.  This remains true even after an involutive
  head autoreduction with all those polynomials $\tilde f\in\tilde{\F}$ that
  are of class greater than $k$, as $x^\nu\notin\lt{\prec}{\F}$.  Hence the
  only possibility to obtain a leading term of class less than or equal to $k$
  consists of an involutive reduction with respect to a polynomial $\tilde
  f\in\tilde{\F}$ with $\cls{f}\leq k$.  But this implies that
  $\cls{\le{\prec}{\tilde f}}>k$.  So we may conclude that after the
  transformation we have at least one polynomial more whose class is greater
  than $k$.  So the coordinates $\xv$ cannot be $\delta$-regular.\qed
\end{proof}

\begin{corollary}\label{cor:jandreg}
  If the coordinates\/ $\xv$ are\/ $\delta$-regular for the Pommaret head
  autoreduced set\/ $\F$, then\/ $\ispan{\F}{J,\prec}=\ispan{\F}{P,\prec}$ for
  any class respecting term order\/ $\prec$.
\end{corollary}

It is important to note that this corollary provides us only with a necessary
but \emph{not} with a sufficient criterion for $\delta$-regularity.  In other
words, even if the Janet and the Pommaret size are equal for a given set
$\F\subset\P$, this does not imply that the used coordinates are
$\delta$-regular.

\begin{example}\label{ex:jansing}
  Let $\F=\bigl\{\underline{z^2}-y^2-2x^2,\ \underline{xz}+xy,\ 
  \underline{yz}+y^2+x^2\bigr\}$.  The underlined terms are the leaders with
  respect to the degree reverse lexicographic order.  One easily checks that
  the Janet and the Pommaret division yield the same multiplicative variables.
  If we perform the transformation $\tilde x=z$, $\tilde y=y+z$ and $\tilde
  z=x$, then after an autoreduction we obtain the set
  $\tilde\F^\triangle=\bigl\{\underline{\tilde z^2}-\tilde x\tilde y,\ 
  \underline{\tilde y\tilde z},\ \underline{\tilde y^2}\bigr\}$.  Again the
  Janet and the Pommaret division yield the same multiplicative variables, but
  $\norm{\tilde\F^\triangle}_{P,\prec}>\norm{\F}_{P,\prec}$.  Thus the
  coordinates $(x,y,z)$ are not $\delta$-regular for $\F$.
  
  The explanation of this phenomenon is very simple.  Obviously our criterion
  depends only on the leading terms of the set $\F$.  In other words, it
  analyses the monomial ideal $\lspan{\lt{\prec}{\F}}$.  Here
  $\lspan{\lt{\prec}{\F}}=\lspan{xz,yz,z^2}$ and one easily verifies that the
  used generating set is already a Pommaret basis.  However, for
  $\I=\lspan{\F}$ the leading ideal is $\lt{\prec}{\I}=\lspan{x^3,xz,yz,z^2}$
  (one obtains a Janet basis for $\I$ by adding the polynomial $x^3$ to $\F$)
  and obviously it does not possess a finite Pommaret basis, as such a basis
  would have to contain all monomials $x^3y^k$ with $k\in\NN$ (or we exploit
  our criterion noting that $y$ is a Janet but not a Pommaret multiplicative
  variable for $x^3$).  Thus we have the opposite situation compared to
  Example \ref{ex:dregalg}: there $\lt{\prec}{\I}$ had a finite Pommaret basis
  but $\lspan{\lt{\prec}{\F}}$ not; here it is the other way round.  We will
  show later in Proposition~\ref{prop:qstabdreg} that whenever
  $\lspan{\lt{\prec}{\F}}$ does not possess a finite Pommaret basis, then
  $\norm{\F}_{J,\prec}>\norm{\F}_{P,\prec}$.\bull
\end{example}

While we have defined $\delta$-regularity as a problem of the Pommaret
division, this example shows that the Janet division ``feels'' it, too.
Usually it yields in $\delta$-singular coordinates a lower number of
multiplicative variables and thus larger involutive bases.
$\delta$-regularity is a typical phenomenon in commutative algebra: many
results hold only in ``generic'' coordinates, i.\,e.\ from a Zariski open set.

Now we are finally in the position to prove that in suitably chosen
coordinates every ideal possesses a finite Pommaret basis for any term order.
The proof exploits again the close relationship of the Janet and the Pommaret
division and a little trick due to Gerdt \cite{vpg:janpom}.

\begin{theorem}\label{thm:expombas}
  Let\/ $\prec$ be an arbitrary term order and the field\/ $\kk$ sufficiently
  large.  Then every ideal\/ $\I\subseteq\P$ possesses a finite Pommaret basis
  for\/ $\prec$ in suitably chosen variables\/ $\xv$.
\end{theorem}

\begin{proof}
  As a first step we show that every ideal has a \emph{Pommaret} head
  autoreduced \emph{Janet} basis.  Indeed, let us apply our polynomial
  completion algorithm (Algorithm~3 of Part I) for the Janet division with one
  slight modification: in the Lines~/1/ and /9/ we perform the involutive head
  autoreductions with respect to the Pommaret division.  It is obvious that if
  the algorithm terminates, the result is a basis with the wanted properties.
  
  The Janet division is Noetherian (Lemma~2.14 of Part I).  Thus without our
  modification the termination is obvious.  With respect to the Janet division
  every set of multi indices is involutively autoreduced.  Hence a Janet head
  autoreduction only takes care that no two elements of a set have the same
  leading exponents.  But in Line /9/ we add a polynomial that is in
  involutive normal form so that no involutive head reductions are possible.
  As the Pommaret head autoreduction may only lead to a larger monoid ideal
  $\le{\prec}{\H_i}$, the Noetherian argument in the proof of the termination
  of the algorithm remains valid after our modification.
  
  Once the ascending ideal chain
  $\lspan{\le{\prec}{\H_1}}\subseteq\lspan{\le{\prec}{\H_2}}\subseteq\cdots$
  has become stationary, the polynomial completion algorithm essentially
  reduces to the ``monomial'' one (Algorithm~2 of Part I).  According to
  Corollary~\ref{cor:janpomred}, the Pommaret head autoreductions may only
  increase the Janet spans $\ispan{\le{\prec}{\H_i}}{J}$.  Thus the
  termination of the monomial completion is not affected by our modification
  and the algorithm terminates for arbitrary input.
  
  By Proposition~\ref{prop:dsing}, the coordinate systems that are
  $\delta$-singular for a given finite set $\F\subset\P$ form a Zariski closed
  set.  As our modified algorithm terminates, it treats only a finite number
  of sets $\H_i$ and the coordinate systems that are $\delta$-singular for at
  least one monomial set $\lt{\prec}{\H_i}$ still form a Zariski closed set.
  Thus generic coordinates are $\delta$-regular for all sets
  $\lt{\prec}{\H_i}$ and by Corollary~\ref{cor:jandreg}\footnote{Note that it
    is not relevant here that the corollary assumes the use of a class
    respecting term order, since we are arguing only with monomial sets.}
  their Janet and their Pommaret spans coincide.  But this implies that the
  result of the modified algorithm is not only a Janet but also a Pommaret
  basis.\qed
\end{proof}

\begin{remark}
  The assumption on the size of the field $\kk$ is obviously trivially
  satisfied for a field of characteristic zero.  In the case of a finite
  field, it may be necessary to enlarge $\kk$ in order to guarantee the
  existence of a Pommaret basis.  This is similar to the situation when one
  tries to put a zero-dimensional ideal in normal $x_n$-position
  \cite{kr:cca1}.\bull
\end{remark}

\begin{remark}
  In Part I we discussed the extension of the Mora normal form to involutive
  basis computation.  Obviously, the above result remains valid, if we
  substitute the ordinary normal form by Mora's version and hence we may also
  apply it to Pommaret bases with respect to semigroup orders.\bull
\end{remark}

Theorem~\ref{thm:dreg} provides us with an effective mean to detect that the
completion to a Pommaret basis does \emph{not} terminate.  We follow the
polynomial completion algorithm with the Pommaret division but check after
each iteration whether $\norm{\H_i}_{J,\prec}=\norm{\H_i}_{P,\prec}$.  If no
finite Pommaret basis exists in the given coordinates, sooner or later the
Janet size is greater than the Pommaret size.  Indeed, by the considerations
above the completion with respect to the Janet division terminates.  Thus
either the result is simultaneously a Pommaret basis or in some iteration our
criterion must detect a variable which is multiplicative only for the Janet
but not for the Pommaret division.  Alternatively, we compute a Pommaret head
autoreduced Janet basis (which always exists by the considerations above) and
check whether it is simultaneously a Pommaret basis.

\begin{example}\label{ex:pomterm}
  Let us apply this approach to the Pommaret completion of the set
  $\F=\bigl\{\underline{z^2}-y^2-2x^2,\ \underline{xz}+xy,\ 
  \underline{yz}+y^2+x^2\bigr\}$ (with respect to the degree reverse
  lexicographic order).  We have seen in Example~\ref{ex:jansing} that the
  coordinates are not $\delta$-regular for $\I$, although the Janet and the
  Pommaret span of $\F$ coincide.  According to our algorithm we must first
  analyse the polynomial $y(xz+xy)$.  Its involutive normal form with respect
  to $\F$ is $-\underline{x^3}$.  If we determine the multiplicative variables
  for the enlarged set, they do not change for the old elements.  For the new
  polynomial the Janet division yields $\{x,y\}$.  But $y$ is obviously not
  multiplicative for the Pommaret division.  Thus our criterion tells us that
  the Pommaret completion will not terminate.  Indeed it is easy to see that
  no matter how often we multiply the new polynomial by $y$, it will never
  become involutively head reducible and no finite Pommaret basis exists.
  
  In this example, the Janet completion (with or without Pommaret
  autoreductions) ends with the addition of this single obstruction to
  involution and we obtain as Janet basis the set
  \begin{equation}\label{eq:pomterm}
    \F_J=\bigl\{\underline{z^2}-y^2-2x^2,\ \underline{xz}+xy,\ 
                \underline{yz}+y^2+x^2,\ \underline{x^3}\bigr\}\;.
  \end{equation}
  In Example~\ref{ex:jansing} we showed that the transformation $\tilde x=z$,
  $\tilde y=y+z$ and $\tilde z=x$ yields after an autoreduction the set
  $\tilde\F^\triangle=\bigl\{\underline{\tilde z^2}-\tilde x\tilde y,\
  \underline{\tilde y\tilde z},\ \underline{\tilde y^2}\bigr\}$.  One easily
  checks that it is a Pommaret and thus also a Janet basis.  We see again that
  the Janet division also ``feels'' $\delta$-singularity in the sense that in
  such coordinates it typically leads to larger bases of higher degree.\bull
\end{example}

Besides being necessary for the mere existence of a finite Pommaret basis, a
second application of $\delta$-regular coordinates is the construction of
$\I$-regular sequences for a homogeneous ideal $\I\subseteq\P$.  Recall that
for any $\P$-module $\M$ a sequence $(f_1,\dots,f_r)$ of polynomials
$f_i\in\P$ is called $\M$-regular, if the polynomials generate a proper ideal,
$f_1$ is a non zero divisor for $\M$ and each $f_i$ is a non zero divisor for
$\M/\lspan{f_1,\dots,f_{i-1}}\M$.  The maximal length of an $\M$-regular
sequence is the \emph{depth} of the module.  While the definition allows for
arbitrary polynomials in such sequences, it suffices for computing the depth
to consider only linear forms $f_i\in\P_1$.  This follows, for example, from
\cite[Cor.~17.7]{de:ca} or \cite[Lemma~4.1]{sw:comb}.  For this reason, the
following proof treats only this case.

\begin{proposition}\label{prop:Iregular}
  Let\/ $\I\subseteq\P$ be a homogeneous ideal and\/ $\H$ a homogeneous
  Pommaret basis of it for a class respecting term order.  Let\/
  $d=\min_{h\in\H}\cls{h}$.  Then the variables\/ $(x_1,\dots,x_d)$ form a
  maximal\/ $\I$-regular sequence and thus\/ $\depth{\I}=d$.
\end{proposition}

\begin{proof}
  A Pommaret basis $\H$ induces a decomposition of $\I$ of the form
  \begin{equation}\label{eq:pomrees}
    \I=\bigoplus_{h\in\H}\kk[x_1,\dots,x_{\cls{h}}]\cdot h\;.
  \end{equation}
  If $d=\min_{h\in\H}\cls{h}$ denotes the minimal class of a generator in
  $\H$, then (\ref{eq:pomrees}) trivially implies that the sequence
  $(x_1,\dots,x_d)$ is $\I$-regular.
  
  Let us try to extend this sequence by a variable $x_k$ with $k>d$.  We
  introduce $\H_d=\{h\in\H\mid\cls{h}=d\}$ and choose an element $\bar
  h\in\H_d$ of maximal degree.  As we use a class respecting order, $\bar
  h\in\lspan{x_1,\dots,x_d}$ by Lemma~A.1 of Part I.  By construction, $x_k$
  is non-multiplicative for $\bar h$ and for each $h\in\H$ a polynomial
  $P_h\in\kk[x_1,\dots,x_{\cls{h}}]$ exists such that $x_k\bar
  h=\sum_{h\in\H}P_hh$.  No polynomial $h$ with $\cls{h}>d$ lies in
  $\lspan{x_1,\dots,x_d}$ (obviously
  $\lt{\prec}{h}\notin\lspan{x_1,\dots,x_d}$).  As the leading terms cannot
  cancel in the sum, $P_h\in\lspan{x_1,\dots,x_d}$ for all
  $h\in\H\setminus\H_d$.  Thus $x_k\bar h=\sum_{h\in\H_d}c_hh+g$ with
  $c_h\in\kk$ and $g\in\lspan{x_1,\dots,x_d}\I$.  As $\I$ is a homogeneous
  ideal and as the degree of $\bar h$ is maximal in $\H_d$, all constants
  $c_h$ must vanish.
  
  It is not possible that $\bar h\in\lspan{x_1,\dots,x_d}\I$, as otherwise
  $\bar h$ would be involutively head reducible by some other element of $\H$.
  Hence we have shown that any variable $x_k$ with $k>d$ is a zero divisor in
  $\I/\lspan{x_1,\dots,x_d}\I$ and the $\I$-regular sequence $(x_1,\dots,x_d)$
  cannot be extended by any $x_k$ with $k>d$.  Obviously, the same argument
  applies to any linear combination of such variables $x_k$.
  
  Finally, assume that $y_1,\dots,y_{d+1}\in\P_1$ form an $\I$-regular
  sequence of length $d+1$.  We extend them to a basis $\{y_1,\dots,y_n\}$ of
  $\P_1$ and perform a coordinate transformation $\xv\mapsto\yv$.  Our basis
  $\H$ transforms into the set $\H_{\yv}$ and after an involutive head
  autoreduction we obtain the set $\H^\triangle_{\yv}$.  In general, the
  coordinates $\yv$ are not $\delta$-regular for $\H^\triangle_{\yv}$.  But
  there exist coordinates $\tilde{\yv}$ of the form $\tilde
  y_k=y_k+\sum_{i=1}^{k-1}a_{ki}y_i$ with $a_{ki}\in\kk$ such that if we
  transform $\H$ to them and perform afterwards an involutive head
  autoreduction, they are $\delta$-regular for the obtained set
  $\tilde\H^\triangle_{\yv}$.
  
  This implies that $\tilde\H^\triangle_{\yv}$ is a Pommaret basis of the
  ideal it generates.\footnote{The involutive sizes of $\H$ and
    $\tilde\H^\triangle_{\yv}$ are equal.  Thus local involution of $\H$
    implies local involution of $\tilde\H^\triangle_{\yv}$ and for the
    Pommaret division as a continuous division local involution is equivalent
    to involution by Corollary 7.3 of Part I.} Thus $\min_{\tilde
    h\in\tilde\H^\triangle_{\yv}}\cls{\tilde h}=d$ and by the same argument as
  above $\tilde y_{d+1}$ is a zero divisor in $\I/\lspan{\tilde
    y_1,\dots,\tilde y_d}\I$.  Because of the special form of the
  transformation from $\yv\mapsto\tilde{\yv}$, we have $\lspan{\tilde
    y_1,\dots,\tilde y_d}=\lspan{y_1,\dots,y_d}$ and $y_{d+1}$ must be a zero
  divisor in $\I/\lspan{y_1,\dots,y_d}\I$.  But this contradicts the
  assumption that $(y_1,\dots,y_{d+1})$ is an $\I$-regular sequence and
  $\depth{\I}=d$.\qed
\end{proof}

\begin{remark}\label{rem:depth}
  One may wonder to what extent this result really requires the Pommaret
  division.  Given an arbitrary involutive basis $\H$ of $\I$, we may
  introduce the set $X_{\I}=\bigcap_{h\in\H}\mult{X}{L,\H,\prec}{h}$;
  obviously, for a Pommaret basis $X_{\I}=\{x_1,\dots,x_d\}$ with
  $d=\min_{h\in\H}{\cls{h}}$.  Again it is trivial to conclude from the
  induced combinatorial decomposition that any sequence formed by elements of
  $X_{\I}$ is $\I$-regular.  But in general we cannot claim that these are
  \emph{maximal} $\I$-regular sequences and there does not seem to exist an
  obvious method to extend them.  Thus only a lower bound for the depth is
  obtained this way.
  
  As a simple example we consider the ideal $\I$ generated by $f_1=z^2-xy$,
  $f_2=yz-wx$ and $f_3=y^2-wz$.  If we set $x_1=w$, $x_2=x$, $x_3=y$ and
  $x_4=z$, then it is straightforward to check that the set
  $\F=\{f_1,f_2,f_3\}$ is a Pommaret basis of $\I$ with respect to the degree
  reverse lexicographic order.  By Proposition~\ref{prop:Iregular}, $(w,x,y)$
  is a maximal $\I$-regular sequence and $\depth{\I}=3$.
  
  If we set $x_1=w$, $x_2=z$, $x_3=y$ and $x_4=x$, then no finite Pommaret
  basis exists; these coordinates are not $\delta$-regular.  In order to
  obtain a Janet basis $\F_J$ of $\I$ (for the degree reverse lexicographic
  order with respect to the new ordering of the variables), we must enlarge
  $\F$ by $f_4=z^3-wx^2$ and $f_5=xz^3-wx^3$.  We find now $X_{\I}=\{w,z\}$, as
  \begin{equation}
    \begin{gathered}
      \mult{X}{J,\F_J,\drl}{f_1}=\{w,z,y,x\}\;,\quad
      \mult{X}{J,\F_J,\drl}{f_3}=\{w,z,y\}\;,\\
      \mult{X}{J,\F_J,\drl}{f_i}=\{w,z\}\quad \mbox{for}\ i=2,4,5\;,
    \end{gathered}
  \end{equation}
  One easily verifies that both $x$ and $y$ are not zero divisors in
  $\I/\lspan{w,z}\I$, so that $X_{\I}$ can be extended to a maximal
  $\I$-regular sequence.  However, the Janet basis gives no indications, how
  this could be done.  One could also conjecture that the minimal number of
  multiplicative variables for a generator gives the depth.  But clearly this
  is also not true for the above Janet basis.  Thus no obvious way seems to
  exist to deduce $\depth{\I}$ from it.\bull
\end{remark}

\section{Combinatorial Decompositions}\label{sec:decomp}

In the proof of Proposition~\ref{prop:Iregular} we could already see the power
of the direct sum decompositions induced by (strong) involutive bases.  In
this section we want to study this aspect in more details.  All results apply
to arbitrary finitely generated polynomial modules.  But for notational
simplicity, we restrict to graded $\kk$-algebras $\A=\P/\I$ with a homogeneous
ideal $\I\subseteq\P$.  If we speak of a basis of the ideal $\I$, we always
assume that it is homogeneous, too.

The main motivation of Buchberger for the introduction of Gr\"obner bases was
to be able to compute effectively in such factor spaces.  Indeed given a
Gr\"obner basis $\G$ of the ideal $\I$, the normal form with respect to $\G$
distinguishes a unique representative in each equivalence class.  Our goal in
this section is to show that Pommaret bases contain in addition much
structural information about the algebra~$\A$.  More precisely, we want to
compute fundamental invariants like the Hilbert polynomial (which immediately
yields the Krull dimension and the multiplicity), the depth or the
Castelnuovo-Mumford regularity (see Section~\ref{sec:cmr}).  Our basic tools
are combinatorial decompositions of the algebra $\A$ into direct sums of
polynomial rings with a restricted number of variables.

\begin{definition}\label{def:sd}
  A\/ \emph{Stanley decomposition} of the graded\/ $\kk$-algebra\/ $\A=\P/\I$
  is an isomorphism of graded\/ $\kk$-linear spaces
  \begin{equation}\label{eq:sd}
    \A\cong\bigoplus_{t\in\T}\kk[X_t]\cdot t
  \end{equation}
  with a finite set\/ $\T\subset\TT$ and sets\/
  $X_t\subseteq\{x_1,\dots,x_n\}$.
\end{definition}

The elements of the set $X_t$ are again called the \emph{multiplicative
  variables} of the generator $t$.  As a first, trivial application of such
decompositions we determine the Hilbert series and the (Krull) dimension.

\begin{proposition}[\cite{rps:hilbert}]\label{prop:hilbert}
  Let the graded algebra\/ $\A$ possess the Stanley decomposition\/
  (\ref{eq:sd}).  Then its Hilbert series is
  \begin{equation}\label{eq:hs}
    \H_\A(\lambda)=\sum_{t\in\T}
        \frac{\lambda^{q_t}}{(1-\lambda)^{k_t}}
  \end{equation}
  where\/ $q_t=\deg{t}$ and\/ $k_t=|X_t|$.  Thus the (Krull) dimension of\/
  $\A$ is given by\/ $D=\max_{t\in\T}{k_t}$ and the multiplicity (or
  degree) by the number of terms\/ $t\in\T$ with\/ $k_t=D$.
\end{proposition}

Vasconcelos \cite[p.~23]{wvv:comp} calls Stanley decompositions \emph{``an
  approach that is not greatly useful computationally but it is often nice
  theoretically''}.  One reason for this assessment is surely that the
classical algorithm for their construction works only for monomial ideals and
uses a recursion over the variables $x_1,\dots,x_n$.  Thus for larger $n$ it
becomes quite inefficient.  For a general ideal $\I$ one must first compute a
Gr\"obner basis of $\I$ for some term order $\prec$ and then, exploiting the
vector space isomorphism $\P/\I\cong\P/\lt{\prec}{\I}$, one determines a
Stanley decomposition.  Its existence is guaranteed by the following result.

\begin{proposition}\label{prop:paradecomp}
  Let\/ $\I\subseteq\Nn$ be a monoid ideal and\/ $\bar\I=\Nn\setminus\I$ its
  complementary set.  There exists a finite set\/ $\bar\N\subset\bar\I$ and
  for each multi index\/ $\nu\in\bar\N$ a set of indices\/
  $N_\nu\subseteq\{1,\dots,n\}$ such that\/\footnote{Recall from Part I the
  notation $\NN^n_N=\bigl\{\nu\in\Nn\mid\forall j\notin N:\nu_j=0\bigr\}$.}
  \begin{equation}\label{eq:paradecomp}
    \bar\I=\bigcup_{\nu\in\bar\N}\bigl(\nu+\NN^n_{N_\nu}\bigr)
  \end{equation}
  and\/ $(\nu+\NN^n_{N_\nu})\cap(\mu+\NN^n_{N_\mu})=\emptyset$ for all\/
  $\mu,\nu\in\bar\N$.
\end{proposition}

A proof of this proposition may be found in the textbook
\cite[pp.~417--418]{clo:iva} (there it is not shown that one can always
construct a disjoint decomposition, but this extension is trivial).  This
proof is not completely constructive, as a certain degree $q_0$ is defined by
a Noetherian argument.  But it is not difficult to see that we may take
$q_0=\max_{\nu\in\N}{\nu_n}$ where $\N$ is the minimal basis of the monoid
ideal $\I$.  Now one can straightforwardly transform the proof into a
recursive algorithm for the construction of Stanley decompositions.  In fact,
one obtains then the algorithm proposed by Sturmfels and White \cite{sw:comb}.

One must stress that the decomposition (\ref{eq:paradecomp}) is not unique and
different decompositions may use sets $\bar\N$ of different sizes.  Given an
involutive basis of the monoid ideal $\I$, it is trivial to determine a
disjoint decomposition of $\I$ itself.  However, there does not seem to exist
an obvious way to obtain a complementary decomposition (\ref{eq:paradecomp}).
The situation is different for Janet bases where already Janet
\cite[\textsection15]{ja:lec} presented a solution of this problem.  Again it
is easy to see that it can straightforwardly be extended to an algorithm.  We
recall his proof in our notations, as we will need it later on.

\begin{proposition}\label{prop:jandecomp}
  Let\/ $\N_J$ be a Janet basis of the monoid ideal\/ $\I\subseteq\Nn$.  Then
  the set\/ $\bar\N\subset\Nn$ in the decomposition\/ (\ref{eq:paradecomp})
  may be chosen such that for all\/ $\nu\in\bar\N$ the equality\/
  $N_\nu=\mult{N}{J,\N_J\cup\{\nu\}}{\nu}$ holds.
\end{proposition}

\begin{proof}\footnote{For an alternative proof see \cite{pr:janet}.}
  We explicitly construct the set $\bar\N$.  Let
  $q_n=\max_{\mu\in\N_J}{\mu_n}$.  Then we put into $\bar\N$ all multi indices
  $\nu=[0,\dots,0,q]$ such that $0\leq q<q_n$ and such that $\N_J$ does not
  contain any multi index $\mu$ with $\mu_n=q$.  We set
  $N_\nu=\{1,\dots,n-1\}$ which obviously is
  $\mult{N}{J,\N_J\cup\{\nu\}}{\nu}$ according to the definition of the Janet
  division.
  
  Let $(d_k,\dots,d_n)\subseteq\N_J$ be one of the subsets considered in the
  assignment of the Janet multiplicative variables to the elements of $\N_J$
  (see Example 2.2 in Part~I) and set
  $q_{k-1}=\max_{\mu\in(d_k,\dots,d_n)}{\mu_{k-1}}$.  We enlarge the set
  $\bar\N$ by all those multi indices $\nu=[0,\dots,0,q,d_k,\dots,d_n]$ where
  $0\leq q<q_{k-1}$ and where $(d_k,\dots,d_n)$ does not contain any multi
  index $\mu$ with $\mu_{k-1}=q$.  The corresponding sets $N_\nu$ of
  multiplicative indices contain the indices $1,\dots,k-1$ and all those
  indices greater than $k$ that are multiplicative for the elements of
  $(d_k,\dots,d_n)$.  Again it is easy to verify that this entails
  $N_\nu=\mult{N}{J,\N_J\cup\{\nu\}}{\nu}$.
  
  We claim that the thus constructed set $\bar\N$ (together with the sets
  $N_\nu$) defines a disjoint decomposition (\ref{eq:paradecomp}) of the
  complementary set $\bar\I$.  We proceed by an induction on $n$.  The case
  $n=1$ is trivial.  Suppose that our assertion is true for $n-1$.  Let
  $\rho\in\Nn$ be an arbitrary multi index and let $q_1<q_2<\cdots<q_r$
  represent all values that occur for the entry $\mu_n$ in the multi indices
  $\mu\in\N_J$.  We distinguish three cases depending on $\rho_n$: (i)
  $\rho_n<q_r$ but $\rho_n\notin\{q_1,\dots,q_{r-1}\}$,
  (ii)~$\rho_n\in\{q_1,\dots,q_{r-1}\}$, and (iii) $\rho_n\geq q_r$.
  
  In the first case, the set $\bar\N$ contains the multi index
  $\nu=[0,\dots,0,\rho_n]$ and, as $N_\nu=\{1,\dots,n-1\}$, obviously
  $\rho\in\nu+\NN^n_{N_\nu}$.  It is easy to see that $\rho$ cannot be an
  element of $\I$ or of the cone $\bar\nu+\NN^n_{N_{\bar\nu}}$ of any other
  multi index $\bar\nu\in\bar\N$.
  
  In the second case, we ``slice'' the two sets $\N_J$ and $\bar\N$ at degree
  $\rho_n$.  For this purpose we must introduce some notations.  Primed
  letters always refer to multi indices $\nu'\in\NN_0^{n-1}$; unprimed ones to
  multi indices $\nu\in\Nn$.  For a given $\nu\in\Nn$, we define
  $\nu'=[\nu_1,\dots,\nu_{n-1}]\in\NN_0^{n-1}$, i.\,e.\ we simply drop the
  last entry.  Finally, we write $[\nu',q]$ for the multi index
  $[\nu_1,\dots,\nu_{n-1},q]\in\Nn$.  If we are given some monoid ideal
  $\I\subseteq\Nn$, we define for every integer $q\in\NN_0$ the monoid ideal
  $\I_q'=\bigl\{\nu'\in\NN_0^{n-1}\mid[\nu',q]\in\I\bigr\}\subseteq\NN_0^{n-1}$.
  
  The ``slicing'' yields the two sets:
  $\N_J'=\bigl\{\mu'\in\NN_0^{n-1}\mid\mu\in\N_J,\ \mu_n=\rho_n\bigr\}$ and
  $\bar\N'=\bigl\{\nu'\in\NN_0^{n-1}\mid\nu\in\bar\N,\ \nu_n=\rho_n\bigr\}$,
  respectively.  If we compute the sets $\mult{N}{J,\N_J'}{\mu'}$ of
  multiplicative indices for the elements $\mu'\in\N'_J$, it is
  straightforward to verify that they are just $\mult{N}{J,\N_J}{\mu}$, as
  $\mu_n$ is not maximal and $\mu\in(\rho_n)$, if and only if $\mu'\in\N'_J$.
  Furthermore, $\N_J'$ is a Janet basis of the monoid ideal $\I_{\rho_n}'$.
  If we apply the procedure above to this ideal in $\NN_0^{n-1}$, we obtain
  $\bar\N'$ as complementary basis and the sets of multiplicative variables
  remain unchanged.
  
  By our inductive hypothesis, $\bar\N'$ defines a disjoint decomposition of
  the sought form of the set $\bar\I_{\rho_n}'$.  If $\rho\in\bar\I$, then
  $\rho'\in\bar\I_{\rho_n}'$ and a multi index $\nu'\in\bar\N'$ exists such
  that $\rho'\in\nu'+\NN^{n-1}_{N_{\nu'}}$.  By construction,
  $\nu=[\nu',\rho_n]\in\bar\N$ and $\rho\in\nu+\NN^n_{N_\nu}$.  If
  $\rho\in\I$, then $\rho'\in\I_{\rho_n}'$, and it is not possible to find a
  multi index $\nu'\in\bar\N'$ such that $\rho'\in\nu'+\NN^{n-1}_{N_{\nu'}}$.
  Hence $\rho\notin\nu+\NN^n_{N_\nu}$ for any $\nu\in\bar\N$.
  
  The third case is very similar to the second one.  For the definition of
  $\N_J'$ and $\bar\N'$ we consider only multi indices where the value of the
  last entry is $q_r$.  Note that for all multi indices $\nu\in\bar\N$ that
  contribute to $\bar\N'$ the index $n$ is multiplicative.  If
  $\rho\in\bar\I$, then $\rho'\in\bar\I_{\rho_n}'$ and a multi index
  $\nu'\in\bar\N'$ exists such that $\rho'\in\nu'+\NN^{n-1}_{N_{\nu'}}$ and we
  conclude as above that $\rho\in\nu+\NN^n_{N_\nu}$ for
  $\nu=[\nu',q_r]\in\bar\N$.  Again it is obvious that for a multi index
  $\rho\in\I$, it is not possible to find such a $\nu\in\bar\N$.\qed
\end{proof}

We have formulated this proposition only for the case that $\N_J$ is a Janet
basis, but the proof yields an even stronger result.  Let $\N_J$ be an
\emph{arbitrary} subset of $\Nn$ and construct the corresponding set $\bar\N$.
Then we may substitute everywhere in the proof $\I$ by the involutive span
$\ispan{\N_J}{J}$ and still obtain that any multi index $\rho\in\Nn$ either
lies in $\ispan{\N_J}{J}$ or in exactly one of the cones defined by $\bar\N$
and the sets $N_\nu$.

\begin{remark}
  Janet did not formulate his algorithm in this algebraic language.  He
  considered the problem of determining a formally well-posed initial value
  problem for an overdetermined system of partial differential equations
  \cite{wms:habil}.  Identifying this system with our ideal $\I$, his problem
  is equivalent to computing a Stanley decomposition of $\P/\I$.  An
  algorithmic approach to formally well-posed initial value problem was also
  presented by Reid \cite{gr:red}.\bull
\end{remark}

According to Corollary~\ref{cor:janpom}, we may apply the construction of
Proposition~\ref{prop:jandecomp} to Pommaret bases, too.  But as the Pommaret
division has such a simple global definition, it is almost trivial to provide
an alternative decomposition depending only on the degree $q$ of a Pommaret
basis of the ideal $\I$ (we will see later in Sect.~\ref{sec:cmr} that this
degree is in fact an important invariant of $\I$).  In general, this
decomposition is larger than the one obtained with the Janet approach but it
has some advantages in theoretical applications.

\begin{proposition}\label{prop:pomdecomp}
  The monoid ideal\/ $\I\subseteq\Nn$ has a Pommaret basis of degree\/ $q$, if
  and only if the sets\/ $\bar\N_0=\{\nu\in\bar\I\mid|\nu|<q\}$ and\/
  $\bar\N_1=\{\nu\in\bar\I\mid|\nu|=q\}$ yield the disjoint decomposition
  \begin{equation}\label{eq:pomdecomp}
    \bar\I=\bar\N_0\cup\bigcup_{\nu\in\bar\N_1}\cone{P}{\nu}\;.
  \end{equation}
\end{proposition}

\begin{proof}
  The definition of the Pommaret division implies the identity
  \begin{equation}
    \left(\Nn\right)_{\geq q}=\bigcup_{|\nu|=q}\cone{P}{\nu}
  \end{equation}
  from which one direction of the proposition follows trivially.  Here
  $\left(\Nn\right)_{\geq q}$ denotes the set of all multi indices of length
  greater than or equal to $q$.  By the definition of an involutive division,
  the union on the right hand side is disjoint.
  
  For the converse, we claim that the set $\H=\{\mu\in\I_q\}$ is a Pommaret
  basis of monoid ideal $\I_{\geq q}$; this immediately implies our assertion
  by Lemma \ref{lem:truncbas}.  Assume that $\mu\in\H$ with $\cls\mu=k$ and
  let $k<j\leq n$ be a non-multiplicative index for it.  We must show that
  $\mu+1_j\in\ispan{\H}{P}$.  But this is trivial: we have
  $\mu+1_j\in\cone{P}{\mu-1_k+1_j}$ and $\mu-1_k+1_j\in\I_q$, as otherwise we
  encounter the contradiction $\mu+1_j\notin\I$ by (\ref{eq:pomdecomp}).\qed
\end{proof}

\begin{example}\label{ex:pomdecomp}
  The decomposition (\ref{eq:pomdecomp}) is usually redundant.  Considering
  for $\bar\N_1$ only multi indices of length $q$ makes the formulation much
  easier but it is not optimal.  Consider the trivial example
  $\N_P=\bigl\{[0,1]\bigr\}$.  According to Proposition~\ref{prop:pomdecomp}
  we should set $\bar\N_0=\bigl\{[0,0]\bigr\}$ and
  $\bar\N_1=\bigl\{[1,0]\bigr\}$.  But obviously $\bar\I=[0,0]+\NN^2_{\{1\}}$.
  Applying the construction used in the proof of
  Proposition~\ref{prop:jandecomp} yields directly this more compact
  form.\bull
\end{example}

If $\J\subseteq\P$ is a polynomial ideal possessing a Pommaret basis for some
term order $\prec$, then applying Proposition~\ref{prop:pomdecomp} to
$\I=\le{\prec}{\J}$ yields a Stanley decomposition of a special type: all sets
$X_t$ are of the form $X_t=\{x_1,x_2,\dots,x_{\cls{t}}\}$ where the number
$\cls{t}$ is called the \emph{class}\footnote{Some authors prefer the term
  \emph{level}.}  of the generator $t$.  One speaks then of a \emph{Rees
  decompositions} of $\A=\P/\J$ \cite{rees:polymod}.  It is no coincidence
that we use here the same terminology as in the definition of the Pommaret
division: if $t=x^\mu$ with $\mu\in\bar\N_1$, then indeed its class is
$\cls{\mu}$.  Reducing in a straightforward manner the redundancy in the
decomposition (\ref{eq:pomdecomp}) leads to the following result.

\begin{corollary}\label{cor:stacon}
  Let\/ $\I\subseteq\P$ be a polynomial ideal which has for some term order\/
  $\prec$ a Pommaret basis\/ $\H$ such that
  $\min_{h\in\H}{\cls{\le{\prec}{h}}}=d$.  Then\/ $\P/\I$ possesses a Rees
  decomposition where the minimal class of a generator is\/ $d-1$.
\end{corollary}

\begin{proof}
  Obviously, it suffices to consider the monomial case and formulate the proof
  therefore in the multi index language of Proposition~\ref{prop:pomdecomp}.
  Furthermore, for $d=1$ there is nothing to be shown so that we assume from
  now on $d>1$.  Our starting point is the decomposition (\ref{eq:pomdecomp}).
  For each $\nu\in\bar\N_1$ with $\cls{\nu}=k<d$ we introduce the multi
  index  $\tilde\nu=\nu-(\nu_k)_k$, i.\,e.\ $\tilde\nu$ arises from
  $\nu$ by setting the $k$th entry to zero.  Obviously, the $k$-dimensional
  cone $C_\nu=\tilde\nu+\NN^n_{\{1,\dots,k\}}$ is still completely contained
  in the complement $\bar\I$ and $\cone{P}{\nu}\subset C_\nu$.

  If we replace in (\ref{eq:pomdecomp}) for any such $\nu$ the cone
  $\cone{P}{\nu}$ by $C_\nu$, then we still have a decomposition of $\bar\I$,
  but no longer a disjoint one.  We now show first that in the thus obtained
  decomposition all cones $C$ with $0<\dim{C}<d-1$ can be dropped without
  loss.  Indeed, for $k<d-1$ we consider the multi index
  $\mu=\tilde\nu+(\nu_k)_{k+1}$. Obviously, $|\mu|=q$ and $\cls{\mu}=k+1$;
  hence under the made assumptions $\mu\in\bar\N_1$.  Furthermore,
  $\tilde\mu=\mu-(\mu_{k+1})_{k+1}$ is a divisor of $\tilde\nu$ (the two multi
  indices can differ at most in their $(k+1)$st entries and
  $\tilde\mu_{k+1}=0$) and thus the inclusion $C_\nu\subset
  C_\mu=\tilde\mu+\NN^n_{\{1,\dots,k+1\}}$ holds. 

  The remaining cones with $\dim{C}\geq d-1$ are all disjoint. This is
  trivially true for all cones with $\dim{C}\geq d$, as these have not been
  changed.  For the other ones, we note that if $\mu$ and $\nu$ are two multi
  indices with $\cls{\mu}=\cls{\nu}=d-1$ and $|\mu|=|\nu|=q$, then they must
  differ at some position $\ell$ with $\ell\geq d$.  But this implies that the
  cones $C_\mu$ and $C_\nu$ are disjoint.
  
  Thus there only remains to study the zero-dimensional cones consisting of
  the multi indices $\nu\in\bar\N_0$.  If we set $\ell=q-|\nu|$, then
  $\mu=\nu+\ell_1\in\bar\N_1$, since we assumed $d>1$, and trivially $\nu\in
  C_\mu=(\mu-(\mu_1)_1)+\NN^n_{\{1\}}$.  By our considerations above the cone
  $C_\mu$ and thus $\nu$ is contained in some $(d-1)$-dimensional cone.
  Therefore we may also drop all zero-dimensional cones and obtain a Rees
  decomposition where all cones are at least $(d-1)$-dimensional. \qed
\end{proof}

Sturmfels et al.\ \cite{stv:degree} introduced the notion of \emph{standard
  pairs} also leading to a kind of combinatorial decomposition, however not a
disjoint one.  They consider pairs $(\nu,N_\nu)$ where $\nu\in\Nn$ is a multi
index and $N_\nu\subseteq\{1,\dots,n\}$ a set of associated indices.  Such a
pair is called \emph{admissible}, if $\supp{\nu}\cap N_\nu=\emptyset$, i.\,e.\ 
$\nu_i=0$ for all $i\in N_v$.  On the set of admissible pairs one defines a
partial order: $(\nu,N_\nu)\leq(\mu,N_\mu)$, if and only if the restricted
cone $\mu+\NN^n_{N_\mu}$ is completely contained in $\nu+\NN^n_{N_\nu}$.
Obviously, this is equivalent to $\nu\mid\mu$ and any index $i$ such that
either $\mu_i>\nu_i$ or $i\in N_\mu$ is contained in $N_\nu$.

\begin{definition}\label{def:stanpair}
  Let\/ $\I\subseteq\Nn$ be an arbitrary monoid ideal.  An admissible pair\/
  $(\nu,N_\nu)$ is called\/ \emph{standard} for\/ $\I$, if\/
  $\nu+\NN^n_{N_\nu}\cap\I=\emptyset$ and\/ $(\nu,N_\nu)$ is minimal with
  respect to\/ $<$ among all admissible pairs with this property.
\end{definition}

Any monoid ideal $\I\subset\Nn$ leads thus automatically to a uniquely
determined set of standard pairs.  These define both a decomposition of the
complementary set $\overline{\I}$ into cones (though these will overlap in
general) and a decomposition of the ideal $\I$ itself as an intersection of
irreducible monomial ideals.

\begin{proposition}[\cite{stv:degree}]\label{prop:stanpair}
  Let\/ $\I\subseteq\Nn$ be an arbitrary monoid ideal and denote the set of
  all associated standard pairs by\/
  $\S_{\I}=\bigl\{(\nu,N_\nu)\mid(\nu,N_\nu) \mbox{\ \upshape{standard for}\ 
  }\I\bigr\}$.  Then the complementary set\/ $\overline{\I}$ of\/ $\I$ can be
  written in the form
  \begin{equation}\label{eq:stanpair1}
    \overline{\I}=\bigcup_{(\nu,N_\nu)\in\S_{\I}}\nu+\NN^n_{N_\nu}
  \end{equation}
  and the ideal\/ $\I$ itself can be decomposed as
  \begin{equation}\label{eq:stanpair2}
    \I=\bigcap_{(\nu,N_\nu)\in\S_{\I}}
        \lspan{(\nu_i+1)_i\mid i\notin N_\nu}\;.
  \end{equation}
\end{proposition}

According to Sturmfels et al. \cite[Lemma 3.3]{stv:degree} the number of
standard pairs of a monomial ideal $\I$ equals the \emph{arithmetic degree} of
$\I$, a refinement of the classical concept of the degree of an ideal
introduced by Bayer and Mumford \cite{bm:cac}.  We further note that the
ideals on the right hand side of (\ref{eq:stanpair2}) are trivially
irreducible, so that (\ref{eq:stanpair2}) indeed represents an irreducible
decomposition of $\I$.

In general, this decomposition is highly redundant.  Let
$N\subseteq\{1,\dots,n\}$ be an arbitrary subset and consider all standard
pairs $(\nu,N_\nu)$ with $N_\nu=N$.  Obviously, among these only the ones with
multi indices $\nu$ which are maximal with respect to divisibility are
relevant for the decomposition (\ref{eq:stanpair2}) and in fact restricting to
the corresponding ideals yields the irredundant irreducible decomposition of
$\I$ (which is unique according to \cite[Thm.~5.27]{ms:cca}).  Their
intersection defines a possible choice for the primary component for the prime
ideal $\pf_N=\lspan{x_i\mid i\notin N}$, so that we can also extract an
irredundant primary decomposition from the standard pairs.  As a trivial
corollary to these considerations the standard pairs immediately yield the set
$\ass{(\P/\I)}$ of associated prime ideals, as it consists of all prime ideals
$\pf_N$ such that a standard pair $(\nu,N)$ exists.

Ho\c{s}ten and Smith \cite{hs:monom} discuss two algorithms for the direct
construction of the set $\S_{\I}$ of all standard pairs given the minimal
basis of $\I$.  Alternatively, $\S_{\I}$ can easily be extracted from any
complementary decomposition, as we show now.  Thus once a Janet basis of $\I$
is known, we may use Janet's algorithm for the construction of a complementary
decomposition and then obtain the standard pairs.

Let the finite set $\T_{\I}=\bigl\{(\nu,N_\nu)\mid\nu\in\Nn,
N_\nu\subset\{1,\dots,n\}\bigr\}$ define a complementary decomposition of
$\I$.  If the pair $(\nu,N_\nu)\in\T_{\I}$ is not admissible, then we
substitute it by the pair $(\bar\nu,N_\nu)$ where $\bar\nu_i=0$ for all $i\in
N_\nu$ and $\bar\nu_i=\nu_i$ else.  Obviously, this operation produces an
admissible pair and the thus obtained set $\overline{\S}_{\I}$ still defines a
(generally no longer disjoint) decomposition of the complementary set
$\overline{\I}$.  Finally, we eliminate all pairs in $\overline{\S}_{\I}$
which are not minimal with respect to the partial order $\leq$ and obtain a
set $\S_{\I}$.

\begin{proposition}\label{prop:standard}
  Let\/ $\T_{\I}$ be a finite complementary decomposition of the monoid
  ideal\/ $\I\subseteq\Nn$.  The thus constructed set\/ $\S_{\I}$ consists of
  all standard pairs of\/ $\I$.
\end{proposition}

\begin{proof}
  It is trivial to see that the set $\overline{\S}_{\I}$ contains only
  admissible pairs and that $\nu+\NN^n_{N_\nu}\subseteq\overline{\I}$ for any
  pair $(\nu,N_\nu)\in\overline{\S}_{\I}$.  Thus there only remains to show
  that all standard pairs are contained in $\overline{\S}_{\I}$.
  
  Let $(\mu,N_\mu)$ be an admissible pair such that
  $\mu+\NN^n_{N_\mu}\subseteq\overline{\I}$.  Since the union of the cones
  $\nu+\NN^n_{N_\nu}$ with $(\nu,N_\nu)\in\overline{\S}_{\I}$ still covers
  $\overline{\I}$, the finiteness of $\overline{\S}_{\I}$ implies the
  existence of a multi index $\overline{\mu}\in\mu+\NN^n_{N_\mu}$ and a pair
  $(\nu,N_\nu)\in\overline{\S}_{\I}$ such that
  $\overline{\mu}+\NN^n_{N_\mu}\subseteq\nu+\NN^n_{N_\nu}$ (obviously, it is
  not possible to cover $\mu+\NN^n_{N_\mu}$ with a finite number of
  lower-dimensional cones).  As both $(\mu,N_\mu)$ and $(\nu,N_\nu)$ are
  admissible pairs, this entails that in fact $(\nu,N_\nu)\leq(\mu,N_\mu)$.
  Hence either $(\mu,N_\mu)\in\overline{\S}_{\I}$ or it is not a standard
  pair.\qed
\end{proof}

\begin{remark}\label{rem:stanpom}
  If we use the decomposition (\ref{eq:pomdecomp}) derived from a Pommaret
  basis of degree $q$, then the determination of the set $\overline{\S}_{\I}$
  is completely trivial.  For all pairs $(\nu,N_\nu)\in\T_{\I}$ with $|\nu|<q$
  we have $N_\nu=\emptyset$ and hence they are trivially admissible.  For all
  other pairs we find that $\supp{\nu}\cap N_\nu=\{\cls{\nu}\}$. Thus none of
  them is admissible, but they become admissible by simply setting the first
  non-vanishing entry of $\nu$ to zero.\bull
\end{remark}

\begin{example}\label{ex:stanpair}
  Consider the ideal $\I=\lspan{z^3,yz^2-xz^2,y^2-xy}\subset\kk[x,y,z]$.  Both
  a Janet and Pommaret basis of $\I$ for the degree reverse lexicographic
  order is given by the set $\H=\{z^3,yz^2-xz^2,y^2z-xyz,y^2-xy\}$.  Following
  the construction in the proof of Proposition~\ref{prop:jandecomp}, we obtain
  the set $\T=\{1,y,z,yz,z^2\}$ and the complementary decomposition
  \begin{equation}\label{eq:hircomp}
    \P/\I\cong
        \kk[x]\oplus\kk[x]\cdot y\oplus\kk[x]\cdot z\oplus
        \kk[x]\cdot yz\oplus\kk[x]\cdot z^2\;.
  \end{equation}
  The term $z^2$ comes from the set $(2)=\bigl\{[0,1,2]\bigr\}$.  The terms
  $yz$ and $z$ arise from $(1)=\bigl\{[0,2,1]\bigr\}$ and, finally, $1$ and
  $y$ stem from $(0)=\bigl\{[0,2,0]\bigr\}$.  The last non-vanishing set
  $(3)=\bigl\{[0,0,3]\bigr\}$ does not contribute to the decomposition.
  
  It follows from (\ref{eq:hircomp}) that a complementary decomposition of the
  corresponding monoid ideal $\le{\prec}{\I}=\lspan{[0,0,3], [0,1,2],
    [0,2,0]}$ is given by
  \begin{equation}\label{eq:hirstanpair}
    \begin{split}
      \S_{\I}=\Bigl\{&\bigl([0,0,0],\{1\}\bigr), 
                      \bigl([0,1,0],\{1\}\bigr),
                      \bigl([0,0,1],\{1\}\bigr),\\
                     &\bigl([0,1,1],\{1\}\bigr),
                      \bigl([0,0,2],\{1\}\bigr)\Bigr\}
    \end{split}
  \end{equation}
  and one easily verifies that these are all standard pairs.
  
  The complementary decomposition constructed via Proposition
  \ref{prop:pomdecomp} is much larger.  Besides many multi indices without any
  multiplicative indices, we obtain the following six multi indices for which
  $1$ is the sole multiplicative index: $[3,0,0]$, $[2,1,0]$, $[2,0,1]$,
  $[1,1,1]$, $[1,2,0]$ and $[1,0,2]$.  After setting the first entry to zero,
  we find precisely the multi indices appearing in (\ref{eq:hirstanpair}) plus
  the multi index $[0,2,0]$.  As $([0,1,0],\{1\})<([0,2,0],\{1\})$, the latter
  pair is not minimal.  The same holds for all pairs corresponding to the
  multi indices without multiplicative indices and hence we also arrive at
  (\ref{eq:hirstanpair}).\bull
\end{example}

As an application of Rees decompositions, we show that given a Pommaret basis
of the ideal $\I$, we can easily read off the dimension and the depth of the
algebra $\A=\P/\I$.  In principle, the determination of the dimension is of
course already settled by Proposition~\ref{prop:hilbert} for arbitrary
involutive bases.   However, in the case of a Pommaret basis a further useful
characterisation of $\dim{\A}$ exists.

\begin{proposition}\label{prop:dim}
  Let\/ $\H$ be a homogeneous Pommaret basis of the homogeneous ideal\/
  $\I\subseteq\P$ with\/ $\deg{\H}=q$ for some term order\/ $\prec$.  Then the
  dimension\/ $D$ of the algebra\/ $\A=\P/\I$ is
  \begin{equation}\label{eq:dimA}
    D=\min{\bigl\{i\mid\lspan{\H,x_1,\dots,x_i}_q=\P_q\bigr\}}\,.
  \end{equation}
\end{proposition}

\begin{proof}
  The Hilbert polynomials of $\A$ and the truncation $\A_{\geq q}$ coincide.
  Thus it suffice to consider the latter algebra.  By
  Lemma~\ref{lem:truncbas}, a Pommaret basis of $\I_{\geq q}$ is given by the
  set $\H_q$ determined in (\ref{eq:truncbas}).  If $D$ is the smallest number
  such that $\lspan{\H_q,x_1,\dots,x_D}_q=\P_q$, then all multi indices $\nu$
  with $|\nu|=q$ and $\cls\nu>D$ lie in $\le{\prec}{\H_q}$ but a multi index
  $\mu$ exists such that $|\mu|=q$, $\cls\mu=D$ and
  $\mu\notin\le{\prec}{\H_q}$.  By Proposition~\ref{prop:pomdecomp}, this
  observation entails that $\mu$ is a generator of class $D$ of the
  complementary decomposition (\ref{eq:pomdecomp}) and that the decomposition
  does not contain a generator of higher class.  But this trivially implies
  that $\dim\A=D$.\qed
\end{proof}

In a terminology apparently introduced by Gr\"obner
\cite[Section~131]{wg:mag}, a subset $X_{\I}\subseteq\{x_1,\dots,x_n\}$ is
called \emph{independent modulo} $\I$, if $\I\cap\kk[X_{\I}]=\{0\}$.  If even
$\lt{\prec}{\I}\cap\kk[X_{\I}]=\{0\}$ for some term order $\prec$, then one
speaks of a \emph{strongly} independent set for $\prec$.  One can show that
the maximal size of either an independent or a strongly independent set
coincides with $\dim{\A}$.  This approach to determining the dimension of an
ideal has been taken up again by Kredel and Weispfenning \cite{kw:dim} using
Gr\"obner bases (see also \cite[Sects 6.3 \&\ 9.3]{bw:groe}).

\begin{corollary}\label{cor:indset}
  Under the assumptions of Proposition~\ref{prop:dim}, $\{x_1,\dots,x_D\}$ is
  the unique maximal strongly independent set modulo the ideal\/ $\I$.
\end{corollary}

\begin{proof}
  As above, it suffices to consider the truncation $\I_{\geq q}$.  Assume that
  a non-zero polynomial $f\in\I_{\geq q}$ exists with
  $\lt{\prec}{f}\in\kk[x_1,\dots,x_D]$.  As $\H_q$ is a Pommaret basis of the
  ideal $\I_{\geq q}$, it must contain a polynomial $h$ such that
  $\le{\prec}{h}\idiv{P}\le{\prec}{f}$; this implies
  $\lt{\prec}{h}\in\kk[x_d,\dots,x_D]$.  Repeated applications of
  Lemma~\ref{lem:stable} to $\le{\prec}{h}$ with $j=D$ yields that
  $\nu=q_D\in\le{\prec}{\H_q}$.  Further applications of
  Lemma~\ref{lem:stable} to $\nu$ with $k=D$ show that in fact any multi index
  of length $q$ and class $D$ is contained in $\le{\prec}{\H_q}$.  But thus
  already $\lspan{\H,x_1,\dots,x_{D-1}}_q=\P_q$ in contradiction to the
  characterisation of $D$ in Proposition~\ref{prop:dim}.
  
  Hence $\{x_1,\dots,x_D\}$ is a strongly independent set.  Any larger
  strongly independent set would have to contain at least one variable $x_{j}$
  with $j>D$.  However, according to the proof of Proposition~\ref{prop:dim}
  any multi index $\mu$ with $\cls{\mu}>D$ and thus in particular $\mu=j_{q}$
  lies in $\le{\prec}{\H_{q}}$.  Therefore no strongly independent set may
  contain such a variable $x_{j}$.\qed
\end{proof}

Applying standard arguments in homological algebra to the exact sequence
$0\rightarrow\I\rightarrow\P\rightarrow\P/\I\rightarrow0$, one can easily show
that $\depth{(\P/\I)}=\depth{\I}-1$.  Hence Proposition \ref{prop:Iregular}
immediately implies the following result (one can also easily prove it
directly along the lines of the proof of Proposition \ref{prop:Iregular}).

\begin{proposition}\label{prop:depth}
  Let\/ $\H$ be a homogeneous Pommaret basis of the homogeneous ideal\/
  $\I\subseteq\P$ for a class respecting term order and\/
  $d=\min_{h\in\H}\cls{h}$.  Then the depth of\/ $\A=\P/\I$ is\/
  $\depth{\A}=d-1$.
\end{proposition}

Since $\{x_1,\dots,x_{d-1}\}$ is trivially a strongly independent set modulo
$\I$, we obviously find that always $D\geq d-1$.  Thus as a trivial corollary
to Propositions \ref{prop:dim} and \ref{prop:depth}, we find the well-known
fact that for any graded algebra $\A=\P/\I$ the inequality
$\depth\A\leq\dim\A$ holds.  In the limit case $\depth\A=\dim\A$, the algebra
is by definition \emph{Cohen-Macaulay} and we obtain the following
characterisation of such algebras.

\begin{theorem}\label{thm:cohmac}
  Let\/ $\H$ be a Pommaret basis of degree\/ $q$ of the homogeneous ideal\/
  $\I\subseteq\P$ for a class respecting term order\/ $\prec$ and set\/
  $d=\min_{h\in\H}\cls{h}$.  The algebra\/ $\A=\P/\I$ is Cohen-Macaulay, if
  and only if\/ $\lspan{\H,x_1,\dots,x_{d-1}}_q=\P_q$.
\end{theorem}

An alternative characterisation, which is more useful for computations, is
based on the existence of a special kind of Rees decomposition; one sometimes
speaks of a \emph{\index*{Hironaka decomposition}}, a terminology introduced
in \cite[Sect.~2.3]{bs:inv}.

\begin{corollary}\label{cor:hironaka}
  $\A=\P/\I$ is a Cohen-Macaulay algebra, if and only if a Rees decomposition
  of\/ $\A$ exists where all generators have the same class.
\end{corollary}

\begin{proof}
  One direction is trivial.  If such a special decomposition exists with $d$
  the common class of all generators, then obviously both the dimension and
  the depth of $\A$ is $d$ and thus $\A$ is Cohen-Macaulay.
  
  For the converse, let us assume that $\A$ is a Cohen-Macaulay algebra and
  that $\dim\A=\depth\A=d$.  Let $\H$ be a Pommaret basis of $\I$ with respect
  to the degree reverse lexicographic order.  By Theorem~\ref{thm:expombas},
  such a basis always exists in $\delta$-regular variables $\xv$.  Proposition
  \ref{prop:depth} implies that $\min_{h\in\H}\cls{h}=d+1$.  We introduce the
  set $\bar\N=\bigl\{\nu\in\Nn\setminus\lspan{\le{\prec}{\H}}\mid\cls\nu>d\}$
  (recall that by convention we defined $\cls{[0,\dots,0]}=n$ so that
  $[0,\dots,0]\in\bar\N$ whenever $\I\neq\P$).  $\bar\N$ is finite, as all its
  elements satisfy $|\nu|<\deg{\H}$ by Theorem \ref{thm:cohmac}, and we claim
  that
  \begin{equation}\label{eq:hironaka}
    \A\cong\bigoplus_{\nu\in\bar\N}\kk[x_1,\dots,x_d]\cdot x^\nu\;.
  \end{equation}
  
  In fact, this is precisely the decomposition obtained by applying the
  algorithm implicitly contained in the proof of
  Proposition~\ref{prop:jandecomp}.  Consider any multi index $\nu\in\bar\N$;
  obviously, it is of the form $\nu=[0,\dots,0,\nu_{d+1},\dots,\nu_n]$ with
  $\sum_{i=d+1}^n\nu_i<q=\deg{\H}$.  If we set $q'=q-\sum_{i=d+2}^n\nu_i$, then
  by Theorem~\ref{thm:cohmac} the multi index
  $[0,\dots,0,q',\nu_{d+2},\dots,\nu_n]$ lies in the monoid ideal
  $\lspan{\le{\prec}{\H}}$.  But this implies the existence of a multi index
  $\nu'\in\le{\prec}{\H}$ with
  $\nu'=[0,\dots,0,\nu_{d+1}',\nu_{d+2},\dots,\nu_n]$ with
  $\nu_{d+1}<\nu_{d+1}'\leq q'$.  Hence the set $(\nu_{d+2},\dots,\nu_n)$ is
  considered in the assignment of multiplicative variables to the elements of
  $\H$ for the Janet division and it consists only of the multi index $\nu'$,
  as $\H$ is involutively head autoreduced (with respect to the Pommaret
  division).  But this implies that our algorithm chooses $\nu$ as an element
  of $\bar\N$ and assigns to it the multiplicative variables $x_1,\dots,x_d$.
  
  The algorithm cannot lead to a larger set $\bar\N$, as any further multi
  index would be of class less than or equal to $d$ and thus be contained in
  $\kk[x_1,\dots,x_d]\cdot 1$.  But we know that the sets are disjoint, so
  that this cannot happen and we get the decomposition
  (\ref{eq:hironaka}).\qed
\end{proof}

\begin{example}\label{ex:hironaka}
  Consider again the ideal $\I=\lspan{z^3,yz^2-xz^2,y^2-xy}\subset\kk[x,y,z]$
  of Example \ref{ex:stanpair}.  It follows from the Pommaret basis given
  there that both the depth and the dimension of $\P/\I$ is $1$.  Hence
  $\A=\P/\I$ is Cohen-Macaulay and indeed (\ref{eq:hircomp}) is a Hironaka
  decomposition.\bull
\end{example}

\section{Noether Normalisation and Primary Decomposition}
\label{sec:noether}

One immediately sees that any strongly independent set is automatically also
an independent set (the converse is generally not true).  Hence Corollary
\ref{cor:indset} entails that $\I\cap\kk[x_1,\dots,x_D]=\{0\}$ and now it is
easy to show that a Pommaret basis of the ideal $\I$ defines a \emph{Noether
  normalisation} \cite[Def.~3.4.2]{gp:singular} of the algebra $\A$ (or that
$x_1,\dots,x_D$ form a \emph{homogeneous system of parameters} for it).

\begin{corollary}\label{cor:noethernorm}
  Under the assumptions of Proposition~\ref{prop:dim}, the restriction of the
  canonical projection\/ $\pi:\P=\kk[x_1,\dots,x_n]\rightarrow\A=\P/\I$ to\/
  $\kk[x_1,\dots,x_D]$ is a Noether normalisation of\/ $\A$.  If\/ $\H$ is a
  Pommaret bases with respect to the (degree) lexicographic order, then it is
  even a general Noether normalisation.
\end{corollary}

\begin{proof}
  It follows immediately from Corollary~\ref{cor:indset} that the restriction
  of the projection $\pi$ to $\kk[x_1,\dots,x_D]$ is injective.  Thus there
  only remains to show that $\A$ is finitely generated as a
  $\kk[x_1,\dots,x_D]$-module.  Proposition~\ref{prop:pomdecomp} gives us the
  complementary decomposition (\ref{eq:pomdecomp}) for $\le{\prec}{\I}$ which
  is defined by a finite set $\N\subset\Nn$.  As for each generator in $\N$
  the associated multiplicative indices form a subset of $\{1,\dots,D\}$ and
  since the complement of $\lt{\prec}{\I}$ is a basis of $\A$ as $\kk$-linear
  space, the finite set $\{\pi(x^\nu)\mid\nu\in\N\}$ generates $\A$ as
  $\kk[x_1,\dots,x_D]$-module.
  
  It follows from the considerations in the proof of
  Proposition~\ref{prop:dim} that for each $D<k\leq n$ the basis $\H_q$
  contains a generator $h_k$ with leading term $x_k^q$.  If $\H$ is a Pommaret
  basis with respect to the (degree) lexicographic order, then $h_k$ must be
  of the form $h_k=x_k^q+\sum_{j=0}^{q-1}P_{k,j}x_k^j$ with polynomials
  $P_{k,j}\in\kk[x_1,\dots,x_{k-1}]$.  Thus we have a general Noether
  normalisation.\qed 
\end{proof}

The existence of a Noether normalisation for every affine algebra $\A=\P/\I$
follows now immediately from Theorem~\ref{thm:expombas} asserting the
existence of a Pommaret basis for $\I$.  Comparing with the existence proof of
Noether normalisations given in \cite{gp:singular}, we see that the search for
$\delta$-regular coordinates corresponds to putting the ideal $\I$ into
Noether position \cite[Def.~2.22]{wvv:comp}.

The converse of Corollary \ref{cor:noethernorm} is generally not true: even if
the variables are chosen such that $\kk[x_1,\dots,x_D]$ defines a Noether
normalisation of $\A$, this is not sufficient for concluding that the ideal
$\I$ possesses a Pommaret basis.  As we will show now, the existence of a
Pommaret basis is equivalent to a stronger property.  Since under the
assumptions of Proposition~\ref{prop:dim} $\kk[x_1,\dots,x_D]$ also defines a
Noether normalisation of $\P/\lt{\prec}{\I}$, it suffices to consider monomial
ideals.

\begin{definition}\label{def:qstab}
  A monomial ideal\/ $\I\subseteq\P$ is called\/ \emph{quasi-stable}, if it
  possesses a finite Pommaret basis.
\end{definition}

The reason for this terminology will become apparent in Section
\ref{sec:minres} when we consider stable ideals.  We now give several
equivalent algebraic characterisations of quasi-stable ideals which are
independent of the theory of Pommaret bases.  They will provide us with a
further criterion for $\delta$-regularity and also lead to a simple
description of an irredundant primary decomposition of such ideals.

\begin{proposition}\label{prop:qstab}
  Let\/ $\I\subseteq\P$ be a monomial ideal with\/ $\dim{\P/\I}=D$.  Then the
  following six statements are equivalent.
  \begin{description}[{\upshape (iii)}]
  \item[{\upshape (i)}] $\I$ is quasi-stable.
  \item[{\upshape (ii)}] The variable\/ $x_1$ is not a zero divisor
    for\/\footnote{See Section \ref{sec:regsat} for a more detailed discussion
      of the saturation $\Isat$.}  $\P/\Isat$ and for all\/ $1\leq j<D$ the
    variable\/ $x_{j+1}$ is not a zero divisor for\/
    $\P/\lspan{\I,x_1,\dots,x_j}^{\mathrm{sat}}$.
  \item[{\upshape (iii)}] We have
    $\I:x_1^\infty\subseteq\I:x_2^\infty\subseteq\cdots
    \subseteq\I:x_D^\infty$ and for all\/ $D<j\leq n$ an exponent\/ $k_j\geq1$
    exists such that\/ $x_j^{k_j}\in\I$.
  \item[{\upshape (iv)}] For all\/ $1\leq j\leq n$ the equality\/
    $\I:x_j^\infty=\I:\lspan{x_j,\dots,x_n}^\infty$ holds.
  \item[{\upshape (v)}] For every associated prime ideal\/
    $\pf\in\ass{(\P/\I)}$ an integer\/ $1\leq j\leq n$ exists such that\/
    $\pf=\lspan{x_j,\dots,x_n}$.
  \item[{\upshape (vi)}] If\/ $x^\mu\in\I$ and\/ $\mu_i>0$ for some\/
    $1\leq i<n$, then for each\/ $0<r\leq\mu_i$ and\/ $i<j\leq n$ an integer\/
    $s\geq0$ exists such that\/ $x^{\mu-r_i+s_j}\in\I$.
  \end{description}
\end{proposition}

\begin{proof}
  The equivalence of the statements (ii)--(v) was proven by Bermejo and
  Gimenez \cite[Prop.~3.2]{bg:scmr}; the equivalence of (iv) and (vi) was
  shown by Herzog et al. \cite[Prop.~2.2]{hpv:ext} (alternatively the
  equivalence of (i) and (vi) is an easy consequence of Lemma
  \ref{lem:stable}).  Bermejo and Gimenez \cite{bg:scmr} called ideals
  satisfying any of these conditions \emph{monomial ideals of nested type};
  Herzog et al. \cite{hpv:ext} spoke of \emph{ideals of Borel type} (yet
  another terminology used by Caviglia and Sbarra \cite{cs:reg} is
  \emph{weakly stable ideals}).\footnote{As usual, one must revert the
    ordering of the variables $x_1,\dots,x_n$ in order to recover the results
    of the given references.}  Thus it suffices to show that these concepts
  coincide with quasi-stability by proving the equivalence of (i) and (iii).
  
  Assume first that the ideal $\I$ is quasi-stable with Pommaret basis $\H$.
  The existence of a term $x_j^{k_j}\in\I$ for all $D<j\leq n$ follows then
  immediately from Proposition~\ref{prop:dim}.  Now consider a term
  $x^\mu\in\I:x_k^\infty\setminus\I$ for some $1\leq k\leq n$.  By definition
  of such a colon ideal, there exists an integer $\ell$ such that $x_k^\ell
  x^\mu\in\I$ and hence a generator $x^\nu\in\H$ such that $x^\nu\idiv{P}
  x_k^\ell x^\mu$.  If $\cls{\nu}>k$, then $\nu$ would also be an involutive
  divisor of $\mu$ contradicting the assumption $x^\mu\notin\I$.  Thus we find
  $\cls{\nu}\leq k$ and $\nu_k>\mu_k$.
  
  Next we consider for arbitrary exponents $m>0$ the non-multiplicative
  products $x_{k+1}^mx^\nu\in\I$.  For each $m$ a generator
  $x^{\rho^{(m)}}\in\H$ exists which involutively divides $x_{k+1}^mx^\nu$.
  By the same reasoning as above, $\cls{x^{\rho^{(m)}}}>k+1$ is not possible,
  as the Pommaret basis $\H$ is by definition involutively autoreduced.  This
  yields the estimate $\cls{\nu}\leq\cls{x^{\rho^{(m)}}}\leq k+1$.
  
  We claim now that there exists an integer $m_0$ such that
  $\rho^{(m)}=\rho^{(m_0)}$ for all $m\geq m_0$ and
  $\cls{x^{\rho^{(m_0)}}}=k+1$.  Indeed, if $\cls{x^{\rho^{(m)}}}<k+1$, then
  we must have $\rho^{(m)}_{k+1}=v_{k+1}+m$, since $x_{k+1}$ is not
  multiplicative for $x^{\rho^{(m)}}$.  Hence $x^{\rho^{(m)}}$ cannot be an
  involutive divisor of $x_{k+1}^{m+1}x^\nu$ and
  $\rho^{(m+1)}\notin\{\rho^{(1)},\dots,\rho^{(m)}\}$.  As the Pommaret basis
  $\H$ is a finite set, $\cls{x^{\rho^{(m_0)}}}=k+1$ for some value $m_0>0$.
  But then $x_{k+1}$ is multiplicative for $x^{\rho^{(m_0)}}$ and thus
  $x^{\rho^{(m_0)}}$ is trivially an involutive divisor of $x_{k+1}^mx^\nu$
  for all values $m\geq m_0$.
  
  Note that, by construction, the generator $x^{\rho^{(m_0)}}$ is also an
  involutive divisor of $x_{k+1}^{m_0}x^\mu$, as $x_k$ is multiplicative for
  it.  Hence this term must lie in $\I$ and consequently
  $x^\mu\in\I:x_{k+1}^\infty$.  Thus we may conclude that
  $\I:x_k^\infty\subseteq\I:x_{k+1}^\infty$ for all $1\leq k<n$.  This proves
  (iii).
  
  For the converse assume that (iii) holds and let $\B$ be the minimal basis
  of the ideal $\I$.  Let $x^\mu\in\B$ be an arbitrary term of class $k$.
  Then $x^\mu/x_k\in\I:x_k^\infty$.  By assumption, this means that also
  $x^\mu/x_k\in\I:x_\ell^\infty$ for any non-multiplicative index $\ell$.
  Hence for each term $x^\mu\in\B$ and for each value $\cls{(x^\mu)}<\ell\leq
  n$ there exists an integer $q_{\mu,\ell}$ such that
  $x_\ell^{q_{\mu,\ell}}x^\mu/x_k\notin\I$ but
  $x_\ell^{q_{\mu,\ell}+1}x^\mu/x_k\in\I$.  For the values
  $1\leq\ell\leq\cls{x^\mu}$ we set $q_{\mu,\ell}=0$.  Observe that if
  $x^\nu\in\B$ is a minimal generator dividing
  $x_\ell^{q_{\mu,\ell}+1}x^\mu/x_k$, then we find for the inverse
  lexicographic order that $x^\nu\prec_{\mathrm{invlex}}x^\mu$, since
  $\cls{(x^\nu)}\geq\cls{(x^\mu)}$ and $\nu_k<\mu_k$.
  
  Consider now the set
  \begin{equation}
    \H=\bigl\{x^{\mu+\rho}\mid x^\mu\in\B\wedge\forall 1\leq\ell\leq n:
                               0\leq\rho_\ell\leq q_{\mu,\ell}\bigr\}\;.
  \end{equation}
  We claim that $\H$ is a weak involutive completion of $\B$ and thus a weak
  Pommaret basis of $\I$.  In order to prove this assertion, we must show that
  each term $x^\lambda\in\I$ lies in the involutive cone of a member of $\H$.
  
  As $x^\lambda$ is assumed to be an element of $\I$, we can factor it as
  $x^\lambda=x^{\sigma^{(1)}} x^{\rho^{(1)}} x^{\mu^{(1)}}$ where
  $x^{\mu^{(1)}}\in\B$ is a minimal generator, $x^{\sigma^{(1)}}$ contains
  only multiplicative variables for $x^{\mu^{(1)}}$ and $x^{\rho^{(1)}}$ only
  non-multiplicative ones.  If $x^{\mu^{(1)}+\rho^{(1)}}\in\H$, then we are
  done, as obviously $\cls{\bigl(x^{\mu^{(1)}+\rho^{(1)}}\big)}=
  \cls{\bigl(x^{\mu^{(1)}}\bigr)}$ and hence all variables contained in
  $x^{\sigma^{(1)}}$ are multiplicative for $x^{\mu^{(1)}+\rho^{(1)}}$, too.
  
  Otherwise there exists at least one non-multiplicative variables $x_\ell$
  such that $\rho^{(1)}_\ell>q_{\mu^{(1)},\ell}$.  Any minimal generator
  $x^{\mu^{(2)}}\in\B$ dividing
  $x_\ell^{q_{\mu^{(1)},\ell}+1}x^{\mu^{(1)}}/x_k$ is also a divisor of
  $x^\lambda$ and we find a second factorisation $x^\lambda=x^{\sigma^{(2)}}
  x^{\rho^{(2)}} x^{\mu^{(2)}}$ where again $x^{\sigma^{(2)}}$ consists only
  of multiplicative and $x^{\rho^{(2)}}$ only of non-multiplicative variables
  for $x^{\mu^{(2)}}$.  If $x^{\mu^{(2)}+\rho^{(2)}}\in\H$, then we are done
  by the same argument as above; otherwise we iterate.
  
  According to the observation made above, the sequence
  $(x^{\mu^{(1)}},x^{\mu^{(2)}},\dots)$ of minimal generators constructed this
  way is strictly descending with respect to the inverse lexicographic order.
  However, the minimal basis $\B$ is a finite set and thus the iteration
  cannot go on infinitely.  As the iteration only stops, if there exists an
  involutive cone containing $x^\lambda$, the involutive span of $\H$ is
  indeed $\I$ and thus the ideal $\I$ quasi-stable.\qed
\end{proof}

\begin{remark}\label{rem:primdec1}
  Note that our considerations about standard pairs and the induced primary
  decomposition in the last section imply a simple direct proof of the
  implication ``(i)${}\Rightarrow{}$(v)'' in Proposition \ref{prop:qstab}.  If
  the ideal $\I$ is quasi-stable, then $\I$ admits a complementary Rees
  decomposition according to Proposition \ref{prop:pomdecomp}.  Together with
  Propositions \ref{prop:stanpair}, \ref{prop:standard} and Remark
  \ref{rem:stanpom}, this observation trivially implies that all associated
  prime ideals are of the form $\pf=\lspan{x_j,\dots,x_n}$.\bull
\end{remark}

\begin{lemma}\label{lem:siquasi}
  Let\/ $\I_1,\I_2\subseteq\P$ be two quasi-stable ideals.  Then the sum\/
  $\I_1+\I_2$, the product\/ $\I_1\cdot\I_2$ and the intersection\/
  $\I_1\cap\I_2$ are quasi-stable, too.  If\/ $\I\subseteq\P$ is a
  quasi-stable ideal, then the quotient\/ $\I:\J$ is again quasi-stable for
  arbitrary monomial ideals\/ $\J\subseteq\P$.
\end{lemma}

\begin{proof}
  For the sum $\I_1+\I_2$ the claim follows immediately from Remark 2.9 of
  Part I which states that the union $\H_1\cup\H_2$ of (weak) Pommaret bases
  $\H_k$ of $\I_k$ is a weak Pommaret basis of the sum $\I_1+\I_2$.
  Similarly, the case of both the product $\I_1\cdot\I_2$ and the intersection
  $\I_1\cap\I_2$ was settled in Remark 6.5 of Part I where for both ideals
  weak Pommaret bases were constructed.
  
  For the last assertion we use Part (vi) of Proposition~\ref{prop:qstab}.  If
  $\J$ is minimally generated by the monomials $m_1,\dots,m_r$, then
  $\I:\J=\bigcap_{k=1}^r\I:m_k$ and thus it suffice to consider the case that
  $\J$ is a principal ideal with generator $x^\nu$.  Assume that
  $x^\mu\in\I:x^\nu$ and that $\mu_i>0$.  Since $x^{\mu+\nu}$ lies in the
  quasi-stable ideal $\I$, we find for each $0<r\leq\mu_i$ and $i<j\leq n$ and
  integer $s\geq0$ exists such that $x^{\mu+\nu-r_i+s_j}\in\I$.  As
  $r\leq\mu_i$, this trivially implies that $x^{\mu-r_i+s_j}\in\I:x^\nu$.\qed
\end{proof}
  
\begin{remark}\label{rem:qstabcrit}
  Alternative proofs for Lemma \ref{lem:siquasi} were given by Cimpoea{\c s}
  \cite{mc:borel}.  There it was also noted that its final statement trivially
  implies Part (v) of Proposition~\ref{prop:qstab}, as any associated prime
  ideal $\pf$ of a quasi-stable ideal $\I$ is of the form $\pf=\I:x^\nu$ for
  some monomial $x^\nu$ and thus is also quasi-stable.  But the only
  quasi-stable prime ideals are obviously the ideals $\lspan{x_j,\dots,x_n}$.
  
  Above we actually proved that Part (iii) of Proposition~\ref{prop:qstab} may
  be replaced by the equivalent statement
  $\I:x_1^\infty\subseteq\I:x_2^\infty\subseteq\cdots \subseteq\I:x_n^\infty$
  which does not require a priori knowledge of $D$ (the dimension $D$ arises
  then trivially as the smallest value $k$ such that $\I:x_k^\infty=\P$,
  i.,\,e, for which $\I$ contains a minimal generator $x_k^\ell$ for some
  exponent $\ell>0$).  In this formulation it is straightforward to verify
  (iii) effectively: bases of the colon ideals $\I:x_k^\infty$ are easily
  obtained by setting $x_k=1$ in any basis of $\I$ and for monomial ideals it
  is trivial to check inclusion, as one must only compare their minimal bases.
  
  We furthermore note that if we have for some value $1\leq k\leq n$ an
  ascending chain $\I:x_1^\infty\subseteq\I:x_2^\infty\subseteq\cdots
  \subseteq\I:x_k^\infty$, then for each $1\leq j\leq k$ the minimal basis
  $\B_j$ of $\I:x_j^\infty$ lies in $\kk[x_{j+1},\dots,x_n]$.  Indeed, no
  element of $\B_j$ can depend on $x_j$.  Now assume that $x^\nu\in\B_j$
  satisfies $\cls{\nu}=\ell<j$.  Then $x_j^m x^\nu$ is a minimal generator of
  $\I$ for some suitable exponent $m\in\NN_0$.  This in turn implies that
  $x_j^m x^\nu/x_\ell^{\nu_\ell}\in \I:x_\ell^\infty \subseteq\I:x_j^\infty$
  and hence $x^\nu/x_\ell^{\nu_\ell}\in\I:x_j^\infty$ which contradicts our
  assumption that $x^\nu$ was a minimal generator.\bull
\end{remark}

Proposition~\ref{prop:qstab} provides us with a simple criterion for
$\delta$-regularity.  Actually, this new criterion is closely related to our
previous one based on a comparison of the Janet and Pommaret multiplicative
variables, as the proof of the following converse to Theorem~\ref{thm:dreg}
demonstrates.

\begin{proposition}\label{prop:qstabdreg}
  Let\/ $\I\subseteq\P$ be a monomial ideal and\/ $\B$ a finite, Pommaret
  autoreduced monomial basis of it.  If the ideal\/ $\I$ is not quasi-stable,
  then\/ $\norm{\B}_{J}>\norm{\B}_{P}$, i.\,e.\ for at least one generator in
  the basis\/ $\B$ a variable exists which is Janet but not Pommaret
  multiplicative.
\end{proposition}

\begin{proof}
  By Proposition \ref{prop:janpom} we have $\norm{\B}_{J}\geq\norm{\B}_{P}$.
  As $\I$ is not quasi-stable, there exists a minimal value $k$ such that
  $\I:x_k^\infty\nsubseteq\I:x_{k+1}^\infty$.  Let $x^\mu$ be a minimal
  generator of $\I:x_k^\infty$ which is not contained in $\I:x_{k+1}^\infty$.
  Then for a suitable exponent $m\in\NN_0$ the term $x^{\bar\mu}=x_k^mx^\mu$
  is a minimal generator of $\I$ and hence contained in $\B$.
  
  We claim now that $\B$ contains a generator for which $x_{k+1}$ is Janet but
  not Pommaret multiplicative.  If $x_{k+1}\in\mult{X}{J,\B}{x^{\bar\mu}}$,
  then we are done, as according to Remark \ref{rem:qstabcrit}
  $\cls{\bar\mu}=k$ and hence $x_{k+1}\notin\mult{X}{P}{x^{\bar\mu}}$.
  Otherwise $\B$ contains a term $x^\nu$ such that $\nu_\ell=\mu_\ell$ for
  $k+1<\ell\leq n$ and $\nu_{k+1}>\mu_{k+1}$.  If several generators with this
  property exist in $\B$, we choose one for which $\nu_{k+1}$ takes a maximal
  value so that we have $x_{k+1}\in\mult{X}{J,\B}{x^\nu}$ by definition of the
  Janet division.  If $\cls{\nu}<k+1$, we are again done, as then
  $x_{k+1}\notin\mult{X}{P}{x^\nu}$.  Now assume that $\cls{\nu}=k+1$ and
  consider the term $x^\rho=x^\nu/x_{k+1}^{\nu_{k+1}}$.  Obviously,
  $x^\rho\in\I:x_{k+1}^\infty$ contradicting our assumption
  $x^\mu\notin\I:x_{k+1}^\infty$ since $x^\rho\mid x^\mu$.  Hence this case
  cannot arise.\qed
\end{proof}

Proposition~\ref{prop:qstab} (v) furthermore allows us to formulate for
monomial ideals a converse to Corollary \ref{cor:noethernorm}.  It is
equivalent to \cite[Prop.~3.6]{bg:scmr} and shows that a Pommaret basis of a
monomial ideal $\I$ induces \emph{simultaneous} Noether normalisations of all
primary components of the ideal.

\begin{corollary}\label{cor:noether}
  Let\/ $\I\subseteq\P$ be a monomial ideal with\/ $\dim{\P/\I}=D$.
  Furthermore, let\/ $\I=\qf_1\cap\cdots\cap\qf_r$ be an irredundant monomial
  primary decomposition with\/ $D_j=\dim{\P/\qf_j}$ for\/ $1\leq j\leq r$.
  The ideal\/ $\I$ is quasi-stable, if and only if\/ $\kk[x_1,\dots,x_D]$
  defines a Noether normalisation of\/ $\P/\I$ and\/ $\kk[x_1,\dots,x_{D_j}]$
  one of\/ $\P/\qf_j$ for each primary component\/ $\qf_j$.
\end{corollary}

\begin{proof}
  By assumption, each $\qf_j$ is a monomial primary ideal.  This implies that
  $\kk[x_1,\dots,x_{D_j}]$ defines a Noether normalisation of $\P/\qf_j$, if
  and only if the associated prime ideal is
  $\sqrt{\qf_j}=\lspan{x_{D_j+1},\dots,x_n}$.  Now the assertion follows
  immediately from Proposition~\ref{prop:qstab} (v).\qed
\end{proof}

We may also exploit Proposition~\ref{prop:qstab} for actually deriving an
irredundant primary decomposition $\I=\qf_1\cap\cdots\cap\qf_t$ with monomial
ideals $\qf_j$ for an arbitrary quasi-stable ideal $\I$.\footnote{The
  following construction is joint work with M. Hausdorf and M. Sahbi and has
  already appeared in \cite{wms:delta}.}  Bermejo and Gimenez
\cite[Rem.~3.3]{bg:scmr} noted that their proof of the implication
``(v)${}\Rightarrow{}$(iv)'' in Proposition \ref{prop:qstab} has some simple
consequences for the primary ideals $\qf_j$.  Let again $D=\dim{\P/\I}$.  Then
$\pf=\lspan{x_{D+1},\dots,x_n}$ is the unique minimal prime ideal associated
to $\I$ and the corresponding unique primary component is given by
$\I:x_D^\infty$ (if $D=0$, then obviously $\I$ is already a primary ideal).
More generally, we find for any $1\leq k\leq D$ that
\begin{equation}
  \I:x_k^\infty=\bigcap_{\pf_j\subseteq\lspan{x_{k+1},\dots,x_n}} \qf_j
\end{equation}
where $\pf_j=\sqrt{\qf_j}$ is the corresponding associated prime ideal.  Based
on these observations, an irredundant primary decomposition can be constructed
by working backwards through the sequence
$\I\subseteq\I:x_1^\infty\subseteq\I:x_2^\infty\subseteq\cdots
\subseteq\I:x_n^\infty$.

Let $d=\depth{\P/\I}$, i.\,e.\ $d+1$ is the minimal class of a generator in
the Pommaret basis $\H$ of $\I$ according to Proposition
\ref{prop:depth}.\footnote{Note that for determining the depth $d$ in the case
  of a quasi-stable ideal, it is not necessary to compute the Pommaret basis:
  since multiplication with a non-multiplicative variable never decreases the
  class, $d+1$ is also the minimal class of a minimal generator.}  For $1\leq
k\leq D$ we set $s_k=\min{\{s\mid \I:x_k^s=\I:x_k^{s+1}\}}$, i.\,e.\ $s_k$ is
the highest $x_k$-degree of a minimal generator of $\I$.  Then we introduce
the ideals $\J_k=\I+\lspan{x_{k+1}^{s_{k+1}},\dots,x_D^{s_D}}$ and
\begin{equation}\label{eq:qk}
  \qf_k=\J_k:x_k^\infty=\I:x_k^\infty+
        \lspan{x_{k+1}^{s_{k+1}},\dots,x_D^{s_D}}\;.
\end{equation}
It is easy to see that all the ideals $\J_k$ are again quasi-stable provided
the ideal $\I$ is quasi-stable (this follows immediately from
Proposition~\ref{prop:qstab} and the fact that in this case
$(\I:x_i^\infty):x_j^\infty=\I:x_j^\infty$ for $i<j$).  For notational
simplicity we formally define $\I:x_0^\infty=\I$ and
$\qf_0=\J_0=\I+\lspan{x_1^{s_1},\dots,x_D^{s_D}}$.  Since obviously
$\dim{\P/\J_k}=k$ for $0\leq k\leq D$, it follows from the considerations
above that $\qf_k$ is an $\lspan{x_{k+1},\dots,x_n}$-primary ideal.

\begin{proposition}\label{prop:primdec}
  Let\/ $\I\subseteq\P$ be a quasi-stable ideal.  Then\/
  $\I=\bigcap_{k=d}^D\qf_k$ is a primary decomposition.  Eliminating all
  ideals\/ $\qf_k$ where\/ $\I:x_k^\infty=\I:x_{k+1}^\infty$ makes it an
  irredundant decomposition.
\end{proposition}

\begin{proof}
  We first show that the equality $\I:x_k^\infty=\bigcap_{\ell=k}^D\qf_\ell$
  holds or equivalently that $\I:x_k^\infty=\qf_k\cap(\I:x_{k+1}^\infty)$ for
  $0\leq k\leq n$; for $k=d$ this represents the first statement of the
  proposition, since obviously $\I:x_0^\infty=\cdots=\I:x_d^\infty=\I$.  By
  definition of the value $s_{k+1}$, we have that \cite[Lemma
  3.3.6]{gp:singular}
  \begin{equation}
    \I:x_k^\infty=\bigl(\I:x_k^\infty+\lspan{x_{k+1}^{s_{k+1}}}\bigr)\cap
                  \bigl((\I:x_k^\infty):x_{k+1}^\infty\bigr)\;.
  \end{equation}
  The second factor obviously equals $\I:x_{k+1}^\infty$.  To the first one we
  apply the same construction and decompose
  \begin{equation}
    \begin{split}
      \I:x_k^\infty&+\lspan{x_{k+1}^{s_{k+1}}}=\\
      &=\bigl(\I:x_k^\infty+\lspan{x_{k+1}^{s_{k+1}},
                                   x_{k+2}^{s_{k+2}}}\bigr)\cap
          \bigl((\I:x_k^\infty+
                 \lspan{x_{k+1}^{s_{k+1}}}):x_{k+2}^\infty\bigr)\\
      &=\bigl(\I:x_k^\infty+\lspan{x_{k+1}^{s_{k+1}},
                                   x_{k+2}^{s_{k+2}}}\bigr)\cap
        \bigl(\I:x_{k+2}^\infty+\lspan{x_{k+1}^{s_{k+1}}}\bigr)\;.
    \end{split}
  \end{equation}
  Continuing in this manner, we arrive at a decomposition
  \begin{equation}
    \I:x_k^\infty=\qf_k\cap\dots\cap(\I:x_{k+1}^\infty)
  \end{equation}
  where the dots represent factors of the form
  $\I:x_\ell^\infty+\lspan{x_{k+1}^{s_{k+1}},\dots,x_{\ell-1}^{s_{\ell-1}}}$
  with $\ell\geq k+2$.  Since we assume that $\I$ is quasi-stable,
  $\I:x_{k+1}^\infty$ is contained in each of these factors and we may omit
  them which proves our claim.
  
  In the thus obtained primary decomposition of $\I$ the radicals of all
  appearing primary ideals are pairwise different.  Furthermore, it is obvious
  that $\qf_k$ is redundant whenever $\I:x_k^\infty=\I:x_{k+1}^\infty$.  Thus
  there only remains to prove that all the other primary ideals $\qf_k$ are
  indeed necessary.  Assume that $\I:x_k^\infty\subsetneq\I:x_{k+1}^\infty$
  (which is in particular the case for $k<d$).  Then there exists a minimal
  generator $x^\mu$ of $\I:x_{k+1}^\infty$ which is not contained in
  $\I:x_k^\infty$.  Consider the monomial $x_k^{s_k}x^\mu$.  It cannot lie in
  $\I:x_k^\infty$, as otherwise already $x^\mu\in\I:x_k^\infty$, and thus it
  also cannot be contained in $\qf_k$ (since we showed above that
  $\I:x_k^\infty=\qf_k\cap(\I:x_{k+1}^\infty)$).  On the other hand we find
  that $x_k^{s_k}x^\mu\in\qf_\ell$ for all $\ell>k$ since then
  $\I:x_{k+1}^\infty\subseteq\qf_\ell$ and for all $\ell<k$ since then
  $\lspan{x_k^{s_k}}\subseteq\qf_\ell$.  Hence $\qf_k$ is not redundant.\qed
\end{proof}

According to Lemma \ref{lem:siquasi}, the quotient ideals $\I:x_k^\infty$ are
again quasi-stable.  It is straightforward to obtain Pommaret bases for them.
We disjointly decompose the monomial Pommaret basis
$\H=\H_1\cup\cdots\cup\H_n$ where $\H_k$ contains all generators of class $k$.
Furthermore, we write $\H'_k$ for the set obtained by setting $x_k=1$ in each
generator in $\H_k$.

\begin{lemma}
  For any\/ $1\leq k\leq n$ the set\/
  $\H'=\H'_k\cup\bigcup_{\ell=k+1}^n\H_\ell$ is a weak Pommaret basis of the
  colon ideal\/ $\I:x_k^\infty$.
\end{lemma}
  
\begin{proof}
  We first show that $\H'$ is an involutive set.  By definition of the
  Pommaret division, it is obvious that the subset
  $\bigcup_{\ell=k+1}^n\H_\ell$ is involutive. Thus there only remains to
  consider the non-multiplicative products of the members of $\H'_k$.  Take
  $x^\mu\in\H'_k$ and let $x_\ell$ be a non-multiplicative variable for it.
  Obviously, there exists an $m>0$ such that $x_k^mx^\mu\in\H_k$ and hence a
  generator $x^\nu\in\bigcup_{\ell=k}^n\H_\ell$ such that $x_\ell x_k^mx^\mu$
  lies in the involutive cone $\cone{P}{x^\nu}$.  Writing $x_\ell
  x_k^mx^\mu=x^{\rho+\nu}$, we distinguish two cases.  If $\cls{\nu}>k$, then
  $\rho_k=m$ and we can divide by $x_k^m$ in order to obtain an involutive
  standard representation of $x_\ell x^\mu$ with respect to $\H'$.  If
  $\cls{\nu}=k$, then the multi index $\rho$ is of the form $r_k$, i.\,e.\ 
  only the $k$th entry is different from zero, and we even find that $x_\ell
  x^\mu=x^\nu/x_k^r\in\H'_k$.
  
  Thus there only remains to prove that $\H'$ is actually a generating set for
  $\I:x_k^\infty$. For this we first note that the Pommaret basis of a
  quasi-stable ideal contains a generator of class $k$ only, if there is a
  minimal generator of class $k$, as applying the monomial completion
  Algorithm 2 of Part I to the minimal basis adds only non-multiplicative
  multiples of the minimal generators (and these are trivially of the same
  class).  By Remark \ref{rem:primdec1}, all minimal generators of
  $\I:x_k^\infty$ have at least class $k+1$.  Thus setting $x_k=1$ in any
  member of $\bigcup_{\ell=1}^{k-1}\H_\ell$ can never produce a minimal
  generator of $\I:x_k^\infty$ and thus $\H'$ is a weak involutive completion
  of the minimal basis of $\I:x_k^\infty$.  According to Proposition 2.8 of
  Part I, an involutive autoreduction yields a strong basis.\qed
\end{proof}
  
The ideals $\lspan{x_{k+1}^{s_{k+1}},\dots,x_D^{s_D}}$ are obviously
irreducible and for $k\geq d$ exactly of the form that they possess a Pommaret
basis as discussed in Example 2.12 of Part I.  There we also gave an explicit
Pommaret basis for such an ideal.  Since according to Remark 2.9 of Part I the
union of two (weak) Pommaret bases of two monomial ideals $\I_1$, $\I_2$
yields a weak Pommaret basis of $\I_1+\I_2$, we obtain this way easily weak
Pommaret bases for all primary ideals $\qf_k$ appearing in the irredundant
decomposition of Proposition \ref{prop:primdec}.

Thus the crucial information for obtaining an irredundant primary
decomposition of a quasi-stable ideal $\I$ is where ``jumps'' are located,
i.\,e.\ where $\I:x_k^\infty\subsetneq\I:x_{k+1}^\infty$.  Since these ideals
are quasi-stable, the positions of the jumps are determined by their depths.
A chain with all the jumps is obtained by the following simple recipe: set
$\I_0=\I$ and define $\I_{k+1}=\I_k:x_{d_k}^\infty$ where $d_k=\depth{\I_k}$.
This leads to the so-called \emph{sequential chain} of $\I$:
\begin{equation}\label{eq:seqchain}
  \I_0=\I\subsetneq\I_1\subsetneq\cdots\subsetneq\I_r=\P;.
\end{equation}

\begin{remark}
  With the help of the sequential chain (\ref{eq:seqchain}) one can also show
  straightforwardly that any quasi-stable ideal is \emph{sequentially
    Cohen-Macaulay} \cite[Cor.~2.5]{hpv:ext} (recall that the algebra
  $\A=\P/\I$ is sequentially Cohen-Macaulay \cite{rps:cca}, if a chain
  $\I_0=\I\subset\I_1\subset\cdots\subset\I_r=\P$ exists such that all
  quotients $\I_{k+1}/\I_k$ are Cohen-Macaulay and their dimensions are
  ascending: $\dim{(\I_k/\I_{k-1})}<\dim{(\I_{k+1}/\I_k)}$).
  
  Indeed, consider the ideal $\J_k=\I_k\cap\kk[x_{d_k},\dots,x_n]$.  By
  Remark~\ref{rem:qstabcrit}, the minimal generators of $\J_k$ are the same as
  the minimal ones of $\I_k$; furthermore, by Proposition \ref{prop:qstab}(iv)
  $\J_k^{\mathrm{sat}}=\J_k:x_{d_k}^\infty$.  Hence we find that
  $\I_{k+1}=\lspan{\J_k^{\mathrm{sat}}}_{\P}$ and
  \begin{equation}
    \I_{k+1}/\I_k\cong (\J_k^{\mathrm{sat}}/\J_k)[x_1,\dots,x_{d_k-1}]\;.
  \end{equation}
  Since the factor ring $\J_k^{\mathrm{sat}}/\J_k$ is trivially finite (as a
  $\kk$-linear space), the quotient $\I_{k+1}/\I_k$ is thus a
  $(d_k-1)$-dimensional Cohen-Macaulay module.\bull
\end{remark}

\section{Syzygies for Involutive Bases}\label{sec:syz}

Gr\"obner bases are a very useful tool in syzygy theory.  A central result is
\emph{Schreyer's Theorem} \cite{al:gb,fos:syz} that the standard
representations of the $S$-polynomials between the elements of a Gr\"obner
basis determine directly a Gr\"obner basis of the first syzygy module with
respect to an appropriately chosen term order.  Now we study the use of
involutive bases in this context.

In Part I we introduced involutive bases only for ideals, but the extension to
submodules of free modules $\P^m$ is trivial.  We represent elements of $\P^m$
as vectors $\fv=(f_1,\dots,f_m)$ with $f_\alpha\in\P$.  The standard basis of
$\P^m$ consists of the unit vectors $\ev_\alpha$ with
$e_{\alpha\beta}=\delta_{\alpha\beta}$ and $1\leq\alpha\leq m$; thus
$\fv=f_1\ev_1+\cdots+f_m\ev_m$.  Now a term $\tv$ is a vector of the form
$\tv=t\ev_\alpha$ for some $\alpha$ and with $t\in\TT$ a term in $\P$.  We
denote the set of all terms by $\TT^m$; it is a monoid module over $\TT$.

Let $\H\subset\P^m$ be a finite set, $\prec$ a term order on $\TT^m$ and $L$
an involutive division on $\Nn$.  We divide $\H$ into $m$ disjoint sets
$\H_\alpha=\bigl\{\hv\in\H\mid\lt{\prec}{\hv}=t\ev_\alpha,\ t\in\TT\bigr\}$.
This leads naturally to $m$ sets $\N_\alpha=\bigl\{\mu\in\Nn\mid
x^\mu\ev_\alpha\in\lt{\prec}{\H_\alpha}\bigr\}$.  If $\hv\in\H_\alpha$, we
assign the multiplicative variables $\mult{X}{L,\H,\prec}{\hv}=\bigl\{x_i\mid
i\in\mult{N}{L,\N_\alpha}{\le{\prec}{\hv}}\bigr\}$.  The involutive span
$\ispan{\H}{L,\prec}$ is defined by an obvious generalisation of the old
definition in Part I.

Let $\H=\{\hv_1,\dots,\hv_s\}$ be an involutive basis of the submodule
$\M\subseteq\P^m$.  Take an arbitrary element $\hv_\alpha\in\H$ and choose an
arbitrary non-multiplicative variable
$x_k\in\nmult{X}{L,\H,\prec}{\hv_\alpha}$ of it.  By the results of Part I, we
can determine for each generator $\hv_\beta\in\H$ with an involutive normal
form algorithm a unique polynomial
$P^{(\alpha;k)}_\beta\in\kk[\mult{X}{L,\H,\prec}{\hv_\beta}]$ such that
$x_k\hv_\alpha=\sum_{\beta=1}^sP^{(\alpha;k)}_\beta\hv_\beta$.  To this
relation corresponds the syzygy
\begin{equation}\label{eq:invsyz}
  \Sv_{\alpha;k}=x_k\ev_\alpha-
                 \sum_{\beta=1}^sP^{(\alpha;k)}_\beta\ev_\beta\in\P^s\;.
\end{equation}
We denote the set of all thus obtained syzygies by
\begin{equation}\label{eq:Hsyz}
  \H_{\mathrm{Syz}}=\left\{\Sv_{\alpha;k}\mid 1\leq\alpha\leq s;\ 
    x_k\in\nmult{X}{L,\H,\prec}{\hv_\alpha}\right\}\;.
\end{equation}

\begin{lemma}\label{lem:multsyz}
  Let\/ $\H$ be an involutive basis for the involutive division\/ $L$ and the
  term order\/ $\prec$.  If\/ $\Sv=\sum_{\beta=1}^s S_\beta\ev_\beta$ is an
  arbitrary syzygy in the module\/ $\syz{}{\H}$ with\/
  $S_\beta\in\kk[\mult{X}{L,\H,\prec}{\hv_\beta}]$ for all\/ $1\leq\beta\leq
  s$, then\/ $\Sv=0$.
\end{lemma}

\begin{proof}
  By definition of a syzygy, $\sum_{\beta=1}^s S_\beta\hv_\beta=0$.  As the
  involutive basis $\H$ is involutively head autoreduced, each element
  $\fv\in\lspan{\H}$ possesses a unique involutive normal form.  In
  particular, this holds for $0\in\lspan{\H}$.  Thus either $\Sv=0$ or
  $S_\beta\notin\kk[\mult{X}{L,\H,\prec}{\hv_\beta}]$ for at least one
  $\beta$.\qed
\end{proof}

A fundamental ingredient of Schreyer's Theorem is the term order $\prec_{\F}$
on $\TT^s$ induced by an arbitrary finite set
$\F=\{\fv_1,\dots,\fv_s\}\subset\P^m$ and an arbitrary term order $\prec$ on
$\TT^m$: given two terms $\sv=s\ev_\sigma$ and $\tv=t\ev_\tau$, we set
$\sv\prec_{\F}\tv$, if either
$\lt{\prec}{(s\fv_\sigma)}\prec\lt{\prec}{(t\fv_\tau)}$ or
$\lt{\prec}{(s\fv_\sigma)}=\lt{\prec}{(t\fv_\tau)}$ and $\tau<\sigma$.

\begin{corollary}\label{cor:gensyz}
  If\/ $\H\subset\P$ is an involutive basis, then the set\/
  $\H_{\mathrm{Syz}}$ generates the syzygy module\/ $\syz{}{\H}$.
\end{corollary}

\begin{proof}
  Let $\Sv=\sum_{\beta=1}^s S_\beta\ev_\beta$ by an arbitrary non-vanishing
  syzygy in $\syz{}{\H}$.  By Lemma \ref{lem:multsyz}, at least one of the
  coefficients $S_\beta$ must contain a term $x^\mu$ with a non-multiplicative
  variable $x_j\in\nmult{X}{L,\H,\prec}{\hv_\beta}$.  Let $cx^\mu\ev_\beta$ be
  the maximal such term with respect to the term order $\prec_{\H}$ and $j$
  the maximal non-multiplicative index with $\mu_j>0$.  Then we eliminate this
  term by computing $\Sv'=\Sv-cx^{\mu-1_j}\Sv_{\beta;j}$.  If $\Sv'\neq0$, we
  iterate.  Since all new terms introduced by the subtraction are smaller than
  the eliminated term with respect to $\prec_{\H}$, we must reach zero after a
  finite number of steps.  Thus this computation leads to a representation of
  $\Sv$ as a linear combination of elements of $\H_{\mathrm{Syz}}$.\qed
\end{proof}

Let $\H=\{\hv_1,\cdots,\hv_s\}$ be an involutive basis and thus a Gr\"obner
basis for the term order $\prec$.  Without loss of generality we may assume
that $\H$ is a monic basis.  Set $\tv_\alpha=\lt{\prec}{\hv_\alpha}$ and
$\tv_{\alpha\beta}=\lcm{\tv_\alpha}{\tv_\beta}$.  We have for every
$S$-polynomial a standard representation
$\Spolym{\prec}{\hv_\alpha}{\hv_\beta}= \sum_{\gamma=1}^s
f_{\alpha\beta\gamma}\hv_\gamma$ where the polynomials
$f_{\alpha\beta\gamma}\in\P$ satisfy
$\lt{\prec}{\bigl(\Spolym{\prec}{\hv_\alpha}{\hv_\beta}\bigr)}
\succeq\lt{\prec}{(f_{\alpha\beta\gamma}\hv_\gamma)}$ for $1\leq\gamma\leq s$.
Setting $\fv_{\alpha\beta}=\sum_{\gamma=1}^s f_{\alpha\beta\gamma}\ev_\gamma$,
we introduce for $\alpha\neq\beta$ the syzygy
\begin{equation}\label{eq:schreyersyz}
  \Sv_{\alpha\beta}=\frac{\tv_{\alpha\beta}}{\tv_\alpha}\ev_\alpha-
           \frac{\tv_{\alpha\beta}}{\tv_\beta}\ev_\beta-
           \fv_{\alpha\beta}\;.
\end{equation}
Schreyer's Theorem asserts that the set
$\H_{\mathrm{Schreyer}}=\{\Sv_{\alpha\beta}\mid 1\leq\alpha<\beta\leq s\}$ of
all these syzygies is a Gr\"obner basis of the first syzygy module
$\syz{}{\H}$ for the induced term order $\prec_{\H}$.

We denote by
$\tilde{\Sv}_{\alpha\beta}=\tfrac{\tv_{\alpha\beta}}{\tv_\alpha}\ev_\alpha-
\tfrac{\tv_{\alpha\beta}}{\tv_\beta}\ev_\beta$ the syzygy of the leading terms
corresponding to $\Sv_{\alpha\beta}$ and if
$\S\subseteq\H_{\mathrm{Schreyer}}$ is a set of syzygies, $\tilde\S$ contains
the corresponding syzygies of the leading terms.

\begin{lemma}\label{lem:bbc2}
  Let\/ $\S\subseteq\H_{\mathrm{Schreyer}}$ be such that\/ $\tilde\S$
  generates\/ $\syz{}{\lt{\prec}{\H}}$.  Then\/ $\S$ generates\/ $\syz{}{\H}$.
  Assume furthermore that the three pairwise distinct indices\/ $\alpha$,
  $\beta$, $\gamma$ are such that\footnote{If $\alpha>\beta$, then we
    understand that $\Sv_{\beta\alpha}\in\S$ etc.}\/
  $\Sv_{\alpha\beta},\Sv_{\beta\gamma},\Sv_{\alpha\gamma}\in\S$ and\/
  $\tv_\gamma\mid \tv_{\alpha\beta}$.  Then the smaller set\/
  $\S\setminus\{\Sv_{\alpha\beta}\}$ still generates\/ $\syz{}{\H}$.
\end{lemma}

\begin{proof}
  It is a classical result in the theory of Gr\"obner bases that
  $\tilde\S\setminus\{\tilde{\Sv}_{\alpha\beta}\}$ still generates\/
  $\syz{}{\lt{\prec}{\H}}$.  In fact, this is the basic property underlying
  Buchberger's second criterion for avoiding redundant $S$-polynomials.  Thus
  it suffices to show the first assertion; the second one is a simple
  corollary.
  
  Let $\Rv=\sum_{\alpha=1}^s R_\alpha\ev_\alpha\in\syz{}{\H}$ be an arbitrary
  syzygy of the full generators and set
  $\tv_{\Rv}=\max_{\prec}{\bigl\{\lt{\prec}{(R_\alpha\hv_\alpha)}\mid
  1\leq\alpha\leq s\bigr\}}$.  Then
  \begin{equation}
    \tilde{\Rv}=\!\!\sum_{\lt{\prec}{(R_\alpha\hv_\alpha)}=\tv_R}\!\!
    \lt{\prec}{(R_\alpha\hv_\alpha)}\in\syz{}{\lt{\prec}{\H}}\;.
  \end{equation}
  According to our assumption $\tilde\S$ is a generating set of
  $\syz{}{\lt{\prec}{\H}}$, so that we may write
  $\tilde{\Rv}=\sum_{\tilde{\Sv}\in\tilde\S}a_{\tilde{\Sv}}\tilde{\Sv}$ for
  some coefficients $a_{\tilde{\Sv}}\in\P$.  Let us now consider the syzygy
  $\Rv'=\Rv-\sum_{\Sv\in\S}a_{\tilde{\Sv}}\Sv$.  Obviously,
  $\tv_{\Rv'}\prec\tv_{\Rv}$.  By iteration we obtain thus in a finite number
  of steps a representation $\Rv=\sum_{\Sv\in\S}b_{\Sv}\Sv$ and thus $\S$
  generates the module $\syz{}{\H}$.\qed
\end{proof}

As a consequence of this simple lemma, we can now show that each involutive
basis yields immediately a Gr\"obner basis of the first syzygy module.  In
fact, this basis is automatically computed during the determination of the
involutive basis with the completion Algorithm~3 of Part I.  This is
completely analogously to the automatic determination of
$\H_{\mathrm{Schreyer}}$ with the Buchberger algorithm.

\begin{theorem}\label{thm:invschreyer}
  Let\/ $\H$ be an involutive basis for the involutive division\/ $L$ and the
  term order\/ $\prec$.  Then the set\/ $\H_{\mathrm{Syz}}$ is a Gr\"obner
  basis of the syzygy module\/ $\syz{}{\H}$ for the term order\/ $\prec_{\H}$.
\end{theorem}

\begin{proof}
  Without loss of generality, we may assume that $\H$ is a monic basis,
  i.\,e.\ all leading coefficients are $1$.  Let
  $\Sv_{\alpha;k}\in\H_{\mathrm{Syz}}$.  As $\H$ is an involutive basis, the
  unique polynomials $P^{(\alpha;k)}_\beta$ in (\ref{eq:invsyz}) satisfy
  $\lt{\prec}{(P^{(\alpha;k)}_\beta\hv_\beta)}\preceq
  \lt{\prec}{(x_k\hv_\alpha)}$ and there exists only one index $\bar\beta$
  such that $\lt{\prec}{(P^{(\alpha;k)}_{\bar\beta}\hv_{\bar\beta})}=
  \lt{\prec}{(x_k\hv_\alpha)}$.  It is easy to see that we have
  $\Sv_{\alpha;k}=\Sv_{\alpha\bar\beta}$.  Thus
  $\H_{\mathrm{Syz}}\subseteq\H_{\mathrm{Schreyer}}$.
  
  Let $\Sv_{\alpha\beta}\in\H_{\mathrm{Schreyer}}\setminus\H_{\mathrm{Syz}}$
  be an arbitrary syzygy.  We prove first that the set
  $\H_{\mathrm{Schreyer}}\setminus\{\Sv_{\alpha\beta}\}$ still generates
  $\syz{}{\H}$.  Any syzygy in $\H_{\mathrm{Schreyer}}$ has the form
  $\Sv_{\alpha\beta}=x^\mu\ev_\alpha-x^\nu\ev_\beta+\Rv_{\alpha\beta}$.  By
  construction, one of the monomials $x^\mu$ and $x^\nu$ must contain a
  non-multiplicative variable $x_k$ for $\hv_\alpha$ or $\hv_\beta$,
  respectively.  Without loss of generality, we assume that
  $x_k\in\nmult{X}{L,\H,\prec}{\hv_\alpha}$ and $\mu_k>0$.  This implies that
  $\H_{\mathrm{Syz}}$ contains the syzygy $\Sv_{\alpha;k}$.  As shown above, a
  unique index $\gamma\neq\beta$ exists such that
  $\Sv_{\alpha;k}=\Sv_{\alpha\gamma}$.
  
  Let $\Sv_{\alpha\gamma}=x_k\ev_\alpha-x^\rho\ev_\gamma+\Rv_{\alpha\gamma}$.
  By construction, $x^\rho\tv_\gamma=x_k\tv_\alpha$ divides
  $x^\mu\tv_\alpha=\tv_{\alpha\beta}$.  Thus $\tv_\gamma\mid\tv_{\alpha\beta}$
  and by Lemma \ref{lem:bbc2} the set
  $\H_{\mathrm{Schreyer}}\setminus\{\Sv_{\alpha\beta}\}$ still generates
  $\syz{}{\H}$.  If we try to iterate this argument, we encounter the
  following problem.  In order to be able to eliminate $\Sv_{\alpha\beta}$ we
  need both $\Sv_{\alpha\gamma}$ and $\Sv_{\beta\gamma}$ in the remaining set.
  For $\Sv_{\alpha\gamma}\in\H_{\mathrm{Syz}}$, this is always guaranteed.
  But we know nothing about $\Sv_{\beta\gamma}$ and, if it is not an element
  of $\H_{\mathrm{Syz}}$, it could have been removed in an earlier iteration.
  
  We claim that with respect to the term order $\prec_\H$ the term
  $\lt{\prec_{\H}}{\Sv_{\alpha\beta}}$ is greater than both
  $\lt{\prec_{\H}}{\Sv_{\alpha\gamma}}$ and
  $\lt{\prec_{\H}}{\Sv_{\beta\gamma}}$.  Without loss of generality, we may
  assume for simplicity that $\alpha<\beta<\gamma$, as the syzygies
  $\Sv_{\alpha\beta}$ and $\Sv_{\beta\alpha}$ differ only by a sign.  Thus
  $\lt{\prec_{\H}}{\Sv_{\alpha\beta}}=
  \tfrac{\tv_{\alpha\beta}}{\tv_\alpha}\ev_\alpha$ and similarly for
  $\Sv_{\alpha\gamma}$ and $\Sv_{\beta\gamma}$.  Furthermore,
  $\tv_\gamma\mid\tv_{\alpha\beta}$ trivially implies
  $\tv_{\alpha\gamma}\mid\tv_{\alpha\beta}$ and hence
  $\tv_{\alpha\gamma}\prec\tv_{\alpha\beta}$ for any term order $\prec$.
  Obviously, the same holds for $\tv_{\beta\gamma}$.  Now a straightforward
  application of the definition of the term order $\prec_\H$ proves our claim.
  
  Thus if we always remove the syzygy
  $\Sv_{\alpha\beta}\in\H_{\mathrm{Schreyer}}\setminus\H_{\mathrm{Syz}}$ whose
  leading term is maximal with respect to the term order $\prec_{\H}$, it can
  never happen that the syzygy $\Sv_{\beta\gamma}$ required for the
  application of Lemma \ref{lem:bbc2} has already been eliminated earlier and
  $\H_{\mathrm{Syz}}$ is a generating set of $\syz{}{\H}$.
  
  It is a simple corollary to Schreyer's theorem that $\H_{\mathrm{Syz}}$ is
  even a Gr\"obner basis of $\syz{}{\H}$.  Indeed, we know that
  $\H_{\mathrm{Schreyer}}$ is a Gr\"obner basis of $\syz{}{\H}$ for the term
  order $\prec_\H$ and it follows from our considerations above that whenever
  we remove a syzygy $\Sv_{\alpha\beta}$ we still have in the remaining set at
  least one syzygy whose leading term divides
  $\lt{\prec_{\H}}{\Sv_{\alpha\beta}}$.  Thus we find
  \begin{equation}
    \lspan{\lt{\prec_{\H}}{(\H_{\mathrm{Syz}})}}=
    \lspan{\lt{\prec_{\H}}{(\H_{\mathrm{Schreyer}})}}=
    \lt{\prec_{\H}}{\syz{}{\H}}
  \end{equation}
  which proves our assertion.\qed
\end{proof}

This result is not completely satisfying, as it only yields a Gr\"obner and
not an involutive basis of the syzygy module.  The latter seems to be hard to
achieve for arbitrary divisions $L$.  For some divisions it is possible with a
little effort.  The key is that in the order $\prec_{\H}$ the numbering of the
generators in $\H$ is important and we must choose the right one.  For this
purpose we slightly generalise a construction of Plesken and Robertz
\cite{pr:janet} for the special case of a Janet basis.

We associate a graph with each involutive basis $\H$.  Its vertices are given
by the terms in $\lt{\prec}{\H}$.  If $x_j\in\nmult{X}{L,\H,\prec}{\hv}$ for
some generator $\hv\in\H$, then, by definition of an involutive basis, $\H$
contains a unique generator $\bar\hv$ such that $\le{\prec}{\bar\hv}$ is an
involutive divisor of $\le{\prec}{(x_j\hv)}$.  In this case we include a
directed edge from $\le{\prec}{\hv}$ to $\le{\prec}{\bar\hv}$.  The thus
defined graph is called the \emph{$L$-graph} of the basis $\H$.

\begin{lemma}\label{lem:Lgraph}
  If the division\/ $L$ is continuous, then the\/ $L$-graph of any involutive
  set\/ $\H\subset\P$ is acyclic.
\end{lemma}

\begin{proof}
  The vertices of a path in an $L$-graph define a sequence as in the
  definition of a continuous division.  If the path is a cycle, then the
  sequence contains identical elements contradicting the continuity of the
  division.\qed
\end{proof}

We order the elements of $\H$ as follows: whenever the $L$-graph of $\H$
contains a path from $\hv_\alpha$ to $\hv_\beta$, then we must have
$\alpha<\beta$.  Any ordering satisfying this condition is called an
\emph{$L$-ordering}.  Note that by the lemma above for a continuous division
$L$-orderings always exist (although they are in general not unique).

For the Pommaret division $P$ it is easy to describe explicitly a $P$-ordering
without using the $P$-graph: we require that if either
$\cls{\hv_\alpha}<\cls{\hv_\beta}$ or $\cls{\hv_\alpha}=\cls{\hv_\beta}=k$ and
and the last non-vanishing entry of
$\le{\prec}{\hv_\alpha}-\le{\prec}{\hv_\beta}$ is negative, then we must have
$\alpha<\beta$.  Thus we sort the generators $\hv_\alpha$ first by their class
and within each class lexicographically (according to our definition in
Appendix A of Part I).  It is straightforward to verify that this defines
indeed a $P$-ordering.

\begin{example}\label{ex:Pgraph}
  Let us consider the ideal $\I\subset\kk[x,y,z]$ generated by the six
  polynomials $h_1=x^2$, $h_2=xy$, $h_3=xz-y$, $h_4=y^2$, $h_5=yz-y$ and
  $h_6=z^2-z+x$.  One easily verifies that they form a Pommaret basis $\H$ for
  the degree reverse lexicographic order.  The corresponding $P$-graph has the
  following form
  
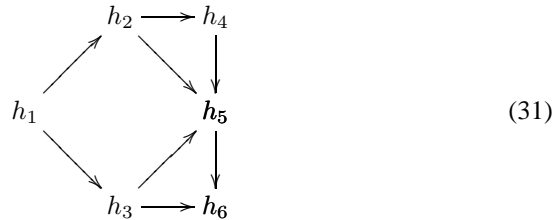
\begin{equation}
    \begin{xy}
      \xygraph{%
        !~:{@{->}}%
        []{h_1}
        (:[ur]{h_2}
         (:[r]{h_4}:[d]{h_5}:[d]{h_6},
          :[dr]{h_5}),
         :[dr]{h_3}
          (:[ur]{h_5},:[r]{h_6}))}
    \end{xy}
  \end{equation}  
  One clearly sees that the generators are already $P$-ordered, namely
  according to the description above.\bull
\end{example}

The decisive observation about an $L$-ordering is that we can now easily
determine the leading terms of all syzygies
$\S_{\alpha;k}\in\H_{\mathrm{Syz}}$ for the Schreyer order $\prec_{\H}$.

\begin{lemma}\label{lem:ltschreyer}
  Let the elements of the involutive basis\/ $\H\subset\P$ be ordered
  according to an\/ $L$-ordering.  Then the syzygies\/ $\S_{\alpha;k}$ satisfy
  $\lt{\prec_{\H}}{\S_{\alpha;k}}=x_k\ev_\alpha$.
\end{lemma}

\begin{proof}
  By the properties of the involutive standard representation, we have in
  (\ref{eq:invsyz}) $\lt{\prec}{(P^{(\alpha;k)}_\beta\hv_\beta)}\preceq
  \lt{\prec}{(x_k\hv_\alpha)}$ for all $\beta$ and only one index $\bar\beta$
  exists for which $\lt{\prec}{(P^{(\alpha;k)}_{\bar\beta}\hv_{\bar\beta})}=
  \lt{\prec}{(x_k\hv_\alpha)}$.  Thus $\le{\prec}{\hv_{\bar\beta}}$ is an
  involutive divisor of $\le{\prec}{(x_k\hv_{\alpha})}$ and the $L$-graph of
  $\H$ contains an edge from $\hv_\alpha$ to $\hv_{\bar\beta}$.  In an
  $L$-ordering, this implies $\alpha<\bar\beta$.  Now the assertion follows
  immediately from the definition of the term order $\prec_{\H}$.\qed
\end{proof}

There remains the problem of controlling the multiplicative variables
associated to these leading terms by the involutive division $L$.  For
arbitrary divisions it does not seem possible to make any statement.  Thus we
simply define a class of involutive divisions with the desired properties and
show afterwards that at least the Janet and the Pommaret division belong to
this class.

\begin{definition}\label{def:schreyerdiv}
  An involutive division\/ $L$ is of\/ \emph{Schreyer type} for the term
  order\/ $\prec$, if for any set\/ $\H$ which is involutive with respect to\/
  $L$ and\/ $\prec$ all sets\/ $\nmult{X}{L,\H,\prec}{\hv}$ with\/ $\hv\in\H$
  are again involutive.
\end{definition}

\begin{lemma}\label{lem:janpomschreyer}
  Both the Janet and the Pommaret division are of Schreyer type for any term
  order\/ $\prec$.
\end{lemma}

\begin{proof}
  For the Janet division any set of variables, i.\,e.\ monomials of degree
  one, is involutive.  Indeed, let $\F$ be such a set and $x_k\in\F$, then
  \begin{equation}
    \mult{X}{J,\F}{x_k}=\{x_i\mid x_i\notin\F\vee i\leq k\}
  \end{equation}
  which immediately implies the assertion.  For the Pommaret division sets of
  non-multiplicative variables are always of the form
  $\F=\{x_k,x_{k+1},\dots,x_n\}$ and such a set is trivially involutive.\qed
\end{proof}

An example of an involutive division which is not of Schreyer type is the
\emph{Thomas division} $T$ defined as follows: let $\N\subset\Nn$ be a finite
set and $\nu\in\N$ an arbitrary element; then $i\in\mult{N}{T,\N}{\nu}$, if
and only if $\nu_i=\max_{\mu\in\N}{\mu_i}$ (obviously, one may consider the
Janet division as a kind of refinement of the Thomas division).  One easily
sees that no set consisting only of variables can be involutive for the Thomas
division so that it cannot be of Schreyer type.

\begin{theorem}\label{thm:pomschreyer}
  Let\/ $L$ be a continuous involutive division of Schreyer type for the term
  order\/ $\prec$ and\/ $\H$ an\/ $L$-ordered involutive basis of the
  polynomial module\/ $\M$ with respect to\/ $L$ and\/ $\prec$.  Then
  $\H_{\mathrm{Syz}}$ is an involutive basis of\/ $\syz{}{\H}$ with respect
  to\/ $L$ and the term order\/ $\prec_{\H}$.
\end{theorem}

\begin{proof}
  By Lemma \ref{lem:ltschreyer}, the leading term of
  $\Sv_{\alpha;k}\in\H_{\mathrm{Syz}}$ is $x_k\ev_\alpha$ and we have one such
  generator for each non-multiplicative variable
  $x_k\in\nmult{X}{L,\H,\prec}{\hv_\alpha}$.  Since we assume that $L$ is of
  Schreyer type for $\prec$, these leading terms form an involutive set.  As
  we know already from Theorem~\ref{thm:invschreyer} that $\H_{\mathrm{Syz}}$
  is a Gr\"obner basis of $\syz{}{\H}$, the assertion follows trivially.\qed
\end{proof}

Note that under the made assumptions it follows immediately from the simple
form of the leading terms that $\H_{\mathrm{Syz}}$ is a minimal Gr\"obner
basis of $\syz{}{\H}$.

\begin{example}\label{ex:syzres}
  We continue with Example \ref{ex:Pgraph}.  As all assumption of Theorem
  \ref{thm:pomschreyer} are satisfied, the eight syzygies
  \begin{subequations}\label{eq:syz1}
  \begin{align}
    \Sv_{1;3}&=z\ev_1-x\ev_3-\ev_2\;,\\
    \Sv_{2;3}&=z\ev_2-x\ev_5-\ev_2\;,\\
    \Sv_{3;3}&=z\ev_3-x\ev_6+\ev_5-\ev_3+\ev_1\;,\\
    \Sv_{4;3}&=z\ev_4-y\ev_5-\ev_4\;,\\
    \Sv_{5;3}&=z\ev_5-y\ev_6+\ev_2\;,\\
    \Sv_{1;2}&=y\ev_1-x\ev_2\;,\\
    \Sv_{2;2}&=y\ev_2-x\ev_4\;,\\
    \Sv_{3;2}&=y\ev_3-x\ev_5+\ev_4-\ev_2
  \end{align}
  \end{subequations}
  form a Pommaret basis of the syzygy module $\syz{}{\H}$ with respect to the
  induced term order $\prec_{\H}$.  Indeed, as
  \begin{subequations}\label{eq:syz2}
  \begin{align}
    z\Sv_{1;2}&=y\Sv_{1;3}-x\Sv_{2;3}+x\Sv_{4;2}+\Sv_{2;2}\;,\\
    z\Sv_{2;2}&=y\Sv_{2;3}-x\Sv_{4;3}+\Sv_{2;2}\;,\\
    z\Sv_{3;2}&=y\Sv_{3;3}-x\Sv_{5;3}-\Sv_{2;3}+\Sv_{4;3}+
                \Sv_{3;2}-\Sv_{1;2}\;,
  \end{align}
  \end{subequations}
  all products of the generators with their non-multiplicative variables
  possess an involutive standard representation.\bull
\end{example}

\section{Free Resolutions I: The Polynomial Case}\label{sec:polyres}

As Theorem \ref{thm:pomschreyer} yields again an involutive basis, we may
apply it repeatedly and construct this way a syzygy resolution for any
polynomial module $\M$ given an involutive basis of it for an involutive
division of Schreyer type.  We specialise now to Pommaret bases where one can
even make a number of statements about the size of the resolution.  In
particular, we immediately obtain a stronger form of Hilbert's Syzygy Theorem
as a corollary (in fact, we will see later that this is the strongest possible
form, as the arising free resolution is always of minimal length).

\begin{theorem}\label{thm:hilsyzpom}
  Let\/ $\H$ be a Pommaret basis of the polynomial module\/ $\M\subseteq\P^m$.
  If we denote by\/ $\bq{0}{k}$ the number of generators\/ $\hv\in\H$ such
  that\/ $\cls{\le{\prec}{\hv}}=k$ and set\/ $d=\min\,\{k\mid\bq{0}{k}>0\}$,
  then\/ $\M$ possesses a finite free resolution
  \begin{equation}\label{eq:pomres}
    0\longrightarrow\P^{r_{n-d}}\longrightarrow\cdots\longrightarrow
    \P^{r_1}\longrightarrow\P^{r_0}\longrightarrow\M\longrightarrow0
  \end{equation}
  of length\/ $n-d$ where the ranks of the free modules are given by
  \begin{equation}\label{eq:resrank}
    r_i=\sum_{k=1}^{n-i}\binom{n-k}{i}\bq{0}{k}\;.
  \end{equation}
\end{theorem}

\begin{proof}
  According to Theorem~\ref{thm:pomschreyer}, $\H_{\mathrm{Syz}}$ is a
  Pommaret basis of $\syz{}{\H}$ for the term order $\prec_{\H}$.  Applying
  the theorem again, we can construct a Pommaret basis of the second syzygy
  module $\syz{2}{\H}$ and so on.  In the proof of
  Theorem~\ref{thm:pomschreyer} we showed that
  $\le{\prec_{\H}}{\Sv_{\alpha;k}}=x_k\ev_\alpha$.  Hence
  $\cls{\Sv_{\alpha;k}}=k>\cls{\hv_\alpha}$ and if $d$ is the minimal class of
  a generator in $\H$, then the minimal class in $\H_{\mathrm{Syz}}$ is $d+1$.
  This yields the length of the resolution (\ref{eq:pomres}), as a Pommaret
  basis with $d=n$ generates a free module.
  
  The ranks of the modules follow from a rather straightforward combinatorial
  calculation.  Let $\bq{i}{k}$ denote the number of generators of class $k$
  of the $i$th syzygy module~$\syz{i}{\H}$.  By definition of the generators
  $\Sv_{\alpha;k}$, we find $\bq{i}{k}=\sum_{j=1}^{k-1}\bq{i-1}{j}$, as each
  generator of class less than $k$ in the Pommaret basis of $\syz{i-1}{\H}$
  contributes one generator of class $k$ to the basis of $\syz{i}{\H}$.  A
  simple induction allows us to express the $\bq{i}{k}$ in terms of the
  $\bq{0}{k}$:
  \begin{equation}\label{eq:bik}
    \bq{i}{k}=\sum_{j=1}^{k-i}\binom{k-j-1}{i-1}\bq{0}{j}\;.
  \end{equation}
  The ranks of the modules in (\ref{eq:pomres}) are given by
  $r_i=\sum_{k=1}^n\bq{i}{k}$; entering (\ref{eq:bik}) yields via a classical
  identity for binomial coefficients (\ref{eq:resrank}).\qed
\end{proof}

\begin{remark}\label{rem:hilsyzgen}
  Theorem \ref{thm:hilsyzpom} remains valid for any involutive basis $\H$ with
  respect to a continuous division of Schreyer type, if we define $\bq{0}{k}$
  (respectively $\bq{i}{k}$ in the proof) as the number of generators with $k$
  multiplicative variables, since Theorem~\ref{thm:pomschreyer} holds for any
  such basis.  Indeed, after the first step we always analyse monomial sets of
  the form $\{x_{i_1},x_{i_2},\dots,x_{i_{n-k}}\}$ with
  $i_1<i_2<\cdots<i_{n-k}$.  By assumption, these sets are involutive and this
  is only possible, if one for the generators possesses $n$ multiplicative
  variables, another one $n-1$ and so on until the last generator which has
  only $n-k$ multiplicative variables (this follows for example from
  Proposition \ref{prop:hilbert} on the form of the Hilbert series).  Hence
  the basic recursion relation $\bq{i}{k}=\sum_{j=1}^{k-1}\bq{i-1}{j}$ and all
  subsequent combinatorial computations remain valid for any division of
  Schreyer type.
  
  For the special case of the Janet division, Plesken and Robertz
  \cite{pr:janet} proved directly the corresponding statement.  Here it is
  straightforward to determine explicitly the multiplicative variables for any
  syzygy: if $h_\alpha$ is a generator in the Janet basis $\H$ with the
  non-multiplicative variables
  $\nmult{X}{J,\H,\prec}{h_\alpha}=\{x_{i_1},x_{i_2},\dots,x_{i_{n-k}}\}$
  where $i_1<i_2<\cdots<i_{n-k}$, then
  \begin{equation}
    \mult{X}{J,\H_{\mathrm{Syz}},\prec}{\Sv_{\alpha;i_j}}=
        \{x_1,\dots,x_n\}\setminus
        \{x_{i_{j+1}},x_{i_{j+2}},\dots,x_{i_{n-k}}\}\;,
  \end{equation}
  as one easily verifies.\bull
\end{remark}

As in general the resolution (\ref{eq:pomres}) is not minimal, the ranks $r_i$
appearing in it cannot be identified with the Betti numbers of the module
$\M$.  However, they obviously represent an upper bound.  With a little bit
more effort one can easily derive similar bounds even for the multigraded
Betti numbers; we leave this as an exercise for the reader.

We may explicitly write the syzygy resolution (\ref{eq:pomres}) as a complex.
Let $\W$ be a free $\P$-module with basis $\{w_1,\dots,w_p\}$, i.\,e.\ its
rank is given by the size of the Pommaret basis $\H$.  Let $\V$ be a further
free $\P$-module with basis $\{v_1,\dots,v_n\}$, i.\,e.\ its rank is
determined by the number of variables in $\P$, and denote by $\Lambda\V$ the
exterior algebra over $\V$.  We set $\C_i=\W\otimes_{\P}\Lambda^i\V$ for
$0\leq i\leq n$.  If $\kv=(k_1,\dots,k_i)$ is a sequence of integers with
$1\leq k_1<k_2<\cdots<k_i\leq n$ and $v_{\kv}$ denotes the wedge product
$v_{k_1}\wedge\cdots\wedge v_{k_i}$, then a basis of this free $\P$-module is
given by the set of all tensor products $w_\alpha\otimes v_{\kv}$.  Finally,
we introduce the submodule $\S_i\subset\C_i$ generated by all those basis
elements where $k_1>\cls\hv_\alpha$.  Note that the rank of $\S_i$ is
precisely $r_i$ as defined by (\ref{eq:resrank}).

We denote the elements of the Pommaret basis of $\syz{i}{\H}$ by
$\Sv_{\alpha;\kv}$ with the inequalities $\cls{\hv_\alpha}<k_1<\cdots<k_i$.
An involutive normal form computation determines for every non-multiplicative
index $n\geq k_{i+1}>k_i=\cls{\Sv_{\alpha;\kv}}$ unique polynomials
$P_{\beta;\ellv}^{(\alpha;\kv,k_{i+1})}\in \kk[x_1,\dots,x_{\ell_i}]$ such
that
\begin{equation}\label{eq:nmsyz}
  x_{k_{i+1}}\Sv_{\alpha;\kv}=
      \sum_{\beta=1}^p\sum_{\ellv} P_{\beta;\ellv}^{(\alpha;\kv,k_{i+1})}
             \Sv_{\beta;\ellv}
\end{equation}
where the second sum is over all integer sequences
$\ellv=(\ell_1,\dots,\ell_i)$ satisfying
$\cls{\hv_\beta}<\ell_1<\cdots<\ell_i\leq n$.  Now we define the $\P$-module
homomorphisms $\epsilon:\S_0\rightarrow\M$ and
$\delta:\S_{i+1}\rightarrow\S_i$ by $\epsilon(w_\alpha)=\hv_\alpha$ and
\begin{equation}\label{eq:pomdiff}
  \delta(w_\alpha\otimes v_{\kv,k_{i+1}})=
     x_{k_{i+1}} w_\alpha\otimes v_{\kv}-
     \sum_{\beta,\ellv} P_{\beta;\ellv}^{(\alpha;\kv,k_{i+1})}
            w_\beta\otimes v_{\ellv}\;.     
\end{equation}
We extend the differential $\delta$ to a map $\C_{i+1}\rightarrow\C_i$ as
follows.  If $k_i\leq\cls{\hv_\alpha}$, then we set $\delta(w_\alpha\otimes
v_{\kv})=0$.  Otherwise let $j$ be the smallest value such that
$k_j>\cls{\hv_\alpha}$ and set (by slight abuse of notation)
\begin{equation}\label{eq:pomdiff2}
  \delta(w_\alpha\otimes v_{k_1}\wedge\cdots\wedge v_{k_i})=
      v_{k_1}\wedge\cdots\wedge v_{k_{j-1}}\wedge
      \delta(w_\alpha\otimes v_{k_j}\wedge\cdots\wedge v_{k_i})\;.
\end{equation}
Thus the factor $v_{k_1}\wedge\cdots\wedge v_{k_{j-1}}$ remains simply
unchanged and does not affect the differential.  This definition makes, by
construction, $(\C_*,\delta)$ to a complex and $(\S_*,\delta)$ to an exact
subcomplex which (augmented by the map $\epsilon:\S_0\rightarrow\M$) is
isomorphic to the syzygy resolution (\ref{eq:pomres}).

\begin{example}\label{ex:syzrescont}
  We continue with the ideal of Example \ref{ex:Pgraph} and \ref{ex:syzres},
  respectively.  As here $d=1$, the resolution has length $2$ in this case.
  Using the notation introduced above, the module $\S_0$ is then generated by
  $\{w_1,\dots,w_6\}$, the module $\S_1$ by the eight elements $\{w_1\otimes
  v_3,\dots, w_5\otimes v_3, w_1\otimes v_2,\dots, w_3\otimes v_2\}$ (the
  first three generators in the Pommaret basis $\H$ are of class $1$, the next
  two of class $2$ and the final one of class $3$) and the module $\S_2$ by
  $\{w_1\otimes v_2\wedge v_3,\dots, w_3\otimes v_2\wedge v_3\}$ corresponding
  to the three first syzygies of class $2$.  It follows from the expressions
  (\ref{eq:syz1}) and (\ref{eq:syz2}), respectively, for the first and second
  syzygies that the differential $\delta$ is here defined by the relations
  \begin{subequations}
    \begin{align}
      \delta(w_1\otimes v_3)&=zw_1-xw_3-w_2\;,\\
      \delta(w_2\otimes v_3)&=zw_2-xw_5-w_2\;,\\
      \delta(w_3\otimes v_3)&=zw_3-xw_6+w_5-w_3+w_1\;,\\
      \delta(w_4\otimes v_3)&=zw_4-yw_5-w_4\;,\\
      \delta(w_5\otimes v_3)&=zw_5-yw_6+w_2\;,\\
      \delta(w_3\otimes v_2)&=yw_3-xw_5+w_4-w_2\;,\\
      \delta(w_2\otimes v_2)&=yw_2-xw_4\;,\\
      \delta(w_1\otimes v_2)&=yw_1-xw_2\;,\\
      \delta(w_1\otimes v_2\wedge v_3)&=
          \begin{aligned}[t]
             &zw_1\otimes v_2-yw_1\otimes v_3+xw_2\otimes v_3-{}\\
             &xw_3\otimes v_2-w_2\otimes v_2\;.
          \end{aligned}\\
      \delta(w_2\otimes v_2\wedge v_3)&=zw_2\otimes v_2-yw_2\otimes v_3+
           xw_4\otimes v_3-w_2\otimes v_2\;,\\
      \delta(w_3\otimes v_2\wedge v_3)&=
          \begin{aligned}[t]
            &zw_3\otimes v_2-yw_3\otimes v_3+xw_5\otimes v_3+{}\\
            &w_2\otimes v_3-w_4\otimes v_3-w_3\otimes v_2+w_1\otimes v_2\;,
          \end{aligned}
    \end{align}
  \end{subequations}
  It is straightforward to verify explicitly the exactness of the complex
  $(\S_*,\delta)$.\bull
\end{example}

In the case that $m=1$ and thus $\M$ is actually an ideal in $\P$, it is
tempting to try to equip the complex $(\C_*,\delta)$ with the structure of a
differential algebra.  We first introduce a multiplication $\times$ on $\W$.
If $h_\alpha$ and $h_\beta$ are two elements of the Pommaret basis $\H$, then
their product possesses a unique involutive standard representation $h_\alpha
h_\beta=\sum_{\gamma=1}^p P_{\alpha\beta\gamma} h_\gamma$ and we define
\begin{equation}\label{eq:wmult}
  w_\alpha\times w_\beta=\sum_{\gamma=1}^p P_{\alpha\beta\gamma} w_\gamma
\end{equation}
and continue $\P$-linearly on $\W$.  This multiplication can be extended to
the whole complex $\C_*$ by defining for arbitrary elements $w,\bar w\in\W$
and $\omega,\bar\omega\in\Lambda\V$
\begin{equation}\label{eq:cmult}
  (w\otimes\omega)\times(\bar w\otimes\bar\omega)=
      (w\times\bar w)\otimes(\omega\wedge\bar\omega)\;.
\end{equation}

The distributivity of $\times$ is obvious from its definition.  For obtaining
a differential algebra, the product $\times$ must furthermore be associative
and satisfy the graded Leibniz rule $\delta(a\times b)=\delta(a)\times
b+(-1)^{\norm{a}}a\times\delta(b)$ where $\norm{a}$ denotes the form degree of
$a$.  While in general both conditions are not met, a number of special
situations exist where one indeed obtains a differential algebra.

Let us first consider the associativity.  It suffices to study it at the level
of $\W$ where we find that
\begin{subequations}
\begin{align}
  w_\alpha\times(w_\beta\times w_\gamma)&= \sum_{\delta,\epsilon=1}^p
  P_{\beta\gamma\delta} P_{\alpha\delta\epsilon} w_\epsilon\;,\\
  (w_\alpha\times w_\beta)\times w_\gamma&= \sum_{\delta,\epsilon=1}^p
  P_{\alpha\beta\delta} P_{\gamma\delta\epsilon} w_\epsilon\;.
\end{align}
\end{subequations}
One easily checks that both $\sum_{\delta,\epsilon=1}^p P_{\beta\gamma\delta}
P_{\alpha\delta\epsilon} h_\epsilon$ and $\sum_{\delta,\epsilon=1}^p
P_{\alpha\beta\delta} P_{\gamma\delta\epsilon} h_\epsilon$ are standard
representations of the product $h_\alpha h_\beta h_\gamma$ for the Pommaret
basis $\H$.  However, we cannot conclude that they are involutive standard
representations, as we do not know whether $P_{\beta\gamma\delta}$ and
$P_{\alpha\beta\delta}$, respectively, are multiplicative for $h_\epsilon$.
If this was the case, the associativity would follow immediately from the
uniqueness of involutive standard representations.

For the graded Leibniz rule the situation is similar but more involved.  In
the next section we will discuss it in more details for the monomial case.  In
the end, it boils down to analysing standard representations for products of
the form $x_kh_\alpha h_\beta$.  Again there exist two different ways for
obtaining them and a sufficient condition for the satisfaction of the Leibniz
rule is that both lead always to the unique involutive standard
representation.

\begin{example}\label{ex:diffalg1}
  Let us analyse the by now familiar ideal $\I\subset\kk[x,y,z]$ generated by
  $h_1=y^2-z$, $h_2=yz-x$ and $h_3=z^2-xy$.  We showed already in Part~I
  (Example~5.10) that these polynomials form a Pommaret basis of $\I$ for the
  degree reverse lexicographic term order.  The Pommaret basis of the first
  syzygy module consists of $\Sv_{1;3}=z\ev_1-y\ev_2+\ev_3$ and
  $\Sv_{2;3}=z\ev_2-y\ev_3-x\ev_1$.  As both generators are of class $3$, this
  is a free module and the resolution stops here.
  
  In a straightforward calculation one obtains for the multiplication $\times$
  the following defining relations:
  \begin{subequations}
    \begin{gather}
      w_1^2=w_3-yw_2+y^2w_1\;,\quad w_1\times w_2=-yw_3+y^2w_2-xw_1\;,\\
      w_1\times w_3=(y^2-z)w_3\;,\quad w_2^2=y^2w_3-xw_2+xyw_1\;,\\
      w_2\times w_3=(yz-x)w_3\;,\quad w_3^2=(z^2-xy)w_3\;.
    \end{gather}
  \end{subequations}
  Note that all coefficients of $w_1$ and $w_2$ are contained in $\kk[x,y]$
  and are thus multiplicative for all generators.  This immediately implies
  that our multiplication is associative, as any way to evaluate the product
  $w_\alpha\times w_\beta\times w_\gamma$ leads to the unique involutive
  standard representation of $h_\alpha h_\beta h_\gamma$.
  
  As furthermore in the only two non-multiplicative products $zh_2=yh_3+xh_1$
  and $zh_1=yh_2+h_3$ all coefficients on the right hand sides lie in
  $\kk[x,y]$, too, it follows from the same line of reasoning that the
  differential satisfies the Leibniz rule and we have a differential
  algebra.\bull
\end{example}

The situation is not always as favourable as in this example.  The next
example shows that in general we cannot expect to obtain a differential
algebra (in fact, not even an associative algebra).

\begin{example}\label{ex:diffalg2}
  Let us continue with the ideal of Examples \ref{ex:Pgraph}, \ref{ex:syzres}
  and \ref{ex:syzrescont}.  Evaluating the defining relation (\ref{eq:wmult})
  is particularly simple for the products of the form $w_i\times w_6=h_iw_6$,
  as all variables are multiplicative for the generator $h_6$.  Two further
  products are $w_5^2=y^2w_6-yw_5-xw_4$ and $w_3\times w_5=xyw_6-yw_5-xw_2$.
  In a straightforward computation one finds
  \begin{equation}
    (w_3\times w_5)\times w_5-w_3\times w_5^2=x^2w_4-xyw_2\;,
  \end{equation}
  so that the multiplication is not associative.  Note that the difference
  corresponds to the syzygy $x^2h_4-xyh_2=0$.  This is not surprising, as it
  encodes the difference between two standard representations of $h_3h_5^2$.
  The reason for the non-associativity lies in the coefficient $y$ of $w_5$ in
  the power $w_5^2$; it is non-multiplicative for $h_2$ and the generator
  $w_2$ appears in the product $w_3\times w_5$.  Hence computing $w_3\times
  w_5^2$ does not lead to an involutive standard representation of $h_3h_5^2$
  whereas the product $(w_3\times w_5)\times w_5$ does.\bull
\end{example}

\section{Free Resolutions II: The Monomial Case}\label{sec:monres}

In the special case of monomial modules, stronger results can be obtained.  In
particular, it is possible to obtain a closed form of the differential
(\ref{eq:pomdiff}) based only on the set $\H$ and to characterise those
modules for which (\ref{eq:pomres}) is a minimal resolution.  These are
generalisations of results by Eliahou and Kervaire \cite{ek:res}.

For a monomial module the existence of a Pommaret basis is a non-trivial
assumption, as the property of being a monomial module is not invariant under
coordinate transformations.  Therefore we always assume in the sequel that we
are dealing with a quasi-stable submodule $\M\subseteq\P^{m}$.  Let
$\H=\{\hv_1,\dots,\hv_p\}$ with $\hv_\alpha\in\TT^m$ be its monomial Pommaret
basis (by Proposition 2.11 of Part I, it is unique).  Furthermore, we
introduce the function $\Delta(\alpha,k)$ determining the unique generator in
the Pommaret basis $\H$ such that
$x_k\hv_\alpha=t_{\alpha,k}\hv_{\Delta(\alpha,k)}$ with a term
$t_{\alpha,k}\in\kk[\mult{X}{P}{\hv_{\Delta(\alpha,k)}}]$.

\begin{lemma}\label{lem:delta}
  The function\/ $\Delta$ and the terms\/ $t_{\alpha,k}$ satisfy the following
  relations.
  \begin{description}
  \item[\phantom{i}{\upshape (i)}] The inequality\/
    $\cls{\hv_\alpha}\leq\cls{\hv_{\Delta(\alpha,k)}}\leq k$ holds for all
    non-multiplicative indices\/ $k>\cls{\hv_\alpha}$.
  \item[{\upshape (ii)}] Let\/ $k_2>k_1> \cls{\hv_\alpha}$ be two
    non-multiplicative indices.  If\/ $\cls{\hv_{\Delta(\alpha,k_2)}}\geq
    k_1$, then\/ $\Delta\bigl(\Delta(\alpha,k_1),k_2\bigr)=\Delta(\alpha,k_2)$
    and\/ $x_{k_1}t_{\alpha,k_2}=t_{\alpha,k_1}t_{\Delta(\alpha,k_1),k_2}$.
    Otherwise we have the two equations\/
    $\Delta\bigl(\Delta(\alpha,k_1),k_2\bigr)=
    \Delta\bigl(\Delta(\alpha,k_2),k_1\bigr)$ and\/
    $t_{\alpha,k_1}t_{\Delta(\alpha,k_1),k_2}=
    t_{\alpha,k_2}t_{\Delta(\alpha,k_2),k_1}$.
  \end{description}
\end{lemma}

\begin{proof}
  Part (i) is trivial.  The inequality
  $\cls{\hv_\alpha}\leq\cls{\hv_{\Delta(\alpha,k)}}$ follows from the
  definition of $\Delta$ and the Pommaret division.  If
  $\cls{\hv_{\Delta(\alpha,k)}}>k$, then $\hv_{\Delta(\alpha,k)}$ would be an
  involutive divisor of $\hv_\alpha$ which contradicts the fact that any
  involutive basis is involutively head autoreduced.
  
  For Part (ii) we compute the involutive standard representation of
  $x_{k_1}x_{k_2}\hv_\alpha$.  There are two ways to do it.  We may either
  write
  \begin{equation}
    x_{k_1}x_{k_2}\hv_\alpha=x_{k_2}t_{\alpha,k_1}\hv_{\Delta(\alpha,k_1)}=
      t_{\alpha,k_1}t_{\Delta(\alpha,k_1),k_2}
      \hv_{\Delta(\Delta(\alpha,k_1),k_2)}\;,
  \end{equation}
  which is an involutive standard representation by Part (i), or start with
  \begin{equation}
    x_{k_1}x_{k_2}\hv_\alpha=x_{k_1}t_{\alpha,k_2}\hv_{\Delta(\alpha,k_2)}
  \end{equation}
  which requires a case distinction.  If $\cls{\hv_{\Delta(\alpha,k_2)}}\geq
  k_1$, we have already an involutive standard representation and its
  uniqueness implies our claim.  Otherwise we must rewrite multiplicatively
  $x_{k_1}\hv_{\Delta(\alpha,k_2)}=t_{\Delta(\alpha,k_2),k_1}
  \hv_{\Delta(\Delta(\alpha,k_2),k_1)}$ in order to obtain the involutive
  standard representation.  Again our assertion follows from its
  uniqueness.\qed
\end{proof}

Using this lemma, we can now provide a closed form for the differential
$\delta$ which does not require involutive normal form computations in the
syzygy modules $\syz{i}{\H}$ (which are of course expensive to perform) but is
solely based on information already computed during the determination of $\H$.
For its proof we must introduce some additional notations and conventions.  If
again $\kv=(k_1,\dots,k_i)$ is an integer sequence with $1\leq
k_1<\cdots<k_i\leq n$, then we write $\kv_j$ for the same sequence of indices
but with $k_j$ eliminated.  Its first entry is denoted by $(\kv_j)_1$; hence
$(\kv_j)_1=k_1$ for $j>1$ and $(\kv_j)_1=k_2$ for $j=1$.  The syzygy
$\Sv_{\alpha;\kv}$ is only defined for $\cls{\hv_\alpha}<k_1$.  We extend this
notation by setting $\Sv_{\alpha;\kv}=0$ for $\cls{\hv_\alpha}\geq k_1$.  This
convention will simplify some sums in the sequel.

\begin{theorem}\label{thm:pomdiffmon}
  Let\/ $\M\subseteq\P^m$ be a quasi-stable submodule and
  $\kv=(k_1,\dots,k_i)$.  Then the differential\/ $\delta$ of the complex\/
  $\C_*$ may be written in the form
  \begin{equation}\label{eq:pomdiffmon}
    \delta(w_\alpha\otimes v_{\kv})=
      \sum_{j=1}^i(-1)^{i-j}
          \bigl(x_{k_j}w_\alpha-t_{\alpha,k}w_{\Delta(\alpha,k)}\bigr)
          \otimes v_{\kv_j}\;.
  \end{equation}
\end{theorem}

\begin{proof}
  Note that all summands where $k_j$ is multiplicative for $\hv_\alpha$
  vanish.  This implies trivially (\ref{eq:pomdiff2}), so that we can restrict
  to the case that $\cls{\hv_\alpha}<k_1$.  Then our theorem is equivalent to
  \begin{equation}\label{eq:syzmon}
    \Sv_{\alpha;\kv}=
        \sum_{j=1}^i(-1)^{i-j}\bigl(x_{k_j}\Sv_{\alpha;\kv_j}-
                    t_{\alpha,k_j}\Sv_{\Delta(\alpha,k);\kv_j}\bigr)\;.
  \end{equation}
  Some of the terms $\Sv_{\Delta(\alpha,k);\kv_j}$ might vanish by our above
  introduced convention.  The equation (\ref{eq:syzmon}) is trivial for $i=1$
  (with $\Sv_\alpha=\hv_\alpha$) and a simple corollary to Lemma
  \ref{lem:delta} (ii) for $i=2$.
  
  For $i>2$ things become messy.  We proceed by induction on $i$.  In our
  approach, the syzygy $\Sv_{\alpha;\kv}$ arises from the non-multiplicative
  product $x_{k_i}\Sv_{\alpha;\kv_i}$.  Thus we must compute now the
  involutive normal form of this product.  By our induction hypothesis we may
  write
  \begin{equation}\label{eq:hypsyz}
    x_{k_i}\Sv_{\alpha;\kv_i}=
        \sum_{j=1}^{i-1}(-1)^{i-1-j}\bigl(x_{k_j}x_{k_i}\Sv_{\alpha;\kv_{ji}}-
        x_{k_i}t_{\alpha,k_j}\Sv_{\Delta(\alpha,k_j);\kv_{ji}}\bigr)\;.
  \end{equation}
  As $x_{k_i}$ is always non-multiplicative, using again the induction
  hypothesis, each summand may be replaced by the corresponding syzygy---but
  only at the expense of the introduction of many additional terms.  The main
  task in the proof will be to show that most of them cancel.  However, the
  cancellations occur in a rather complicated manner with several cases, so
  that no simple way for proving (\ref{eq:syzmon}) seems to exist.  We obtain
  the following lengthy expression:
  \begin{equation}\label{eq:exsyz}
    \begin{aligned}
      x_{k_i}&\Sv_{\alpha;\kv_i}=
         \sum_{j=1}^{i-1}(-1)^{i-1-j}\Bigl[x_{k_j}\Sv_{\alpha;\kv_j}\cnum{1}-
          t_{\alpha,k_j}\Sv_{\Delta(\alpha,k_j);\kv_j}\cnum{2}\Bigr]\\
      {}+{}&\sum_{j=1}^{i-1}x_{k_j}\Biggl[\,
          \sum_{\ell=1}^{j-1}(-1)^{\ell+j+1}
              x_{k_\ell}\Sv_{\alpha;\kv_{\ell j}}\cnum{3}
         -\sum_{\ell=j+1}^{i-1}(-1)^{\ell+j+1}
              x_{k_\ell}\Sv_{\alpha;\kv_{j\ell}}\cnum{4}\Biggr]\\
      {}-{}&\sum_{j=1}^{i-1}
            \begin{aligned}[t]
            \Biggl[
            &\sum_{\ell=1}^{j-1}(-1)^{\ell+j+1}x_{k_j}t_{\alpha,k_\ell}
                \Sv_{\Delta(\alpha,k_\ell);\kv_{\ell j}}\cnum{5}-{}\\
            &\sum_{\ell=j+1}^{i-1}(-1)^{\ell+j+1}x_{k_j}t_{\alpha,k_\ell}
                \Sv_{\Delta(\alpha,k_\ell);\kv_{j\ell}}\cnum{6}\Biggr]
            \end{aligned}\\
      {}+{}&\sum_{j=1}^{i-2}(-1)^{i-1-j}x_{k_j}t_{\alpha,k_i}
                \Sv_{\Delta(\alpha,k_i);\kv_{ji}}\cnum{7}+
            x_{k_{i-1}}t_{\alpha,k_i}
                \Sv_{\Delta(\alpha,k_i);\kv_{i-1,i}}\cnum{8}\\
      {}-{}&\sum_{j=1}^{i-1}t_{\alpha,k_j}
            \begin{aligned}[t]
              \Biggl[
              &\sum_{\ell=1}^{j-1}(-1)^{\ell+j+1}x_{k_\ell}
                   \Sv_{\Delta(\alpha,k_j);\kv_{\ell j}}\cnum{9}-{}\\
              &\sum_{\ell=j+1}^{i-1}(-1)^{\ell+j+1}x_{k_\ell}
                   \Sv_{\Delta(\alpha,k_j);\kv_{j\ell}}\cnum{10}\Biggr]
            \end{aligned}\\
      {}+{}&\sum_{j=1}^{i-1}t_{\alpha,k_j}
            \begin{aligned}[t]
              \Biggl[
              &\sum_{\ell=1}^{j-1}(-1)^{\ell+j+1}t_{\Delta(\alpha,k_j),k_\ell}
                   \Sv_{\Delta(\Delta(\alpha,k_j),k_\ell);\kv_{\ell j}}
                   \cnum{11}-{}\\
              &\sum_{\ell=j+1}^{i-1}(-1)^{\ell+j+1}
                   t_{\Delta(\alpha,k_j),k_\ell}
                   \Sv_{\Delta(\Delta(\alpha,k_j),k_\ell);\kv_{j\ell}}
                   \cnum{12}\Biggr]
            \end{aligned}\\
      {}-{}&\sum_{j=1}^{i-1}(-1)^{i-1-j}
                t_{\alpha,k_j}t_{\Delta(\alpha,k_j),k_i}
                \Sv_{\Delta(\Delta(\alpha,k_j),k_i);\kv_{ji}}
                \cnum{13}\;.
    \end{aligned}
  \end{equation}
  Note that the terms $\cnumt{7}$, $\cnumt{8}$ and $\cnumt{13}$, respectively,
  correspond to the special case $\ell=i$ (and $j=i-1$) in the sums
  $\cnumt{6}$ and $\cnumt{12}$, respectively.  We list them separately, as
  they must be treated differently.  The existence of any summand where the
  coefficient contains a term $t_{\cdot,\cdot}$ is bound on conditions.
  
  With the exception of the coefficient $x_{k_{i-1}}$ in the term $\cnumt{8}$,
  all coefficients are already multiplicative.  Thus this term must be further
  expanded using the induction hypothesis for the last time:
  \begin{equation}
    \begin{aligned}
      x_{k_{i-1}}t_{\alpha,k_i}&\Sv_{\Delta(\alpha,k_i);\kv_{i-1,i}}=
      t_{\alpha,k_i}\Sv_{\Delta(\alpha,k_i);\kv_i}\cnum{14}\\
      {}-{}&\sum_{j=1}^{i-2}(-1)^{i-1-j}x_{k_j}t_{\alpha,k_i}
                \Sv_{\Delta(\alpha,k_i);\kv_{ji}}\cnum{15}\\
      {}+{}&\sum_{j=1}^{i-1}(-1)^{i-1-j}t_{\alpha,k_i}
                t_{\Delta(\alpha,k_i),k_j}
                \Sv_{\Delta(\Delta(\alpha,k_i),k_j);\kv_{ji}}\cnum{16}\;.
    \end{aligned}
  \end{equation}
  
  The left hand side of (\ref{eq:exsyz}) and the terms $\cnumt{1}$,
  $\cnumt{2}$ and $\cnumt{14}$ represent the syzygy $\Sv_{\alpha,\kv}$ we are
  looking for.  We must thus show that all remaining terms vanish.  In order
  to simplify the discussion of the double sums, we swap $j$ and $\ell$ in
  $\cnumt{3}$, $\cnumt{5}$, $\cnumt{9}$ and $\cnumt{11}$ so that everywhere
  $j<\ell$.  It is now easy to see that $\cnumt{3}$ and $\cnumt{4}$ cancel;
  each summand of $\cnumt{3}$ also appears in $\cnumt{4}$ but with the
  opposite sign.  Note, however, that the same argument does not apply to
  $\cnumt{11}$ and $\cnumt{12}$, as the existence of these terms is bound to
  different conditions!
  
  For the other cancellations, we must distinguish several cases depending on
  the classes of the generators in the Pommaret basis $\H$.  We first study
  the double sums and thus assume that $1\leq j<i$.
  \begin{itemize}
  \item If $\cls{\hv_{\Delta(\alpha,k_j)}}<(\kv_{j})_1$, the terms $\cnumt{5}$
    and $\cnumt{10}$ are both present and cancel each other.  We must now make
    a second case distinction on the basis of $\hv_{\Delta(\alpha,k_\ell)}$.
    \begin{itemize}
    \item If $\cls{\hv_{\Delta(\alpha,k_\ell)}}<(\kv_{j})_1$, then the terms
      $\cnumt{6}$ and $\cnumt{9}$ are also present and cancel each other.
      Furthermore, both $\cnumt{11}$ and $\cnumt{12}$ exist and cancel due to
      the second case of Lemma \ref{lem:delta} (ii).
    \item If $\cls{\hv_{\Delta(\alpha,k_\ell)}}\geq(\kv_{j})_1$, then none of
      the four terms $\cnumt{6}$, $\cnumt{9}$, $\cnumt{11}$ and $\cnumt{12}$
      exists.  For the latter two terms, this fact is a consequence of the
      first case of Lemma \ref{lem:delta} (ii).
    \end{itemize}
  \item If $\cls{\hv_{\Delta(\alpha,k_j)}}\geq(\kv_{j})_1$, then neither
    $\cnumt{5}$ nor $\cnumt{10}$ nor $\cnumt{12}$ exists.  For the remaining
    double sums, we must again consider the class of
    $\hv_{\Delta(\alpha,k_\ell)}$.
    \begin{itemize}
    \item If $\cls{\hv_{\Delta(\alpha,k_\ell)}}<(\kv_{j})_1$, then the terms
      $\cnumt{6}$ and $\cnumt{9}$ exist and cancel each other.  The term
      $\cnumt{11}$ does not exist, as Lemma \ref{lem:delta} implies the
      inequalities $\cls{\hv_{\Delta(\Delta(\alpha,k_\ell),k_j)}}=
      \cls{\hv_{\Delta(\Delta(\alpha,k_j),k_\ell)}}\geq
      \cls{\hv_{\Delta(\alpha,k_j)}} \geq(\kv_{j})_1$.
    \item If $\cls{\hv_{\Delta(\alpha,k_\ell)}}\geq(\kv_{j})_1$, then neither
      $\cnumt{6}$ nor $\cnumt{9}$ exist and the term $\cnumt{11}$ is not
      present either; this time the application of Lemma \ref{lem:delta} (ii)
      yields the chain of inequalities
      $\cls{\hv_{\Delta(\Delta(\alpha,k_\ell),k_j)}}\geq
      \cls{\hv_{\Delta(\alpha,k_\ell)}} \geq(\kv_{j})_1$.
    \end{itemize}
  \end{itemize}
  
  For the remaining terms everything depends on the class of
  $\hv_{\Delta(\alpha,k_i)}$ controlling in particular the existence of the
  term $\cnumt{8}$.
  \begin{itemize}
  \item If $\cls{\hv_{\Delta(\alpha,k_i)}}<k_1\leq(\kv_j)_1$, then the term
    $\cnumt{8}$ exists and generates the terms $\cnumt{15}$ and $\cnumt{16}$.
    Under this condition, the term $\cnumt{7}$ is present, too, and because of
    Lemma \ref{lem:delta} (ii) it cancels $\cnumt{15}$.  Again by Lemma
    \ref{lem:delta} (ii), the conditions for the existence of $\cnumt{13}$ and
    $\cnumt{16}$ are identical and they cancel each other.
  \item If $\cls{\hv_{\Delta(\alpha,k_i)}}\geq k_1$, then $\cnumt{8}$ and
    consequently $\cnumt{15}$ and $\cnumt{16}$ are not present.  The analysis
    of $\cnumt{7}$ and $\cnumt{13}$ requires a further case distinction.
    \begin{itemize}
    \item Under the made assumption, the case
      $\cls{\hv_{\Delta(\alpha,k_i)}}<(\kv_j)_1$ can occur only for $j=1$ as
      otherwise $(\kv_j)_1=k_1$.  Because of Lemma \ref{lem:delta} (ii), the
      terms $\cnumt{7}$ and $\cnumt{13}$ exist for $j=1$ and cancel each
      other.
    \item If $\cls{\hv_{\Delta(\alpha,k_i)}}\geq(\kv_j)_1$, then $\cnumt{7}$
      does not exist.  The term $\cnumt{13}$ is also not present, but there
      are two different possibilities: depending on which case of Lemma
      \ref{lem:delta} (ii) applies, we either find
      $\cls{\hv_{\Delta(\Delta(\alpha,k_j),k_i)}}=
      \cls{\hv_{\Delta(\alpha,k_i)}}$ or
      $\cls{\hv_{\Delta(\Delta(\alpha,k_j),k_i)}}=
      \cls{\hv_{\Delta(\Delta(\alpha,k_i),k_j)}}\geq
      \cls{\hv_{\Delta(\alpha,k_i)}}$; but in any case the class is too high.
    \end{itemize}
  \end{itemize}
  
  Thus we have shown that indeed all terms vanish with the exception of
  $\cnumt{1}$, $\cnumt{2}$ and $\cnumt{14}$ which are needed for the syzygy
  $\Sv_{\alpha,\kv}$.  This proves our claim.\qed
\end{proof}

As in the previous section, we may introduce for monomial ideals, i.\,e.\ for
$m=1$, the product $\times$.  The right hand side of its defining equation
(\ref{eq:wmult}) simplifies for a monomial basis $\H$ to
\begin{equation}\label{eq:wmultmon}
  w_\alpha\times w_\beta=m_{\alpha,\beta}w_{\Gamma(\alpha,\beta)}
\end{equation}
where the function $\Gamma(\alpha,\beta)$ determines the unique generator
$h_{\Gamma(\alpha,\beta)}$ such that $h_\alpha
h_\beta=m_{\alpha,\beta}h_{\Gamma(\alpha,\beta)}$ with a term
$m_{\alpha,\beta}\in\kk[\mult{X}{P}{h_{\Gamma(\alpha,\beta)}}]$.
Corresponding to Lemma \ref{lem:delta}, we obtain now the following result.

\begin{lemma}\label{lem:gamma}
  The function\/ $\Gamma$ and the terms\/ $m_{\alpha,\beta}$ satisfy the
  following relations.
  \begin{description}
  \item[\phantom{ii}{\upshape (i)}]
    $\cls{h_{\Gamma(\alpha,\beta)}}\geq\max{\{\cls{h_\alpha},\cls{h_\beta}\}}$.
  \item[\phantom{i}{\upshape (ii)}]
    $\Gamma\bigl(\Gamma(\alpha,\beta),\gamma\bigr)=
    \Gamma\bigl(\alpha,\Gamma(\beta,\gamma)\bigr)$ and\/
    $m_{\alpha,\beta}m_{\Gamma(\alpha,\beta),\gamma}=
    m_{\beta,\gamma}m_{\Gamma(\beta,\gamma),\alpha}$.
  \item[{\upshape (iii)}] $\Gamma\bigl(\Delta(\alpha,k),\beta\bigr)=
    \Delta\bigl(\Gamma(\alpha,\beta),k\bigr)$ and\/
    $t_{\alpha,k}m_{\Delta(\alpha,k),\beta}=
    t_{\Gamma(\alpha,\beta),k}m_{\alpha,\beta}$.
  \end{description}
\end{lemma}

\begin{proof}
  Part (i) is obvious from the definition of the function $\Gamma$.  Part (ii)
  and (iii), respectively, follow from the analysis of the two different ways
  to compute the involutive standard representation of $h_\alpha h_\beta
  h_\gamma$ and $x_k h_\alpha h_\beta$, respectively.  We omit the details, as
  they are completely analogous to the proof of Lemma \ref{lem:delta}.\qed
\end{proof}

\begin{theorem}
  Let\/ $\H$ be the Pommaret basis of the quasi-stable ideal\/
  $\I\subseteq\P$.  Then the product\/ $\times$ defined by (\ref{eq:wmultmon})
  makes the complex\/ $(\C_*,\delta)$ to a differential algebra.
\end{theorem}

\begin{proof}
  This is a straightforward consequence of Lemma \ref{lem:gamma}.  Writing out
  the relations one has to check, one easily finds that Part (ii) ensures the
  associativity of $\times$ and Part (iii) the satisfaction of the graded
  Leibniz rule.\qed
\end{proof}

\section{Minimal Resolutions and Projective Dimension}\label{sec:minres}

Recall that for a graded polynomial module $\M$ a graded free resolution is
\emph{minimal}, if all entries of the matrices corresponding to the maps
$\phi_i:\P^{r_i}\rightarrow\P^{r_{i-1}}$ are of positive degree, i.\,e.\ no
constant coefficients appear.  Up to isomorphisms, the minimal resolution is
unique and its length is an important invariant, the \emph{projective
dimension} $\pd{\M}$ of the module.  If the module $\M$ is graded, the
resolution (\ref{eq:pomres}) is obviously graded, too.  However, in general,
it is not minimal.  We introduce now a class of monomial modules for which it
is always the minimal resolution.

\begin{definition}\label{def:stable}
  A (possibly infinite) set\/ $\N\subseteq\Nn$ is called\/ \emph{stable}, if
  for each multi index\/ $\nu\in\N$ all multi indices\/ $\nu-1_k+1_j$ with\/
  $k=\cls{\nu}<j\leq n$ are also contained in\/ $\N$.  A monomial submodule\/
  $\M\subseteq\P^m$ is\/ \emph{stable}, if each of the sets\/
  $\N_\alpha=\bigl\{\mu\mid x^\mu\ev_\alpha\in\M\bigr\}\subseteq\Nn$ with\/
  $1\leq\alpha\leq m$ is stable.
\end{definition}

\begin{remark}\label{rem:borel}
  The stable modules are of considerable interest, as they contain as a subset
  the \emph{Borel-fixed} modules, i.\,e.\ modules $\M\subseteq\P^m$ which
  remain invariant under the natural action of the Borel
  group.\footnote{Classically, the Borel group consists of upper triangular
    matrices.  In our ``inverse'' conventions we must take lower triangular
    matrices.}  Indeed, one can show that (for a ground field of
  characteristic $0$) a module is Borel-fixed, if and only if it can be
  generated by a set $\S$ of monomials such that whenever $x^\nu\ev_j\in\S$
  then also $x^{\nu-1_k+1_j}\ev_j\in\S$ for all $\cls{\nu}\leq k<j\leq n$
  \cite[Thm.~15.23]{de:ca}.  Generically, the leading terms of any polynomial
  module form a Borel-fixed module \cite[Thm.~15.20]{de:ca}.  Note that while
  stability is obviously independent of the characteristic of the base field,
  the same is not true for the notion of a Borel-fixed module.\bull
\end{remark}

Any monomial submodule has a unique minimal basis.  For stable submodules it
must coincide with its Pommaret basis.  This result represents a very simple
and effective characterisation of stable submodules.  Furthermore, it shows
that any stable submodule is trivially quasi-stable and thus explains the
terminology introduced in Definition \ref{def:qstab}.

\begin{proposition}\label{prop:invstab}
  Let\/ $\M\subseteq\P^m$ be a monomial submodule.  $\M$ is stable, if and
  only if its minimal basis\/ $\H$ is simultaneously a Pommaret basis.
\end{proposition}

\begin{proof}
  Let us assume first that $\M$ is stable; we have to show that
  $\ispan{\H}{P}=\M$.  For every term $\sv\in\M$ a unique term $\tv_1\in\H$
  exists such that $\sv=s_1\tv_1$ for some term $s_1\in\TT$.  If
  $s_1\in\kk[X_P(\tv_1)]$, we are done.  Otherwise there exists an index
  $j>k=\cls{\tv_1}$ such that $x_j\mid s_1$ and we rewrite
  \begin{equation}
    \sv=\left(\frac{x_k}{x_j}s_1\right)\left(\frac{x_j}{x_k}\tv_1\right)\;.
  \end{equation}
  Since, by assumption, $\M$ is stable, $(x_j/x_k)\tv_1\in\M$.  Thus a term
  $\tv_2\in\H$ exists such that $(x_j/x_k)\tv_1=s_2\tv_2$ for some term
  $s_2\in\TT$.  We are done, if $\tv_2\idiv{P}\sv$.  Otherwise we iterate this
  construction.  By the continuity of the Pommaret division, this cannot go on
  infinitely, i.\,e.\ after a finite number of steps we must reach a term
  $\tv_N\in\H$ such that $\tv_N\idiv{P}\sv$ and thus $\sv\in\ispan{\H}{P}$.
  
  For the converse, assume that the minimal basis $\H$ is a Pommaret basis.
  Then a unique $\tv\in\H$ exists for each $\sv\in\M$ such that
  $\tv\idiv{P}\sv$.  We must show that with $k=\cls{\sv}\leq\cls{\tv}$ for all
  $i>k$ the terms $(x_i/x_k)\sv$ are also elements of $\M$.  We distinguish
  two cases.  If $\sv=\tv$, a $\bar{\tv}\in\H$ exists with
  $\bar{\tv}\idiv{P}(x_i\tv)$.  As $\H$ is a minimal basis, it cannot be that
  $\bar{\tv}=x_i\tv$.  Instead we must have that
  $\bar{\tv}\mid(x_i/x_k)\tv=(x_i/x_k)\sv$ and we are done.  If $\sv\neq\tv$,
  we write $\sv=s\tv$ with $s\in\TT$.  If $k<\cls{\tv}$, then $x_k\mid s$
  which implies that we can divide by $x_k$ and thus $(x_i/x_k)\sv\in\M$.
  Otherwise, $\cls{\sv}=\cls{\tv}$ and we know from the first case that
  $(x_i/x_k)\tv\in\M$.  But $(x_i/x_k)\sv=(x_i/x_k)s\tv$.\qed
\end{proof}

This result is due to Mall \cite{mall:pom} (but see also the remark after
Lemma~1.2 in \cite{ek:res}).  He used it to give a direct proof of the
following interesting result which for us is an immediate corollary to our
definition of an involutive basis.

\begin{theorem}\label{thm:redgbpom}
  Let\/ $\M\subseteq\P^m$ be a graded submodule in generic position and\/ $\G$
  the reduced Gr\"obner basis of\/ $\M$ for an arbitrary term order\/ $\prec$.
  If\/ $\chr{\kk}=0$, then\/ $\G$ is also the minimal Pommaret basis of\/
  $\I$ for\/ $\prec$.
\end{theorem}

\begin{proof}
  The generic initial module of any module is Borel fixed \cite{gall:weier}
  (see also the discussion in \cite[Section~15.9]{de:ca}) which implies in
  particular that it is stable for $\chr{\kk}=0$.  By definition of a reduced
  Gr\"obner basis, $\le{\prec}{\G}$ is the minimal basis of $\le{\prec}{\M}$.
  But as $\le{\prec}{\M}$ is generically stable, this implies that
  $\le{\prec}{\G}$ is a Pommaret basis of $\le{\prec}{\M}$ and thus $\G$ is a
  Pommaret basis of $\M$.  As the Pommaret division is global, $\G$ is
  automatically the unique minimal basis.\qed
\end{proof}

\begin{remark}\label{rem:stable}
  It follows from Lemma~\ref{lem:stable}, that if $\H$ is a Pommaret basis of
  the submodule $\M\subseteq\P^m$ of degree $q$, then $(\lt{\prec}{\M})_{\geq
    q}$ is a stable monomial submodule, as $\lt{\prec}{\H_q}$ is obviously its
  minimal basis and simultaneously a Pommaret basis.\bull
\end{remark}

\begin{theorem}\label{thm:minres}
  Let\/ $\M\subseteq\P^m$ be a quasi-stable module.  Then the syzygy resolution
  given by\/ (\ref{eq:pomres}) is minimal, if and only if\/ $\M$ is stable.
\end{theorem}

\begin{proof}
  We denote again the elements of the Pommaret basis of $\syz{i}{\H}$ as above
  by $\Sv_{\alpha;\kv}$ with the inequalities
  $\cls{\hv_\alpha}<k_1<\cdots<k_i$ and the basis elements of the free module
  $\P^{r_i}$ by $\ev_{\alpha;\kv}$ with the same inequalities.  For a monomial
  module our syzygies are, by construction, of the form
  \begin{equation}\label{eq:isyz}
    \Sv_{\alpha;\kv}=x_{k_i}\ev_{\alpha;\kv_i}-
        \sum_{\beta=1}^p\sum_{\ellv} t_{\beta;\ellv}^{(\alpha;\kv)}
             \ev_{\beta;\ellv}
  \end{equation}
  with terms $t_{\beta;\ellv}^{(\alpha;\kv)}\in \kk[x_1,\dots,x_{\ell_{i-1}}]$
  and where the second sum is over all sequences
  $\ellv=(\ell_1,\dots,\ell_{i-1})$ with
  $\cls{\hv_\beta}<\ell_1<\cdots<\ell_{i-1}$.  By definition, the resolution
  is minimal, if and only if all non-vanishing
  $t_{\beta;\ellv}^{(\alpha;\kv)}\notin\kk$.
  
  We proceed by an induction over the resolution.  If $\M$ is stable, the
  Pommaret basis $\H$ is the minimal basis of $\M$ by
  Proposition~\ref{prop:invstab}.  For each element $\hv_\alpha\in\H$ and each
  non-multiplicative index $k>\cls{\hv_\alpha}$ the product $x_k\hv_\alpha$
  has the involutive standard representation
  $t_{\alpha,k}\hv_{\Delta(\alpha,k)}$.  As $\H$ is the minimal basis of $\M$,
  it is not possible that $t_{\alpha,k}\in\kk$, as otherwise
  $\hv_\alpha\mid\hv_{\Delta(\alpha,k)}$.
  
  Assume that $t_{\beta;\ellv}^{(\alpha;\kv)}\notin\kk$ for all possible
  values of the scripts.  For the induction we must show that all
  non-vanishing coefficients $t_{\beta;\ellv,\ell_i}^{(\alpha;\kv,k_{i+1})}$
  with $\ell_i>\ell_{i-1}$ and $k_{i+1}>k_i$ do not lie in $\kk$.  For this
  purpose, we rewrite
  \begin{equation}\label{eq:i1syz}
      x_{k_{i+1}}\Sv_{\alpha;\kv}=
          x_{k_i}(x_{k_{i+1}}\ev_{\alpha;\kv_i})-
          \sum_{\beta,\ellv} t_{\beta;\ellv}^{(\alpha;\kv)}
              (x_{k_{i+1}}\ev_{\beta;\ellv})
  \end{equation}
  and compute a normal form with respect to the Pommaret basis of
  $\syz{i+1}{\H}$.  We solve each syzygy (\ref{eq:isyz}) for its initial term,
  the first term on the right hand side, and substitute the result in
  (\ref{eq:i1syz}) whenever there is a non-multiplicative product of a basis
  vector $\ev_{\beta;\ellv}$.  As all non-vanishing coefficients
  $t_{\beta;\ellv}^{(\alpha;\kv)}\notin\kk$, such a substitution can never
  decrease the degree of a coefficient.  Thus
  $t_{\beta;\ellv,\ell_i}^{(\alpha;\kv,k_{i+1})}\notin\kk$ for all possible
  values of the scripts, too.
  
  This proves that for a stable module the resolution (\ref{eq:pomres}) is
  minimal.  Assume for the converse that $\M$ is not stable.  By
  Proposition~\ref{prop:invstab} the Pommaret basis $\H$ cannot be the minimal
  basis.  In this case, $\H$ contains elements $\hv_\alpha$, $\hv_\beta$ with
  $\hv_\beta=x_k\hv_\alpha$ for some $x_k\in\nmult{X}{P}{\hv_\alpha}$.  So the
  syzygy $\Sv_{\alpha;k}=x_k\ev_\alpha-\ev_\beta$ contains a constant
  coefficient and the resolution is not minimal.\qed
\end{proof}

\begin{example}\label{ex:stableres}
  One might be tempted to conjecture that this result extended to polynomial
  modules, i.\,e.\ that (\ref{eq:pomres}) was minimal for polynomial modules
  $\M$ with stable leading module $\lt{\prec}{\M}$.  Unfortunately, this is
  not true.  Consider the homogeneous ideal $\I\subset\kk[x,y,z]$ generated by
  $h_1=z^2+xy$, $h_2=yz-xz$, $h_3=y^2+xz$, $h_4=x^2z$ and $h_5=x^2y$.  One
  easily checks that these elements form a Pommaret basis $\H$ for the degree
  reverse lexicographic term order and that $\lt{\prec}{\M}$ is a stable
  module.  A Pommaret basis of $\syz{}{\H}$ is given by
  \begin{subequations}
    \begin{align}
      \Sv_{2;3}&=z\ev_2+(x-y)\ev_1+x\ev_3-\ev_4-\ev_5\;,\\
      \Sv_{3;3}&=z\ev_3-x\ev_1-(x+y)\ev_2-\ev_4+\ev_5\;,\\
      \Sv_{4;3}&=z\ev_4-x^2\ev_1+x\ev_5\;,\\
      \Sv_{5;3}&=z\ev_5-x^2\ev_2-x\ev_4\;,\\
      \Sv_{4;2}&=(y-x)\ev_4-x^2\ev_2\;,\\
      \Sv_{5;2}&=y\ev_5-x^2\ev_3+x\ev_4\;.
    \end{align}
  \end{subequations}
  As the first two generators show, the resolution (\ref{eq:pomres}) is not
  minimal.\bull
\end{example}

Given an arbitrary graded free resolution, it is a standard task to reduce it
to the minimal resolution using just some linear algebra (see for example
\cite[Chapt.~6, Theorem~3.15]{clo:uag} for a detailed discussion).  Thus for
any concrete module $\M$ it is straightforward to obtain from
(\ref{eq:pomres}) the minimal resolution.  However, even in the monomial case
it seems highly non-trivial to find a closed form description of the outcome
of the minimisation process.  Nevertheless, the resolution (\ref{eq:pomres})
contains so much structure that certain statements are possible.

\begin{theorem}\label{thm:pd}
  Let\/ $\H$ be a Pommaret basis of the graded module\/ $\M\subseteq\P^m$ for
  a class respecting term order and set\/ $d=\min_{\hv\in\H}\cls{\hv}$.
  Then\/ $\pd{\M}=n-d$.
\end{theorem}

\begin{proof}
  Consider the resolution (\ref{eq:pomres}) which is of length $n-d$.  The
  last map in it is defined by the syzygies $\Sv_{\alpha;(d+1,\dots,n)}$
  originating in the generators $\hv_\alpha\in\H$ with $\cls{\hv_\alpha}=d$.
  Choose now among these generators an element $\hv_{\gamma}$ of maximal
  degree (recall that the same choice was crucial in the proof of Proposition
  \ref{prop:Iregular}).  Then the syzygy $\Sv_{\gamma;(d+1,\dots,n)}$ cannot
  contain any constant coefficient, as the coefficients of all basis vectors
  $\ev_{\beta;\kv}$ where the last entry of $\kv$ is $n$ must be contained in
  $\lspan{x_1,\dots,x_{n-1}}$ and the coefficients of the basis vectors
  $\ev_{\alpha;(d+1,\dots,n-1)}$ cannot be constant for degree reasons.
  
  If we start now a minimisation process at the end of the resolution, then it
  will never introduce a constant term into the syzygy
  $\Sv_{\gamma;(d+1,\dots,n)}$ and thus it will never be eliminated.  It is
  also not possible that it is reduced to zero, as the last map in a free
  resolution is obviously injective.  This implies that the last term of the
  resolution will not vanish during the minimisation and the length of the
  minimal resolution, i.\,e.\ $\pd{\M}$, is still $n-d$.\qed
\end{proof}

The graded form of the \emph{Auslander-Buchsbaum formula}
\cite[Exercise~19.8]{de:ca} is now a trivial corollary to this theorem and
Proposition \ref{prop:Iregular} on the depth.  Note that, in contrast to other
proofs, our approach is constructive in the sense that we automatically have
an explicit regular sequence of maximal length and an explicit free resolution
of minimal length.

\begin{corollary}[Auslander-Buchsbaum]
  Let\/ $\M\subseteq\P^m$ be a graded polynomial module with\/
  $\P=\kk[x_1,\dots,x_n]$.  Then\/ $\depth{\M}+\pd{\M}=n$.
\end{corollary}

As for a monomial module we do not need a term order, we obtain as further
simple corollary the following relation between $\pd{\M}$ and
$\pd{(\lt{\prec}{\M})}$.

\begin{corollary}\label{cor:pdlt}
  Let\/ $\M\subseteq\P^m$ be a graded module and\/ $\prec$ an arbitrary term
  order for which\/ $\M$ possesses a Pommaret basis.  Then\/
  $\pd{\M}\leq\pd{(\lt{\prec}{\M})}$.
\end{corollary}

\begin{proof}
  Let $\H$ be the Pommaret basis of the module $\M$ for the term order $\prec$
  and set $d=\min_{\hv\in\H}\cls{(\lt{\prec}{\hv})}$.  Then it follows
  immediately from Theorem \ref{thm:pd} that $\pd{(\lt{\prec}{\M})}=n-d$.  On
  the other hand, Theorem \ref{thm:hilsyzpom} guarantees the existence of the
  free resolution (\ref{eq:pomres}) of length $n-d$ for $\M$ so that this
  value is an upper bound for $\pd{\M}$.\qed
\end{proof}

\section{Castelnuovo-Mumford Regularity}\label{sec:cmr}

For notational simplicity we restrict again to ideals instead of submodules.
In many situations it is of interest to obtain a good estimate on the degree
of an ideal basis.  Up to date, no satisfying answer is known to this
question.  Somewhat surprisingly, the stronger problem of bounding not only
the degree of a basis of $\I$ but also of its syzygies can be treated
effectively.

\begin{definition}\label{def:cmreg}
  Let\/ $\I\subseteq\P$ be a homogeneous ideal.  $\I$ is called\/
  \emph{$q$-regular}, if its\/ $i$th syzygy module is generated by elements of
  degree less than or equal to\/ $q+i$.  The\/ \emph{Castelnuovo-Mumford
    regularity} $\reg{\I}$ is the least\/ $q$ for which\/ $\I$ is\/
  $q$-regular.
\end{definition}

Among other applications, the Castelnuovo-Mumford regularity $\reg{\I}$ is a
useful measure for the complexity of Gr\"obner basis computations
\cite{bm:cac}.  The question of effectively computing $\reg{\I}$ has recently
attracted some interest.  In this section we show that $\reg{\I}$ is trivially
determined by a Pommaret basis with respect to the degree reverse
lexicographic order and provide alternative proofs to some characterisations
of the Castelnuovo-Mumford regularity proposed in the literature.

\begin{theorem}\label{thm:pomreg}
  Let\/ $\I\subseteq\P$ be a homogeneous ideal.  The Castelnuovo-Mumford
  regularity of\/ $\I$ is\/ $q$, if and only if\/ $\I$ has in some coordinates
  a homogeneous Pommaret basis of degree\/ $q$ with respect to the degree
  reverse lexicographic order.
\end{theorem}

\begin{proof}
  Let $\xv$ be some $\delta$-regular coordinates for the ideal $\I$ so that it
  possess a Pommaret basis $\H$ of degree $q$ with respect to the degree
  reverse lexicographic order in these coordinates.  Then the $i$th module of
  the syzygy resolution (\ref{eq:pomres}) induced by the basis $\H$ is
  obviously generated by elements of degree less than or equal to $q+i$.  Thus
  we have the trivial estimate $\reg{\I}\leq q$ and there only remains to show
  that it is in fact an equality.
  
  For this purpose, consider a generator $h_\gamma\in\H$ of degree $q$ which
  is of minimal class among all elements of this maximal degree $q$ in $\H$.
  If $\cls{h_\gamma}=n$, then $h_\gamma$ cannot be removed from $\H$ without
  loosing the basis property, as the leading term of no other generator of
  class $n$ can divide $\lt{\prec}{h_\gamma}$ and, since the degree reverse
  lexicographic order is class respecting, all other generators do not contain
  any terms of class $n$.  Hence we trivially find $\reg{\I}=q$ in this case.
  
  If $\cls{h_\gamma}=n-i$ for some $i>0$, then the resolution
  (\ref{eq:pomres}) contains at the $i$th position the syzygy
  $\Sv_{\gamma;(n-i+1,\dots,n)}$ of degree $q+i$.  Assume now that we minimise
  the resolution step by step starting at the end.  We claim that the syzygy
  $\Sv_{\gamma;(n-i+1,\dots,n)}$ is not eliminated during this process.
  
  There are two possibilities how $\Sv_{\gamma;(n-i+1,\dots,n)}$ could be
  removed during the minimisation.  The first one is that a syzygy at the next
  level of the resolution contained the term $\ev_{\gamma;(n-i+1,\dots,n)}$
  with a constant coefficient.  Any such syzygy is of the form
  (\ref{eq:nmsyz}) with $\cls{h_\alpha}<n-i$ and
  $\cls{h_\alpha}<k_1<\cdots<k_i<n$ and its leading term is
  $x_{k_{i+1}}\ev_{\alpha;\kv}$ with $k_{i+1}>k_i$.  However, since
  $\cls{(x_{k_1}\cdots x_{k_{i+1}}h_\alpha)}<n-i$ and $\cls{(x_{n-i+1}\cdots
    x_nh_\gamma)}=n-i$, it follows from our use of the degree reverse
  lexicographic order (since we assume that everything is homogeneous, both
  polynomial have the same degree) and the definition of the induced Schreyer
  term orders, that the term $\ev_{\gamma;(n-i+1,\dots,n)}$ is greater than
  the leading term $x_{k_{i+1}}\ev_{\alpha;\kv}$ of any syzygy
  $\Sv_{\alpha;(k_1,\dots,k_{i+1})}$ at the level $i+1$ and thus cannot
  appear.
  
  The second possibility is that $\Sv_{\gamma;(n-i+1,\dots,n)}$ itself
  contained a constant coefficient at some vector $\ev_{\beta;\ellv}$.
  However, this required $\deg{h_\beta}=\deg{h_\gamma}+1$ which is again a
  contradiction.\footnote{For later use we note the following fact about this
    argument.  If $\ev_{\beta;\ellv}$ is a constant term in the syzygy
    $\Sv_{\gamma;(n-i+1,\dots,n)}$, then it must be smaller than the leading
    term and hence $\lt{\prec}{(x_{\ell_1}\cdots
      x_{\ell_i}h_\beta)}\prec\lt{\prec}{(x_{k_1}\cdots x_{k_{i+1}}h_\alpha)}$
    implying that $\cls{h_\beta}\leq\cls{h_\gamma}$.  Thus it suffices, if
    $h_\gamma$ is of maximal degree among all generators $h_\beta\in\H$ with
    $\cls{h_\beta}\leq\cls{h_\gamma}$.  For the special case that $h_\gamma$
    is of minimial class, we exploited this observation already in the proof
    of Theorem \ref{thm:pd}.}  As the minimisation process never introduces
  new constant coefficients, the syzygy $\Sv_{\gamma;(n-i+1,\dots,n)}$ may
  only be modified but not eliminated.  Furthermore, the modifications cannot
  make $\Sv_{\gamma;(n-i+1,\dots,n)}$ to the zero syzygy, as otherwise a basis
  vector of the next level was in the kernel of the differential.  However,
  this is not possible, as we assume that the tail of the resolution is
  already minimised and by the exactness of the sequence any kernel member
  must be a linear combination of syzygies.  Hence the final minimal
  resolution will contain at the $i$th position a generator of degree $q+i$
  and $\reg{\I}=q$.\qed
\end{proof}

To some extent this result was to be expected.  By Theorem~\ref{thm:redgbpom},
the reduced Gr\"obner basis is generically also a Pommaret basis and according
to Bayer and Stillman \cite{bs:mreg} this basis has for the degree reverse
lexicographic order generically the degree $\reg{\I}$.  Thus the only surprise
is that Theorem \ref{thm:pomreg} does not require that the leading ideal is
stable and the Pommaret basis $\H$ is not necessarily a reduced Gr\"obner
basis (if the ideal $\I$ has a Pommaret basis of degree $q$, then the
truncated ideal $(\le{\prec}{\I})_{\geq q}$ is always stable by
Remark~\ref{rem:stable} and thus the set $\H_q$ defined by (\ref{eq:truncbas})
is the reduced Gr\"obner basis of $\I_{\geq q}$).

Note furthermore that Theorem \ref{thm:pomreg} implies a remarkable fact:
given arbitrary coordinates~$\xv$, the ideal $\I$ either does not possess a
finite Pommaret basis for the degree reverse lexicographic order or, if such a
basis exists, it is of degree $\reg{\I}$.  Hence using Pommaret bases, it
becomes trivial to determine the Castelnuovo-Mumford regularity: it is just
the degree of the basis.

\begin{remark}\label{rem:cmrpos}
  The proof of Theorem \ref{thm:pomreg} also provides us with information
  about the positions at which in the minimal resolution the maximal degree is
  attained.  We only have to look for all elements of maximal degree in the
  Pommaret basis; their classes correspond to these positions.\bull
\end{remark}

\begin{remark}\label{rem:cmrjan}
  Recall from Remark \ref{rem:hilsyzgen} that Theorem \ref{thm:hilsyzpom}
  remains valid for any involutive basis $\H$ with respect to a continuous
  division of Schreyer type (with an obvious modification of the definition of
  the numbers $\bq{0}{k}$) and that it is independent of the used term order.
  It follows immediately from the form of the resolution (\ref{eq:pomres}),
  i.\,e.\ from the form of the maps in it given by the respective involutive
  bases according to Theorem~\ref{thm:pomschreyer}, that always the estimate
  $\reg{\I}\leq\deg{\H}$ holds and thus any such basis provides us with a
  bound for the Castelnuovo--Mumford regularity.
  
  This observation also implies that an involutive basis with respect to a
  division of Schreyer type and an arbitrary term order can never be of lower
  degree than the Pommaret basis for the degree reverse lexicographic order.
  The latter one is thus in this sense optimal.  As a concrete example
  consider again the ideal mentioned in Remark \ref{rem:depth}: in ``good''
  coordinates a Pommaret basis of degree $2$ exists for it and after a simple
  permutation of the variables its Janet basis is of degree $4$.\bull
\end{remark}

In analogy to Corollary \ref{cor:pdlt} comparing the projective dimension of a
module~$\M$ and its leading module $\lt{\prec}{\M}$ with respect to an
arbitrary term order $\prec$, we may derive a similar estimate for the
Castelnuovo-Mumford regularity.

\begin{corollary}\label{cor:reglt}
  Let\/ $\I\subseteq\P$ be a homogeneous ideal and\/ $\prec$ an arbitrary term
  order such that a Pommaret basis of\/ $\I$ exists.  Then\/
  $\reg{\I}\leq\reg{(\lt{\prec}{\I})}$.
\end{corollary}

\begin{proof}
  Let $\H$ be the Pommaret basis of $\I$ for the term order $\prec$ and set
  $q=\deg{\H}$.  Then it follows immediately from Theorem \ref{thm:pomreg}
  that $\reg{(\lt{\prec}{\I})}=q$.  On the other hand, the form of the free
  resolution (\ref{eq:pomres}) implies trivially that $\reg{\I}\leq q$.\qed
\end{proof}

\begin{remark}
  Bayer et al. \cite{bcp:betti} introduced a refinement of the
  Castelnuovo--Mumford regularity: the \emph{extremal Betti numbers}.  Recall
  that the (graded) Betti number $\beta_{ij}$ of the ideal $I$ is defined as
  the number of minimal generators of degree $i+j$ of the $i$th module in the
  minimal free resolution of $\I$ (thus $\reg{\I}$ is the maximal value $j$
  such that $\beta_{i,i+j}>0$ for some $i$).  A Betti number $\beta_{ij}>0$ is
  called extremal, if $\beta_{k\ell}=0$ for all $k\geq i$ and $\ell>j$.  There
  always exists a least one extremal Betti number: if we take the maximal
  value $i$ for which $\beta_{i,i+\reg{\I}}>0$, then $\beta_{i,i+\reg{\I}}$ is
  extremal.  In general, there may exist further extremal Betti numbers.
  Bayer et al. \cite[Thm.~1.6]{bcp:betti} proved that for any ideal $\I$ both
  the positions and the values of the extremal Betti numbers coincides with
  those of its generic initial ideal with respect to the degree reverse
  lexicographic order.
  
  Our proof of Theorem \ref{thm:pomreg} allows us to make the same statement
  for the ordinary initial ideal for $\prec_{\mbox{\scriptsize
      degrevlex}}$---provided the coordinates are $\delta$-regular.
  Furthermore, it shows that the extremal Betti numbers of $\I$ can be
  immediately read off the Pommaret basis $\H$ of $\I$.  Finally, if we
  introduce ``pseudo-Betti numbers'' for the (in general non-minimal)
  resolution (\ref{eq:pomres}), then the positions and values of the extremal
  numbers coincide with the true extremal Betti numbers of $\I$.
  
  Take the generator $h_\gamma$ used in the proof of Theorem \ref{thm:pomreg}.
  If $\cls{h_\gamma}=n-i_1$ and $\deg{h_\gamma}=q_1$, then the considerations
  in the proof imply immediately that $\beta_{i_1,q_1+i_1}$ is an extremal
  Betti number and its value is given by the number of generators of degree
  $q_1$ and class $n-i_1$ in the Pommaret basis $\H$.  If $i_1=\depth{\I}$,
  then this is the only extremal Betti number.  Otherwise, let $q_2$ be the
  maximal degree of a generator $h\in\H$ with $\cls{h}<n-i_1$ and assume that
  $n-i_2$ is the minimal class of such a generator.  Then the arguments used
  in the proof of Theorem \ref{thm:pomreg} show that $\beta_{i_2,q_2+i_2}$ is
  also an extremal Betti number and that its value is given by the number of
  generators of degree $q_2$ and class $n-i_2$ in the Pommaret basis $\H$.
  Continuing in this manner, we obtain all extremal Betti numbers.  Since all
  our considerations depend only on the leading terms of the generators, we
  find for the leading ideal $\lt{\prec}{\I}$ exactly the same situation.\bull
\end{remark}

Combining the above results with Remark \ref{rem:stable} and Proposition
\ref{prop:invstab} immediately implies the following generalisation of a
result by Eisenbud, Reeves and Totaro \cite[Prop.~10]{ert:veronese} for
Borel-fixed monomial ideals.

\begin{proposition}
  Let\/ $\I$ be a quasi-stable ideal generated in degrees less than or equal
  to\/ $q$.  The ideal\/ $\I$ is\/ $q$-regular, if and only if the
  truncation\/ $\I_{\geq q}$ is stable.
\end{proposition}

Bayer and Stillman \cite{bs:mreg} gave the following characterisation of
$q$-regularity for which we now provide a new proof.  Note the close
relationship between their first condition and the idea of assigning
multiplicative variables.

\begin{theorem}\label{thm:cmreg}
  Let\/ $\I\subseteq\P$ be a homogeneous ideal which can be generated
  by elements of degree less than or equal to\/ $q$.  Then\/ $\I$ is\/
  $q$-regular, if and only if for some value\/ $0\leq d\leq n$ linear
  forms\/ $y_1,\dots,y_d\in\P_1$ exist such that
  \begin{subequations}\label{eq:bscrit}
    \begin{gather}
      \bigl(\lspan{\I,y_1,\dots,y_{j-1}}:y_j\bigr)_q=
          \lspan{\I,y_1,\dots,y_{j-1}}_q\;,
          \quad 1\leq j\leq d\;,\label{eq:bscrit1}\\
      \lspan{\I,y_1,\dots,y_d}_q=\P_q\;.\label{eq:bscrit2}
    \end{gather}
  \end{subequations}
\end{theorem}

\begin{proof}
  Assume first that the conditions (\ref{eq:bscrit}) are satisfied for some
  linear forms $y_1,\dots,y_d\in\P_1$ and choose variables $\xv$ such that
  $x_i=y_i$ for $1\leq i\leq d$.  Let the finite set $\H_q$ be a basis of
  $\I_q$ as a vector space in triangular form with respect to the degree
  reverse lexicographic order, i.\,e.\ $\lt{\prec}{h_1}\neq\lt{\prec}{h_2}$
  for all $h_1,h_2\in\H_q$.  We claim that $\H_q$ is a Pommaret basis of the
  truncation $\I_{\geq q}$ implying that the full ideal $\I$ possesses a
  Pommaret basis of degree $q'\leq q$ and hence by Theorem \ref{thm:pomreg}
  that $\reg{\I}\leq q$.
  
  Let us write $\H_q=\bigl\{h_{k,\ell}\mid 1\leq k\leq n,\
  1\leq\ell\leq\ell_k\bigr\}$ where $\cls{h_{k,\ell}}=k$.  A basis of the
  vector space $\lspan{\I,x_1,\dots,x_j}_q$ is then given by all $h_{k,\ell}$
  with $k>j$ and all terms in $\lspan{x_1,\dots,x_j}_q$.  We will now show
  that
  \begin{equation}\label{eq:prolbas}
    \H_{q+1}=\bigl\{x_j h_{k,\ell}\mid 1\leq j\leq k,\ 1\leq k\leq n,\ 
                                       1\leq\ell\leq\ell_k\bigr\}
  \end{equation}
  is a basis of $\I_{q+1}$ as a vector space.  This implies that $\H_q$ is
  locally involutive for the Pommaret division and thus involutive by
  Corollary 7.3 of Part I.  Since, by assumption, $\I$ is generated in degrees
  less than or equal to $q$, we have furthermore $\lspan{\H_q}=\I_{\geq q}$ so
  that indeed $\H_q$ is a Pommaret basis of the ideal $\I_{\geq q}$.
  
  Let $f\in\I_{q+1}$ and $\cls{f}=j$.  By the properties of the degree reverse
  lexicographic order this implies that $f=x_j\hat f+g$ with $\hat
  f\in\bigl(\kk[x_j,\dots,x_n]\setminus\{0\}\bigr)_q$ and
  $g\in\bigl(\lspan{x_1,\dots,x_{j-1}}\bigr)_{q+1}$ (cf.\ Lemma~A.1 of Part
  I).  We distinguish two cases.  The condition (\ref{eq:bscrit2}) implies
  that $\bigl(\lspan{\I,x_1,\dots,x_d}\bigr)_q=\P_q$.  Thus if $j>d$, we may
  write $\hat f=\sum_{k=d+1}^n\sum_{\ell=1}^{\ell_k}c_{k,\ell}h_{k,\ell}+\hat
  g$ with $c_{k,\ell}\in\kk$ and $\hat
  g\in\bigl(\lspan{x_1,\dots,x_d}\bigr)_q$.  We set $\hat
  f_0=\sum_{k=j}^n\sum_{\ell=1}^{\ell_k}c_{k,\ell}h_{k,\ell}$ and $\hat
  f_1=\sum_{k=d+1}^{j-1}\sum_{\ell=1}^{\ell_k}c_{k,\ell}h_{k,\ell}+\hat g$.
  Obviously, $\hat f\in\bigl(\lspan{\I,x_1,\dots,x_{j-1}}:x_j\bigr)_q$.  If
  $j\leq d$, then the condition (\ref{eq:bscrit1}) implies that actually $\hat
  f\in\lspan{\I,x_1,\dots,x_{j-1}}_q$.  Hence in this case we may
  decompose $\hat f=\hat f_0+\hat f_1$ with $\hat
  f_0=\sum_{k=j}^n\sum_{\ell=1}^{\ell_k}c_{k,\ell}h_{k,\ell}$ and $\hat
  f_1\in\bigl(\lspan{x_1,\dots,x_{j-1}}\bigr)_q$.
  
  It is trivial that $\lspan{\H_{q+1}}\subseteq\I_{q+1}$ (here we mean the
  linear span over $\kk$ and not over $\P$).  We show by an induction over $j$
  that $\I_{q+1}\subseteq\lspan{\H_{q+1}}$.  If $j=1$, then $f=x_1\hat f$ with
  $\hat f\in\I_q$.  Thus $f\in\lspan{\H_{q+1}}$.  If $j>1$, we write
  $f=f_0+f_1$ with $f_0=x_j\hat f_0$ and $f_1=x_j\hat f_1+g$ where $\hat f_0$
  and $\hat f_1$ have been defined above.  By construction,
  $f_0\in\lspan{\H_{q+1}}$, as $x_j$ is multiplicative for all generators
  contained in $\hat f_0$, and $f_1\in\I_{q+1}$ with $\cls{f_1}<j$.  According
  to our inductive hypothesis this implies that $f_1\in\lspan{\H_{q+1}}$, too.
  Hence $\lspan{\H_{q+1}}=\I_{q+1}$.
  
  Assume conversely that the ideal $\I$ is $q$-regular.  Then, by Theorem
  \ref{thm:pomreg}, it possesses a Pommaret basis $\H$ of degree $\reg{\I}\leq
  q$ with respect to the degree reverse lexicographic order.  We set
  $d=\dim{\P/\I}$ and claim that for the choice $y_i=x_i$ for $1\leq i\leq d$
  the conditions (\ref{eq:bscrit}) are satisfied.  For the second equality
  (\ref{eq:bscrit2}) this follows immediately from Proposition \ref{prop:dim}
  which shows that it actually holds already at degree $\reg{\I}\leq q$.
  
  For the first equality (\ref{eq:bscrit1}) take a polynomial
  $f\in\bigl(\lspan{\I,x_1,\dots,x_{j-1}}:x_j\bigr)_q$.  By definition, we
  have then $x_jf\in\lspan{\I,x_1,\dots,x_{j-1}}$.  If
  $f\in\lspan{x_1,\dots,x_{j-1}}$, then there is nothing to prove.  Otherwise,
  a polynomial $g\in\lspan{x_1,\dots,x_{j-1}}$ exists such that $x_jf-g\in\I$
  and obviously $\cls{(x_jf-g)}=j$.  If we introduce the set $\H_{\geq
  j}=\{h\in\H\mid\cls{h}\geq j\}$, the involutive standard representation of
  $x_jf-g$ induces an equation $x_jf=\sum_{h\in\H_{\geq j}}P_hh+\bar g$ where
  $\bar g\in x_j\lspan{x_1,\dots,x_{j-1}}$ and $P_h\in\lspan{x_j}$ (this is
  trivial if $\cls{h}>j$ and follows from $\deg{h}\leq q$ if $\cls{h}=j$).
  Thus we can divide by $x_j$ and find that already
  $f\in\lspan{\I,x_1,\dots,x_{j-1}}_q$.\qed
\end{proof}

Bayer and Stillman \cite{bs:mreg} further proved that in \emph{generic}
coordinates it is not possible to find a Gr\"obner basis of degree less than
$\reg{\I}$ and that this estimate is sharp, as it is realised by bases with
respect to the degree reverse lexicographic order.  The restriction to the
generic case is here essential, as for instance most monomial ideals are
trivial counter examples.  Hence their result is only of limited use for the
actual computation of the Castelnuovo-Mumford regularity, as one never knows
whether one works with generic coordinates.

\begin{example}\label{ex:pbcomplex}
  Consider the homogeneous ideal
  \begin{equation}
    \I=\lspan{z^{8}-wxy^{6},\ y^{7}-x^{6}z,\  yz^{7}-wx^{7}}
       \subset\QQ[w,x,y,z]\;.
  \end{equation}
  The given basis of degree $8$ is already a Gr\"obner basis for the degree
  reverse lexicographic term order.  If we perform a permutation of the
  variables and consider $\I$ as an ideal in $\QQ[w,y,x,z]$, then we obtain
  for the degree reverse lexicographic term order the following Gr\"obner
  basis of degree $50$:
  \begin{multline}
    \bigl\{y^{7}-x^{6}z,\ yz^{7}-wx^{7},\ z^{8}-wxy^{6},\ 
           y^{8}z^{6}-wx^{13},\\
           y^{15}z^{5}-wx^{19},\ y^{22}z^{4}-wx^{25},\ y^{29}z^{3}-wx^{31},\\
           y^{36}z^{2}-wx^{37},\ y^{43}z-wx^{43},\ y^{50}-wx^{49}\bigr\}\;.
  \end{multline}
  Unfortunately, neither coordinate system is generic: as $\reg{\I}=13$, one
  yields a basis of too low degree and the other one one of too high degree.
  
  With a Pommaret basis it is no problem to determine the Castelnuovo-Mumford
  regularity, as the first coordinate system is $\delta$-regular.  A Pommaret
  basis of $\I$ for the degree reverse lexicographic term order is obtained by
  adding the polynomials $z^k(y^{7}-x^{6}z)$ for $1\leq k\leq 6$ and thus the
  degree of the basis is indeed $13$.\bull
\end{example}

Yet another characterisation of $q$-regularity is due to Eisenbud and Goto
\cite{eg:freeres}.  Again it can be obtained as an easy corollary to Theorem
\ref{thm:pomreg}.

\begin{theorem}\label{thm:linres}
  The homogeneous ideal\/ $\I\subseteq\P$ is\/ $q$-regular, if and only if its
  truncation\/ $\I_{\geq q}$ admits a linear free resolution, i.\,e.\ 
  $\I_{\geq q}$ is generated by elements of degree\/ $q$ and all maps in the
  resolution are linear in the sense that the entries of the matrices
  describing them are zero or homogeneous polynomials of degree $1$.
\end{theorem}

\begin{proof}
  If $\I$ is $q$-regular, then by Theorem \ref{thm:pomreg} it possesses in
  suitable coordinates a Pommaret basis $\H$ of degree $\reg{\I}\leq q$.  The
  set $\H_q$ defined by (\ref{eq:truncbas}) is a Pommaret basis of the
  truncated ideal $\I_{\geq q}$ according to Lemma~\ref{lem:truncbas}.  Now it
  follows easily from Theorem~\ref{thm:pomschreyer} that $\I_{\geq q}$
  possesses a linear free resolution, as all syzygies in the resolution
  (\ref{eq:pomres}) derived from $\H_q$ are necessarily homogeneous of degree
  $1$.
  
  The converse is trivial.  The existence of a linear resolution for $\I_{\geq
    q}$ immediately implies that $\reg{\I_{\geq q}}=q$.  Hence $\I_{\geq q}$
  possesses a Pommaret basis of degree $q$ by Theorem \ref{thm:pomreg}
  entailing the existence of a Pommaret basis for $\I$ of degree $q'\leq q$.
  Hence, again by Theorem \ref{thm:pomreg}, $\reg{\I}=q'\leq q$.\qed
\end{proof}

\section{Regularity and Saturation}\label{sec:regsat}

Already in the work of Bayer and Stillman \cite{bs:mreg} on the
Castelnuovo-Mumford regularity the \emph{saturation} $\Isat$ of a homogeneous
ideal $\I\subseteq\P$ plays an important role.  Recall that by definition
\begin{equation}\label{eq:Isat}
  \Isat=\I:\P_+^\infty=
      \bigl\{f\in\P\mid\exists k\in\NN_0:f\cdot\P_k\subset\I\bigr\}\;.
\end{equation}
An ideal $\I$ such that $\I=\Isat$ is called \emph{saturated}.  We show now
first how $\Isat$ can be effectively determined from a Pommaret basis of
$\I$.\footnote{It seems to be folklore that for Gr\"obner bases the
  construction in Proposition \ref{prop:Isat} yields a Gr\"obner basis of
  $\I:x_1^\infty$; in \cite[Prop.~5.1.11]{ms:pc} this observation is
  attributed (without reference) to Bayer.  In our case we do not only get a
  Pommaret basis but it also turns out that here $\Isat=\I:x_1^\infty$ (see
  the remarks below).}

\begin{proposition}\label{prop:Isat}
  Let\/ $\H$ be a Pommaret basis of the homogeneous ideal\/ $\I$ for a class
  respecting term order.  We introduce the sets\/
  $\H_1=\{h\in\H\mid\cls{h}=1\}$ and
  $\bar\H_1=\bigl\{h/x_1^{\deg_{x_1}{\lt{\prec}{h}}}\mid h\in\H_1\bigr\}$.
  Then\/ $\bar\H=(\H\setminus\H_1)\cup\bar\H_1$ is a weak Pommaret basis of
  the saturation\/ $\Isat$.
\end{proposition}

\begin{proof}
  Recall that for terms of the same degree any class respecting term order
  coincides with the reverse lexicographic order.  Hence of all terms in a
  generator $h\in\H_1$ the leading term $\lt{\prec}{h}$ has the lowest
  $x_1$-degree.  This implies in particular that $\bar\H_1$ is well-defined
  and does not contain a generator of class $1$ anymore.
  
  We first show that indeed $\bar\H_1\subset\Isat$.  Let
  $d_1=\max_{h\in\H_1}{\{\deg_{x_1}{\lt{\prec}{h}}\}}$ and
  $\Delta=d_1+\max_{h\in\H_1}{\{\deg{h}\}}-\min_{h\in\H_1}{\{\deg{h}\}}$.  We
  claim that $\P_\Delta\cdot\bar h\subset\I$ for all $\bar h\in\bar\H_1$.
  Thus let $x^\mu\in\P_\Delta$ and choose $k\in\NN_0$ such that $x_1^k\bar
  h\in\H_1$; obviously, we have $k\leq d_1$.  Since the polynomial $x^\mu
  x_1^k\bar h$ lies in $\I$, it possesses an involutive standard
  representation of the form
  \begin{equation}\label{eq:hbarisr}
    x^\mu x_1^k\bar h=\sum_{h\in\H\setminus\H_1}P_hh+
                       \sum_{h\in\H_1}Q_hh
  \end{equation}
  with $P_h\in\kk[x_1,\dots,x_{\cls{h}}]$ and $Q_h\in\kk[x_1]$.  
  
  The left hand side of this equation is contained in $\lspan{x_1^k}$ and thus
  also the right hand side.  Analysing an involutive normal form computation
  leading to the representation (\ref{eq:hbarisr}), one immediately sees that
  this implies that all coefficients $P_h$ (since here $\cls{h}>1$) and all
  summands $Q_hh$ lie in $\lspan{x_1^k}$.  As a first consequence of this
  representation we observe that for any monomial $x^\mu$ (not necessarily of
  degree $\Delta$) we may divide (\ref{eq:hbarisr}) by $x_1^k$ and then obtain
  an involutive standard representation of $x^\mu\bar h$ with respect to the
  set $\bar\H$; hence this set is indeed weakly involutive for the Pommaret
  division and the given term order.
  
  If $x^\mu\in\P_\Delta$, then we find for any $h\in\H_1$ that $|\deg{\bar
    h}-\deg{h}|\leq\Delta$ and hence $\deg{Q_h}=\deg{\bigl(x^\mu x_1^k\bar
    h\bigr)}-\deg{h}\geq k$.  Since $Q_h\in\kk[x_1]$, this implies that under
  the made assumption on $x^\mu$ already the coefficient $Q_h$ lies in
  $\lspan{x_1^k}$ so that the product $x^\mu\bar h$ possesses an involutive
  standard representation with respect to $\H$ and thus is contained in the
  ideal $\I$ as claimed.
  
  Now we show that every polynomial $f\in\Isat$ may be decomposed into an
  element of $\I$ and a linear combination of elements of $\bar\H_1$.  We may
  write $f=\tilde f+g$ where $\tilde f$ is the involutive normal form of $f$
  with respect to $\H$ and $g\in\I$.  If $\tilde f=0$, then already $f\in\I$
  and nothing is to be shown.  Hence we assume that $\tilde f\neq0$.  By
  definition of the saturation $\Isat$, there exists a $k\in\NN_0$ such that
  $\tilde f\cdot\P_k\subset\I$, hence in particular $x_1^k\tilde f\in\I$.
  This implies that $\lt{\prec}{(x_1^k\tilde f)}\in\ispan{\lt{\prec}{\H}}{P}$.
  Therefore a unique generator $h\in\H$ exists with
  $\lt{\prec}{h}\idiv{P}\lt{\prec}{(x_1^k\tilde f)}$.
  
  So let $\lt{\prec}{(x_1^k\tilde f)}=x^\mu\lt{\prec}{h}$ and assume first
  that $\cls{h}>1$.  Since the term on the left hand side is contained in
  $\lspan{x_1^k}$, we must have $\mu_1\geq k$ so that we can divide by
  $x_1^k$.  But this implies that already $\lt{\prec}{\tilde
    f}\in\ispan{\lt{\prec}{\H}}{P}$ contradicting our assumption that $\tilde
  f$ is in involutive normal form.  Hence we must have $\cls{h}=1$ and by the
  same argument as above $\mu_1<k$.
  
  Division by $x_1^k$ shows that $\lt{\prec}{\tilde
    f}\in\ispan{\lt{\prec}{\bar\H_1}}{P}$.  Performing the corresponding
  involutive reduction leads to a new element $f_1\in\Isat$.  We compute again
  its involutive normal form $\tilde f_1$ and apply the same argument as
  above, if $\tilde f_1\neq0$.  After a finite number of such reductions we
  obtain an involutive standard representation of $f$ with respect to the set
  $\bar\H$ proving our assertion.\qed
\end{proof}  

By Proposition~5.7 of Part~I, an involutive head autoreduction of the set
$\bar\H$ yields a strong Pommaret basis for the saturation $\Isat$.  As a
trivial consequence of the considerations in the proof above we find that in
$\delta$-regular coordinates $\Isat$ is simply given by the quotient
$\I:x_1^\infty$ (in the monomial case this also follows immediately from
Proposition \ref{prop:qstab} (iv)).  This observation in turn implies that for
degrees $q\geq\deg{\H_1}$ we have $\I_q=\Isat_q$.  Hence all ideals with the
same saturation possess also the same Hilbert polynomial and become identical
for sufficiently high degrees; $\Isat$ is the largest among all these ideals.
The smallest value $q_0$ such that $\I_q=\Isat_q$ for all $q\geq q_0$ is often
called the \emph{satiety} $\sat{\I}$ of the ideal $\I$.

\begin{corollary}\label{cor:sat}
  Let\/ $\H$ be a Pommaret basis of the ideal\/ $\I\subseteq\P$.  Then $\I$ is
  saturated, if and only if\/ $\H_1=\emptyset$.  If\/ $\I$ is not saturated,
  then\/ $\sat{\I}=\deg{\H_1}$.  Independent of the existence of a Pommaret
  basis, we have for any homogeneous generating set\/ $\F$ of the socle\/
  $\I:\P_+$ the equality
  \begin{equation}\label{eq:satbg}
    \sat{\I}=1+\max{\{\deg{f}\mid f\in\F\wedge f\notin\I\}}\;.
  \end{equation}
\end{corollary}

\begin{proof}
  Except of the last statement, everything has already been proven in the
  discussion above.  For its proof we may assume without loss of generality
  that the coordinates are $\delta$-regular so that a Pommaret basis $\H$ of
  $\I$ exists, as all quantities appearing in (\ref{eq:satbg}) are invariant
  under linear coordinate transformations.
  
  Let $\tilde h$ be an element of $\H_1$ having maximal degree.  We claim that
  then $\tilde h/x_1\in(\I:\P_+)\setminus\I$.  Indeed, since $x_1$ is always
  multiplicative for the Pommaret division, we cannot have $\tilde h/x_1\in\I$
  (otherwise $\H$ would not be involutively head autoreduced), and if we
  analyse for any $1<\ell\leq n$ the involutive standard representation of
  $x_\ell\tilde h$, then all coefficients of generators $h\in\H\setminus\H_1$
  are trivially contained in $\lspan{x_1}$ and for the coefficients of
  elements $h\in\H_1$ the same holds for degree reasons.  Hence we can divide
  by $x_1$ and find that $x_\ell\tilde h/x_1\in\I$ for all $1\leq\ell\leq n$
  implying $\tilde h/x_1\in\I:\P_+$.
  
  By our previous results, $\sat{\I}=\deg{\tilde h}$.  By assumption,
  homogeneous polynomials $P_f\in\P$ exist such that $\tilde
  h/x_1=\sum_{f\in\F}P_ff$.  If $\deg{P_f}\geq1$, then $P_ff\in\I$ since
  $f\in\I:\P_+$.  Hence for at least one $\tilde f\in\F\setminus\I$, the
  coefficient $P_{\tilde f}$ must be a non-zero constant and thus $\deg{\tilde
    f}=\sat{\I}-1$.\qed
\end{proof}

\begin{remark}\label{rem:satiety}
  The last statement in Corollary \ref{cor:sat} is due to Bermejo and Gimenez
  \cite[Prop.~2.1]{bg:scmr} who proved it in a slightly different way.  For
  \emph{monomial} ideals $\I$, one obtains as further corollary
  \cite[Cor.~2.4]{bg:scmr} that, if $\I:\P_+=\I:x_1$, then $\sat{\I}$ is the
  maximal degree of a minimal generator of $\I$ divisible by $x_1$ (this
  observation generalises a classical result about Borel-fixed ideals
  \cite[Cor.~2.10]{mlg:gin}).  If the considered ideal $\I$ possesses a
  Pommaret basis $\H$, this statement also follows from the fact that under
  the made assumptions all elements of $\H_1$ are minimal generators.  Indeed,
  suppose to the contrary that $\H_1$ contains two elements $h_1\neq h_2$ such
  that $h_1\mid h_2$.  Obviously, the minimality implies
  $\deg_{x_1}{h_1}=\deg_{x_1}{h_2}$ and a non-multiplicative index $1<\ell\leq
  n$ exists such that $x_\ell h_1\mid h_2$.  Without loss of generality, we
  may assume that $h_2=x_\ell h_1$.  But this immediately entails that $x_\ell
  h_1/x_1=h_2/x_1\notin\I$ and hence
  $h_1/x_1\in(\I:x_1)\setminus(\I:\P_+)$.\bull
\end{remark}

A first trivial consequence of our results is the following well-known formula
relating Castelnuovo-Mumford regularity and saturation.

\begin{corollary}
  Let\/ $\I\subseteq\P$ be an ideal.  Then\/
  $\reg{\I}=\max{\{\sat{\I},\reg{\Isat}\}}$.
\end{corollary}

\begin{proof}
  Without loss of generality, we may assume that we use $\delta$-regular
  coordinates so that $\I$ possesses a Pommaret basis $\H$ with respect to the
  degree reverse lexicographic order.  Now the statement follows immediately
  from Proposition \ref{prop:Isat} and Corollary \ref{cor:sat}.\qed
\end{proof}

Trung \cite{ngt:cmr} proposed the following approach for computing the
regularity of a monomial ideal $\I$ based on evaluations.  Let
$D=\dim{(\P/\I)}$ and introduce for $j=0,\dots,D$ the polynomial
subrings\footnote{Compared with Trung \cite{ngt:cmr}, we revert as usual the
  order of the variables in order to be consistent with our conventions.}
$\P^{(j)}=\kk[x_{j+1},\dots,x_n]$ and within them the elimination ideals
$\I^{(j)}=\I\cap\P^{(j)}$ and their saturations $\tilde\I^{(j)}=
\I^{(j)}:x_{j+1}^\infty$.  A basis of $\I^{(j)}$ is obtained by setting
$x_1=\cdots=x_j=0$ in a basis of $\I$ and for a basis of $\tilde\I^{(j)}$ we
must additionally set $x_{j+1}=1$.  Now define the numbers
\begin{subequations}
  \begin{align}
    c_j&=\sup{\bigl\{q\mid(\tilde\I^{(j)}/\I^{(j)})_q\neq0\bigr\}}+1\;,
         \qquad 0\leq j<D\\
    c_D&=\sup{\bigl\{q\mid(\P^{(D)}/\I^{(D)})_q\neq0\bigr\}}+1\;.
  \end{align}
\end{subequations}
Trung \cite{ngt:cmr} proved that whenever none of these numbers is infinite,
then their maximum is just $\reg{\I}$.  We show now that this genericity
condition is satisfied, if and only if the coordinates are $\delta$-regular
and express the numbers $c_j$ as satieties.

\begin{theorem}\label{thm:trung}
  The numbers\/ $c_0,\dots,c_D$ are all finite, if and only if the monomial
  ideal\/ $\I\subseteq\P$ is quasi-stable.  In this case\/
  $c_j=\sat{\I^{(j)}}$ for\/ $0\leq j\leq D$ and
  \begin{equation}\label{eq:trung}
    \max{\{c_0,\dots,c_D\}}=\reg{\I}\;.
  \end{equation}
  If\/ $d=\depth{\I}$, then it suffices to consider\/ $c_d,\dots,c_D$.
\end{theorem}

\begin{proof}
  We assume first that $\I$ is quasi-stable and thus possesses a Pommaret
  basis which we write $\H=\bigl\{h_{k,\ell}\mid 1\leq k\leq n,\ 
  1\leq\ell\leq\ell_k\bigr\}$ where $\cls{h_{k,\ell}}=k$.  One easily verifies
  that the subset $\H^{(j)}=\bigl\{h_{k,\ell}\in\H\mid k>j\bigr\}$ is the
  Pommaret basis of the ideal $\I^{(j)}$.  If we set
  $a_{k,\ell}=\deg_{x_k}{h_{k,\ell}}$, then the Pommaret basis of
  $\tilde\I^{(j)}$ is
  $\tilde\H^{(j)}=\H^{(j+1)}\cup\bigl\{h_{j+1,\ell}/x_{j+1}^{a_{j+1,\ell}}\mid
  1\leq\ell\leq\ell_{j+1}\bigr\}$.  This immediately implies that
  $c_j=\max{\bigl\{\deg{h_{j+1,\ell}}\mid 1\leq\ell\leq\ell_{j+1}\bigr\}}$.
  By construction, $\dim{(\P^{(D)}/\I^{(D)})}=0$ and
  Proposition~\ref{prop:dim} entails that for $\hat q=\deg{\H^{(D)}}$ the
  equality $\I^{(D)}_{\hat q}=\P^{(D)}_{\hat q}$ holds.  Hence $c_{D}=\hat q$
  (it is not possible that $c_{D}<\hat q$, as otherwise the set $\H$ was not
  involutively autoreduced).
  
  Thus we find that $\max{\{c_0,\dots,c_D\}}=\deg{\H}$ and Theorem
  \ref{thm:pomreg} yields (\ref{eq:trung}).  Furthermore, it follows
  immediately from Corollary \ref{cor:sat} and Proposition \ref{prop:dim},
  respectively, that $c_j=\sat{\I^{(j)}}$ for $0\leq j\leq D$.  Finally,
  Proposition \ref{prop:Iregular} entails that the values $c_0,\dots,c_{d-1}$
  vanish.
  
  Now assume that the ideal $\I$ was not quasi-stable.  By Part (ii) of
  Proposition~\ref{prop:qstab} this entails that for some $0\leq j<D$ the
  variable $x_{j+1}$ is a zero divisor in the ring
  $\P/\lspan{\I,x_1,\dots,x_j}^{\mathrm{sat}}\cong
  \P^{(j)}/(\I^{(j)})^{\mathrm{sat}}$.  Thus a polynomial
  $f\notin(\I^{(j)})^{\mathrm{sat}}$ exists for which
  $x_{j+1}f\in(\I^{(j)})^{\mathrm{sat}}$ which means that we can find for any
  sufficiently large degree $q\gg0$ a polynomial $g\in\P^{(j)}$ with
  $\deg{g}=q-\deg{f}$ such that $fg\notin\I^{(j)}$ but $x_{j+1}fg\in\I^{(j)}$.
  Hence the equivalence class of $fg$ is a non-vanishing element of
  $(\tilde\I^{(j)}/\I^{(j)})_q$ so that for a not quasi-stable ideal $\I$ at
  least one value $c_j$ is not finite.  \qed
\end{proof}

One direction of the proof above uses the same idea as the one of Theorem
\ref{thm:pomreg}: the Castelnuovo-Mumford regularity is determined by the
basis members of maximal degree and their classes give us the positions in the
minimal resolution where it is attained (recall Remark~\ref{rem:cmrpos}; here
these are simply the indices $j$ for which $c_j$ is maximal).  However, while
Theorem \ref{thm:pomreg} holds for arbitrary homogeneous ideals, Trung's
approach can only be applied to monomial ideals.  The formulation using
satieties is at the heart of the method of Bermejo and Gimenez \cite{bg:scmr}
to compute the Castelnuovo-Mumford regularity.  Similar considerations yield
an alternative proof of the following result of Bermejo and Gimenez
\cite[Cor.~17]{bg:scmr} for monomial ideals.

\begin{proposition}\label{prop:regirred}
  Let\/ $\I\subseteq\P$ be a quasi-stable ideal and\/
  $\I=\J_1\cap\cdots\cap\J_r$ its unique irredundant decomposition into
  irreducible monomial ideals.  Then the equality\/
  $\reg{\I}=\max{\{\reg{\J_1},\dots,\reg{\J_r}\}}$ holds.
\end{proposition}

\begin{proof}
  We first note that the Castelnuovo-Mumford regularity of a monomial
  irreducible ideal $\J=\lspan{x_{i_1}^{\ell_1},\dots,x_{i_k}^{\ell_k}}$ is
  easily determined using the considerations in Example 2.12 of Part I.  There
  we showed that any such ideal becomes quasi-stable after a simple
  renumbering of the variables and explicitly gave its Pommaret basis.  Up to
  the renumbering, the unique element of maximal degree in this Pommaret basis
  is the term $x_{i_1}^{\ell_1}x_{i_2}^{\ell_2-1}\cdots x_{i_k}^{\ell_k-1}$
  and thus it follows from Theorem \ref{thm:pomreg} that
  $\reg{\J}=\sum_{j=1}^k\ell_j-k+1$.
  
  Recall from Proposition \ref{prop:stanpair} that an irreducible
  decomposition can be constructed via standard pairs.  As discussed in
  Section \ref{sec:decomp}, the decomposition (\ref{eq:stanpair2}) is in
  general redundant; among all standard pairs $(\nu,N_\nu)$ with $N_\nu=N$ for
  some given set $N$ only those exponents $\nu$ which are maximal with respect
  to divisibility appear in the irredundant decomposition and thus are
  relevant.
  
  If we now determine the standard pairs of $\I$ from a Pommaret basis
  according to Remark \ref{rem:stanpom}, then we must distinguish two cases.
  We have first the standard pairs coming from the terms $x^\mu$ of degree
  $q=\deg{\H}$ not lying in $\I$.  They are of the form
  $\bigl(x^{\nu},\{x_1,\dots,x_k\}\bigr)$ where $k=\cls{\mu}$ and
  $x^{\nu}=x^\mu/x_k^{\mu_k}$.  By Proposition \ref{prop:stanpair}, each
  such standard pair leads to the irreducible ideal
  $\J=\lspan{x_\ell^{\nu_\ell+1}\mid k<\ell\leq n}$.  By the remarks
  above, $\reg{\J}=|\nu|+1\leq|\mu|=q=\reg{\I}$.
  
  The other standard pairs come from the terms $x^\nu\notin\I$ with $|\nu|<q$.
  It is easy to see that among these the relevant ones correspond one-to-one
  to the ``end points'' of the monomial completion process: we call an element
  of the Pommaret basis $\H$ of $\I$ an end point, if each non-multiplicative
  multiple of it has a \emph{proper} involutive divisor in the basis (and thus
  one branch of the completion process ends with this element\footnote{Note
    that an end point may very well be a member of the minimal basis of
    $\I$!}).  If $x^\mu\in\H$ is such an end point, then the corresponding
  standard pair consists of the monomial $x^{\nu}=x^\mu/x_k$ where
  $k=\cls{\mu}$ and the empty set and it yields the irreducible ideal
  $\J=\lspan{x_\ell^{\nu_\ell+1}\mid 1\leq\ell\leq n}$.  Thus we find again
  $\reg{\J}=|\nu|+1=|\mu|\leq q=\reg{\I}$.
  
  This proves the estimate
  $\reg{\I}\geq\max{\{\reg{\J_1},\dots,\reg{\J_r}\}}$.  The claimed equality
  follows from the observation that any element of degree $q$ in $\H$ must
  trivially be an end point and the corresponding standard pair yields then an
  irreducible ideal $\J$ with $\reg{\J}=q$.\qed
\end{proof}

The question of bounding the Castelnuovo-Mumford regularity of a homogeneous
ideal $\I$ in terms of the degree $q$ of an arbitrary generating set has
attracted quite some interest.  Hermann \cite{gh:endl} gave already very early
a doubly exponential bound; much later Mayr and Meyer \cite{mm:word} showed
with explicit examples that this bound is indeed sharp (see \cite{bm:cac} for
a more detailed discussion).

For monomial ideals $\I$ the situation is much more favourable.  It follows
immediately from Taylor's explicit resolution of such ideals \cite{dt:res}
(see \cite{wms:taylor} for a derivation via Gr\"obner bases) that here a
\emph{linear} bound
\begin{equation}
  \reg{\I}\leq n(q-1)+1
\end{equation}
holds where $n$ is again the number of variables.  Indeed, this resolution is
supported by the $\mathrm{lcm}$-lattice of the given basis and the degree of
its $k$th term is thus trivially bounded by $kq$.  By Hilbert's Syzygy
Theorem, it suffices to consider the first $n$ terms which immediately yields
the above bound.  If the ideal $\I$ is even quasi-stable, a simple corollary
to Proposition \ref{prop:regirred} yields an improved bound using the minimal
generators of $\I$.

\begin{corollary}\label{cor:cmrbound}
  Let the monomials\/ $m_1,\dots,m_r$ be the minimal generators of the
  quasi-stable ideal\/ $\I\subseteq\kk[x_1,\dots,x_n]$.  If we set\/
  $x^\lambda=\lcm{m_1,\dots}{m_r}$ and\/
  $d=\min{\{\cls{m_1},\dots,\cls{m_r}\}}$ (i.\,e.\ $d=\depth{\I}$), then the
  Castelnuovo--Mumford regularity of\/ $\I$ satisfies the estimate
  \begin{equation}\label{eq:cmrbound}
    \reg{\I}\leq|\lambda|+d-n
  \end{equation}
  and this bound is sharp.
\end{corollary}

\begin{proof}
  Applying repeatedly the rule
  $\lspan{\F,t_1t_2}=\lspan{\F,t_1}\cap\lspan{\F,t_2}$ for arbitrary
  generating sets $\F$ and coprime monomials $t_1$, $t_2$, one obtains an
  irreducible decomposition of $\I$.  Obviously, in the worst case one of the
  irreducible ideals is $\J=\lspan{x_d^{\lambda_d},\dots,x_n^{\lambda_n}}$.
  As we already know that $\reg{\J}=|\lambda|+d-n$, this value bounds
  $\reg{\I}$ by Proposition \ref{prop:regirred}.\qed
\end{proof}

\begin{remark}
  An alternative direct proof of the corollary goes as follows.  Let $\H$ be
  the Pommaret basis of $\I$.  We claim that each generator $x^\mu\in\H$ with
  $\cls{\mu}=k$ satisfies $\mu_k\leq\lambda_k$ and $\mu_j<\lambda_j$ for all
  $j>k$.  The estimate for $\mu_k$ is obvious, as it follows immediately from
  our completion algorithm that there is a minimal generator $x^\nu\mid x^\mu$
  with $\nu_k=\mu_k$.
  
  Assume for a contradiction that the Pommaret basis $\H$ contains a generator
  $x^\mu$ where $\mu_j>\lambda_j$ for some $j>\cls{\mu}$.  If several such
  generators exist for the same value $j$, choose one for which $\mu_j$ is
  maximal.  Obviously, $j$ is non-multiplicative for $x^\mu$ and hence the
  multiple $x_jx^\mu$ must contain an involutive divisor $x^\nu\in\H$.
  Because of our maximality assumption $\nu_j\leq\mu_j$ and hence $j$ must be
  multiplicative for $x^\nu$ so that $\cls{\nu}\geq j$.  But this fact
  trivially implies that $x^\nu\idiv{P}x^\mu$ contradicting that $\H$ is by
  definition involutively autoreduced.
  
  Now the assertion follows immediately: under the made assumptions
  $\cls{\lambda}=d$ and in the worst case $\H$ contains the generator
  $x_d^{\lambda_d}x_{d+1}^{\lambda_{d+1}-1}\cdots x_n^{\lambda_n-1}$ which is
  of degree $|\lambda|+d-n$.\bull
\end{remark}

\begin{remark}\label{rem:satbound}
  The same arguments together with Proposition \ref{prop:Isat} also yield
  immediately a bound for the satiety of a quasi-stable ideal $\I$.  As
  already mentioned above, a quasi-stable ideal is not saturated, if and only
  if $d=1$.  In this case, we have trivially $\sat{\I}\leq|\lambda|+1-n$.
  Again the bound is sharp, as shown by exactly the same class of irreducible
  ideals as considered above.
  
  The estimate (\ref{eq:cmrbound}) also follows immediately from the results
  in \cite{bg:scmr}.  Yet another derivation is contained in \cite{ah:hilser}.
  \bull
\end{remark}

If one insists on having an estimate involving only the maximal degree $q$ of
the minimal generators and the depth, then the above result yields immediately
the following estimate, variations of which appear in
\cite{aa:reg,cs:reg,mc:reg}.

\begin{corollary}\label{prop:regbound}
  Let\/ $\I\subseteq\P$ be a quasi-stable ideal minimally generated in degrees
  less than or equal to\/ $q$.  If\/ $\depth{\I}=d$, then
  \begin{equation}\label{eq:regbound}
    q\leq\reg{\I}\leq (n-d+1)(q-1)+1
  \end{equation}
  and both bounds are sharp.
\end{corollary}

\begin{proof}
  Under the made assumptions we trivially find that the degree of the least
  common multiple of the minimal generators is bounded by $|\lambda|\leq
  (n-d+1)q$.  Now (\ref{eq:regbound}) follows immediately from
  (\ref{eq:cmrbound}).  The upper bound is realised by the irreducible ideal
  $\I=\lspan{x_1^q,\dots,x_n^q}$.  The lower bound is attained, if $\I$ is
  even stable, as then Proposition \ref{prop:invstab} implies that
  $\reg{\I}=q$ independent of $\depth{\I}$.\qed
\end{proof}

\begin{remark}
  Eisenbud, Reeves and Totaro \cite{ert:veronese} presented a variation of the
  estimate (\ref{eq:regbound}).  They introduced the notion of $s$-stability
  as a generalisation of stability: let $s\geq1$ be an integer; a monomial
  ideal $\I$ is $s$-stable, if for every monomial $x^\mu\in\I$ and every index
  $n\geq j>\cls{\mu}=k$ an exponent $1\leq e\leq s$ exists such that
  $x^{\mu-e_k+e_j}\in\I$.  Then it is easy to see that for an $s$-stable ideal
  generated in degrees less than or equal to $q$ the estimate
  \begin{equation}
    \reg{\I}\leq q+(n-1)(s-1)
  \end{equation}
  holds, as $\I_{\geq q+(n-1)(s-1)}$ is stable (thus any $s$-stable ideal is
  trivially quasi-stable).  However, in general this is an overestimate, as it
  based on the assumption that $\I$ possesses a minimal generator of class $1$
  and degree $q$ which must be multiplied by $x_2^{s-1}x_3^{s-1}\cdots
  x_n^{s-1}$ in order to reach a stable set.
  
  Thus for the $8$-stable ideal $\lspan{x^8,y^8,z^8}$ the estimate is indeed
  sharp (this is exactly the same worst case as in the proof above for an
  ideal of depth $1$); the Pommaret basis contains as maximal degree element
  the monomial $x^8y^7z^7$.  On the other hand, for the also $8$-stable ideal
  $\lspan{x^6,x^2y^4,x^2z^4,y^8,z^8}$ the regularity is only $16$, as now the
  maximal degree element of the Pommaret basis is $x^2y^7z^7$.\bull
\end{remark}

Finally, we recall that, given two quasi-stable ideals $\I,\J\subseteq\P$ and
their respective Pommaret bases, we explicitly constructed in Remarks 2.9 and
6.5, respectively, of Part I weak Pommaret bases for the sum $\I+\J$, the
product $\I\cdot\J$ and the intersection $\I\cap\J$.  They lead to the
following estimates for the regularity of these ideals which were recently
also given by Cimpoea{\c s} \cite{mc:borel,mc:stable}.

\begin{proposition}
  Let\/ $\I,\J\subseteq\P$ be two quasi-stable ideals.  Then the following
  three estimates hold:
  \begin{subequations}
    \begin{eqnarray}
      \reg{(\I+\J)}&\leq&\max{\{\reg{\I},\reg{\J}\}}\;,\\
      \reg{(\I\cdot\J)}&\leq&\reg{\I}+\reg{\J}\;,\\
      \reg{(\I\cap\J)}&\leq&\max{\{\reg{\I},\reg{\J}\}}\;.
    \end{eqnarray}
  \end{subequations}
\end{proposition}

\begin{proof}
  The first two estimates follow immediately from the weak Pommaret bases
  given in the above mentioned remarks and Theorem \ref{thm:pomreg}.  For the
  last estimate the weak Pommaret basis constructed in Remark 6.5 of Part I is
  not good enough; it would also yield $\reg{\I}+\reg{\J}$ as upper bound.
  However, Lemma \ref{lem:truncbas} allows us to improve it significantly.  Let
  $\G$ be the Pommaret basis of $\I$ and $\H$ the one of $\J$.  If we set
  $q=\max{\{\deg{\G},\deg{\H}\}}$, then one easily sees that $\G_q\cap\H_q$ is
  the Pommaret basis of $(\I\cap\J)_{\geq q}$.  Hence, the intersection
  $\I\cap\J$ possesses a Pommaret basis of degree at most $q$.\qed
\end{proof}

\section{Iterated Polynomial Algebras of Solvable Type}\label{sec:iteralg}

In Section 11 of Part I we studied involutive bases in polynomial algebras of
solvable type over rings.  We had to substitute the notion of an involutively
head autoreduced set by the more comprehensive concept of an involutively
$\R$-saturated set.  In a certain sense this was not completely satisfying, as
we had to resort here to classical Gr\"obner techniques, namely computing
normal forms of ideal elements arising from syzygies.  Using the syzygy theory
developed in Section \ref{sec:syz}, we provide now an alternative approach for
the special case that the coefficient ring $\R$ is again a polynomial algebra
of solvable type (over a field).  It is obvious that in this case left ideal
membership in $\R$ can be decided algorithmically and by
Theorem~\ref{thm:invschreyer} it is also possible to construct algorithmically
a basis of the syzygy module.

\begin{remark}\label{rem:ncsyz}
  In Section \ref{sec:syz} we only considered the ordinary commutative
  polynomial ring, whereas now we return to general polynomial algebras of
  solvable type (over a field).  However, it is easy to see that all the
  arguments in the proof of the involutive Schreyer
  Theorem~\ref{thm:invschreyer} depend only on normal form computations and on
  considerations concerning the leading exponents.  The same holds for the
  classical Schreyer theorem, as one may easily check (see also
  \cite{vl:ade,vl:diss} for a non-commutative version).  Thus the in the
  sequel crucial Theorem~\ref{thm:invschreyer} remains valid in the general
  case of non-commutative polynomial algebras.\bull
\end{remark}

We use the following notations in this section:
$\R=(\kk[y_1,\dots,y_m],\star,\prec_y)$ and
$\P=(\R[x_1,\dots,x_n],\star,\prec_x)$.  Furthermore, we are given an
involutive division $L_y$ on $\NN_0^m$ and a division $L_x$ on $\Nn$.  For
simplicity, we always assume in the sequel that at least $L_y$ is Noetherian.
In order to obtain a reasonable theory, we make similar assumptions as in
Section 11 of Part I: both $\R$ and $\P$ are solvable algebras with centred
commutation relations so that both are (left) Noetherian.

We now propose an alternative algorithm for the involutive $\R$-saturation.
Until Line /13/ it is identical with Algorithm~6 of Part~I; afterwards we
perform an involutive completion and multiply in Line /17/ each polynomial in
$\bar\H'_{f,L_x}$ by the non-multiplicative variables of its leading
coefficient.  In the determination of involutive normal forms, we may multiply
each polynomial $h'\in\H'$ only by monomials $rx^\mu$ such that
$x^\mu\in\R[\mult{X}{L_x,\H',\prec_x}{h'}]$ \emph{and}
$r\in\kk\bigl[\mult{Y}{L_y,\lc{\prec_x}{(\bar\H'_{h',L_x})},\prec_y}
{\lc{\prec_x}{h'}}\bigr]$.

\begin{algorithm}
  \caption{Involutive $\R$-saturation (and head autoreduction)\label{alg:rsat}}
  \begin{algorithmic}[1]
    \REQUIRE finite set $\F\subset\P$, involutive divisions $L_y$ on $\NN_0^m$ 
             and $L_x$ on $\Nn$
    \ENSURE involutively $\R$-saturated and head autoreduced set $\H$ with
            $\lspan{\H}=\lspan{\F}$
    \STATE $\H\leftarrow\F$;\quad $\S\leftarrow\F$
    \WHILE{$S\neq\emptyset$}
        \STATE $\nu\leftarrow\max_{\prec_x}{\le{\prec_x}{\S}}$;\quad 
               $\S_\nu\leftarrow\{f\in\H\mid\le{\prec_x}{f}=\nu\}$
        \STATE $\S\leftarrow\S\setminus\S_\nu$;\quad $\H'\leftarrow\H$
        \FORALL{$f\in\S_\nu$}
            \STATE $h\leftarrow\mathtt{HeadReduce}_{L_x,\prec_x}(f,\H)$
            \IF{$f\neq h$}
                \STATE $\S_\nu\leftarrow\S_\nu\setminus\{f\}$;\quad
                       $\H'\leftarrow\H'\setminus\{f\}$
                \IF{$h\neq0$}
                    \STATE $\H'\leftarrow\H'\cup\{h\}$
                \ENDIF
            \ENDIF
        \ENDFOR
        \IF{$\S_\nu\neq\emptyset$}
            \STATE choose $f\in\S_\nu$ and determine the set $\bar\H'_{f,L_x}$
            \REPEAT
                \STATE $\T\leftarrow\bigl\{y_j\star\bar f\mid
                             \bar f\in\bar\H'_{f,L_x},\ 
                             y_j\in\nmult{Y}{L_y,
                                             \lc{\prec_x}{(\bar\H'_{f,L_x})},
                                             \prec_y}
                                         {\lc{\prec_x}{\bar f}}\bigr\}$
                \REPEAT
                    \STATE choose $h'\in\T$ such that 
                           $\le{\prec_y}{(\lc{\prec_x}{h'})}$
                           is minimal
                    \STATE $\T\leftarrow\T\setminus\{h'\}$
                    \STATE $h\leftarrow
                            \mathtt{NormalForm}_{L_x,\prec_x,L_y,\prec_y}
                                   (h',\H')$
                    \IF{$h\neq0$}
                        \STATE $\H'\leftarrow\H'\cup\{h\}$
                    \ENDIF
                \UNTIL{$\T=\emptyset\vee h\neq0$}
            \UNTIL{$\T=\emptyset\wedge h=0$}
        \ENDIF
        \IF{$\H'\neq\H$}
            \STATE $\H\leftarrow\H'$;\quad $\S\leftarrow\H$
        \ENDIF
    \ENDWHILE
    \STATE \algorithmicreturn{$\H$}
  \end{algorithmic}
\end{algorithm}

\begin{proposition}\label{prop:rsat}
  Let\/ $L_y$ be a Noetherian constructive division.  Algorithm~\ref{alg:rsat}
  terminates for any input\/ $\F$ with an involutively\/ $\R$-saturated and
  head autoreduced set\/ $\H$ such that\/ $\lspan{\H}=\lspan{\F}$.
  Furthermore, the sets\/ $\lc{\prec_x}{\bar\H_{h,L_x}}$ form weak\/
  $L_y$-involutive bases of the\/ $\R$-ideals generated by them for each\/
  $h\in\H$.
\end{proposition}

\begin{proof}
  The termination criterion in Line /26/ is equivalent to local involution of
  all the sets $\lc{\prec_x}{\bar\H'_{f,L_x}}$.  Under the made assumptions on
  the division $L_y$ and because of the fact that $\P$ is Noetherian, the
  termination of the algorithm and the assertion about these sets is obvious.
  In general we only obtain weak involutive bases, as no involutive head
  autoreductions of these sets are performed.  The correctness is a
  consequence of Theorem~\ref{thm:invschreyer}: by analysing all
  non-multiplicative products we have taken into account a whole basis of the
  syzygy module $\syz{}{\lc{\prec_x}{\bar\H'_{f,L_x}}}$.  Thus the output $\H$
  is indeed involutively $\R$-saturated.\qed
\end{proof}

\begin{theorem}\label{thm:rsat}
  Let\/ $\P$ satisfy the made assumptions and\/ $L_x$ be a Noetherian
  division.  If in Algorithm 3 of Part I the subalgorithm\/
  $\mathtt{InvHeadAutoReduce}_{L_x,\prec_x}$ is substituted by
  Algorithm~\ref{alg:rsat}, then the completion will terminate with a weak
  involutive basis of\/ $\I=\lspan{\F}$ for any finite input set\/
  $\F\subset\P$ such that the monoid ideal\/ $\le{\prec_x}{\I}$ possesses a
  weak involutive basis.  Furthermore, the sets\/
  $\lc{\prec_x}{\bar\H_{h,L_x}}$ form strong\/ $L_y$-involutive bases of the\/
  $\R$-ideals generated by them for each\/ $h\in\H$.
\end{theorem}

\begin{proof}
  The proof of the termination and of the correctness of the algorithm is as
  in Part~I.  The only new claim is that the sets
  $\lc{\prec_x}{\bar\H_{h,L_x}}$ are strongly $L_y$-involutive.  This is a
  simple consequence of the fact that under the made assumption on the product
  in $\P$ the loop in Lines /5-13/ of Algorithm~\ref{alg:rsat} leads to an
  involutive head autoreduction of these sets.  Hence we indeed obtain strong
  involutive bases.  \qed
\end{proof}

\begin{corollary}
  If\/ $L_x$ is the Janet division, then each polynomial\/ $f\in\I$ possesses
  a unique involutive standard representation\/ $f=\sum_{h\in\H}P_h\star h$
  where\/ $P_h\in\kk[\mult{Y}{L_y,\lc{\prec_x}{(\bar\H_{h,L_x})},\prec_y}
  {\lc{\prec_x}{h}}][\mult{X}{L_x,\H,\prec_x}{h}]$.
\end{corollary}

\begin{proof}
  For the Janet division the only obstruction for $\H$ being a strong
  involutive basis is that some elements of it may have the same leading
  exponents.  More precisely, for any $h\in\H$ we have
  $\H_{h,L_x}=\{h'\in\H\mid\le{\prec_x}{h'}=\le{\prec_x}{h}\}$.  This
  immediately implies furthermore $\bar\H_{h,L_x}=\H_{h,L_x}$.  By
  Theorem~\ref{thm:rsat} the sets $\lc{\prec_x}{\bar\H_{h,L_x}}$ form a strong
  $L_y$-involutive basis of the ideals they generate.  Hence the claimed
  representation must be unique.\qed
\end{proof}

\section{Conclusions}\label{sec:conc}

$\delta$-regularity is often considered as a purely technical nuisance and a
problem specific to the Pommaret division.  Our results here lead to a more
intricate picture.  $\delta$-regularity of a given generating set in the sense
of Definition \ref{def:dregalg} is indeed a technical concept appearing in the
analysis of the termination of the completion Algorithm~3 of Part I.  By
contrast, $\delta$-regularity of an ideal in the sense of Definition
\ref{def:dregid} has an intrinsic meaning.  This can already be seen from the
simple fact that in the case of linear differential operators there is a close
relation to characteristics (see any textbook on partial differential
equations, e.\,g.\ \cite{jo:pde,rr:intro}): a necessary condition for a
coordinate system to be $\delta$-regular is that the hypersurface $x_n=0$ is
non-characteristic.  Indeed, the standard definition of a characteristic
hypersurface may be rephrased that on it one cannot solve for all derivatives
of class $n$.

We have also seen that $\delta$-regularity is related to many regularity
concepts in commutative algebra and algebraic geometry.  Many statements that
are only generically true hold in $\delta$-regular coordinates.  In
particular, in $\delta$-regular coordinates many properties of an affine
algebra $\A=\P/\I$ may already be read off the monomial algebra
$\A'=\P/\lt{\prec}{\I}$ where $\prec$ is the degree reverse lexicographic
order.

For example, it follows immediately from Proposition~\ref{prop:depth} that
$\depth{\A}=\depth{\A'}$ and that $(x_1,\dots,x_d)$ is a maximal regular
sequence for both algebras.  As in the homogeneous case it is also trivial
that $\dim{\A}=\dim{\A'}$, we see that the algebra $\A$ is Cohen-Macaulay, if
and only if $\A'$ is so.  Similarly, it is an easy consequence of
Theorem~\ref{thm:pd} that $\pd{\A}=\pd{\A'}$ and of Theorem~\ref{thm:pomreg}
that $\reg{\A}=\reg{\A'}$.  An exception are the Betti numbers where
Example~\ref{ex:stableres} shows that even in $\delta$-regular coordinates
$\A$ and $\A'$ may have different ones.

These equalities are of course not new; they can already be found in
\cite{bs:mreg} (some even in earlier references).  However, one should note an
important difference: Bayer and Stillman \cite{bs:mreg} work with the
\emph{generic initial ideal}, whereas we assume $\delta$-regularity of the
coordinates.  These are two different genericity concepts, as even in
$\delta$-regular coordinates $\lt{\prec}{\I}$ is not necessarily the generic
initial ideal (in contrast to the former, the latter is always Borel fixed).

When we proved in Corollary \ref{cor:pdlt} and \ref{cor:reglt}, respectively,
the two inequalities $\pd{\A}\leq\pd{\A'}$ and
$\reg{\I}\leq\reg{(\lt{\prec}{\I})}$ for arbitrary term orders $\prec$, we had
to assume the existence of a Pommaret basis of $\I$ for $\prec$.  It is
well-known that these inequalities remain true, if we drop this assumption
(see for example the discussions in \cite{bm:cac,bs:mreg,bc:initial}).  We
included here our alternative proofs because of their great simplicity and
they cover at least the generic case.

In view of these observations, it seems to be of interest for the structure
analysis of polynomial modules to construct explicitly $\delta$-regular
coordinates.  The approaches presented in Sections~\ref{sec:pombas} and
\ref{sec:noether} offer an alternative to the classical expensive methods like
random or generic coordinates.  For coordinate systems not too far away from
$\delta$-regularity, they should be much more efficient and in particular
destroy much less sparsity.  Of course, one may also combine the random
strategy with our criterion of $\delta$-singularity in order to ensure that
the transformation has indeed resulted in a good coordinate system.  We also
note that $\delta$-regularity of the coordinates $x_1,\dots,x_n$ for an ideal
$\I$ is equivalent to their quasi-regularity for the algebra $\P/\I$ in the
sense of Serre (see his letter published as an appendix to the article
\cite{gs:alg} by Guillemin and Sternberg) \cite{wms:delta}.

As already mentioned, many of the results in Sections
\ref{sec:polyres}--\ref{sec:minres} are generalisations of the work of Eliahou
and Kervaire \cite{ek:res}.  They considered exclusively the case of stable
modules where we obtain a minimal resolution.  If one analyses closely their
proofs, it is not difficult to see that implicitly they introduce Pommaret
bases and exploit some of their basic properties.  Our proof of Theorem
\ref{thm:minres} appears so much simpler only because we have already shown
all these properties in Part I.  Furthermore, Eliahou and Kervaire did not
realise that they constructed a syzygy resolution in Schreyer form.  Hence
they had to give a lengthy and rather messy proof that the complex
$(\S_*,\delta)$ is exact, whereas in our approach this is immediate.

We rediscover all their complicated calculations in the proof of Theorem
\ref{thm:pomdiffmon}.  But note that this explicit formula for the
differential is needed neither for proving the minimality of the resolution
nor for its construction, although the latter is of course simplified by it.
Furthermore, the theory of involutive standard representations gives us a
clear guideline how to proceed.

Our results strongly suggest a homological background of the Pommaret
division.  Indeed, most of the quantities like the depth or the
Castelnuovo-Mumford regularity determined by the Pommaret basis of an ideal
$\I$ are of a homological nature; more precisely, they correspond to certain
extremal points in the Betti diagram and thus come from the Koszul homology.
It is a conjecture of us that knowledge of the Pommaret basis (for the degree
reverse lexicographic term order) of $\I$ is equivalent to knowing the full
Koszul homology of $\I$ and that it is possible to construct explicitly one
from the other.  In the special case of monomial ideals, Sahbi showed recently
in his diploma thesis \cite{ms:dipl} how the Koszul homology of a quasi-stable
ideal can be computed from the $P$-graph of its Pommaret basis.

The combination of Corollary \ref{cor:stacon} and Proposition \ref{prop:depth}
allows us to make some statements about the so-called \emph{Stanley
  conjecture}.  It concerns the minimal number of multiplicative variables for
a generator in a Stanley decomposition.  Following Apel
\cite[Def.~1]{apel:stanley1} and Herzog et al. \cite{hsy:stanley} we call this
number the \emph{Stanley depth} of the decomposition and for an ideal
$\I\subseteq\P$ the Stanley depth of the algebra $\A=\P/\I$, written
$\sdepth{\A}$, is defined as the maximal Stanley depth of a complementary
decomposition for $\I$.  In its simplest form the Stanley conjecture claims
that we always have the inequality $\sdepth{\A}\geq\depth{\A}$.  Obviously,
Corollary~\ref{cor:stacon} together with Proposition \ref{prop:depth} (plus
the existence Theorem \ref{thm:expombas} for Pommaret bases) shows that this
inequality holds for arbitrary ideals.

The rigorous formulation of the Stanley conjecture
\cite[Conj.~5.1]{rps:diophant} concerns monomial ideals and requires that all
generators in the decomposition are again monomials.  Furthermore, no
variables transformation is allowed.  Then our results only allow us to
conclude that the Stanley conjecture is true for all quasi-stable ideals.
Some further results on this question have been achieved by Apel
\cite{apel:stanley1,apel:stanley2} with the help of a slightly different
notion of involutive bases.\footnote{\cite{apel:stanley1} considers the
  Stanley conjecture for the ideal $\I$ itself instead of the factor algebra
  $\A=\P/\I$.  Here it follows immediately from the definition of a Pommaret
  basis and Proposition \ref{prop:Iregular} that the Stanley conjecture is
  true in its weak form for arbitrary polynomial ideals and in its strong form
  for all quasi-stable ideals.}

Many of the results mentioned above are quite well-known for Borel-fixed
ideals and thus for generic initial ideals.  However, it appears that for many
purposes it is not necessary to move to this highly special class of ideals;
quasi-stable ideals which are easier to produce algorithmically share many of
their properties.  Thus it is not surprising that quasi-stable ideals have
appeared under different names in quite a number of recent works in
commutative algebra (e.\,g.\ \cite{bg:scmr,cs:reg,hpv:ext}).

The results presented in this article offer two heuristic explanations for the
efficiency of the involutive completion algorithm already mentioned in Part I.
The first one is that according to our proof of Theorem \ref{thm:invschreyer}
the involutive algorithm automatically takes into account many instances of
Buchberger's two criteria for redundant $S$-polynomials.  Whereas a naive
implementation of Buchberger's algorithm without these criteria fails already
for rather small examples, a naive implementation of the involutive completion
algorithm works reasonably for not too large examples.

The second explanation concerns Proposition \ref{prop:hilbert}.  It is
well-known that the so-called ``Hilbert driven'' Buchberger algorithm
\cite{ct:hdba} is often very fast, but it requires a priori knowledge of the
Hilbert polynomial.  The involutive completion algorithm may also be
interpreted as ``Hilbert driven''.  The assignment of the multiplicative
variables defines at each iteration a trial Hilbert function.  This trial
function is the true Hilbert function, if and only if we have already reached
an involutive basis; otherwise it yields too small values.  For continuous
divisions the analysis of the products of the generators with their
non-multiplicative variables represents a simple check for the trial Hilbert
function to be the true one.

While for many ideals the involutive approach is an interesting alternative
for the construction of Gr\"obner bases, there exist some obvious cases where
this is not the case.  The first class are monomial ideals.  Here any basis is
already a Gr\"obner basis, whereas an involutive basis still has to be
constructed.  Another class are toric ideals where recent work by Blinkov and
Gerdt \cite{bg:toric} showed that involutive bases are typically much larger
than Gr\"obner bases.  In both cases, the reason is that these ideals are
rarely in general position and that hence Gr\"obner bases of a much lower
degree than $\reg{\I}$ exist.

An interesting question is whether the results of this second part can be
extended to the polynomial algebras of solvable type introduced in the first
part.  This is trivial only for the determination of Stanley decompositions,
as they are defined as vector space isomorphisms and therefore do not feel the
non-commutativity (note that we always have the commutative product $\cdot$ on
the right hand side of the defining equation (\ref{eq:sd}) of a Stanley
decomposition).  Thus involutive bases are a valuable tool for computing
Hilbert functions even in the non-commutative case.  This also yields
immediately the \emph{Gelfand-Kirillov dimension}
\cite[Sect.~8.1.11]{mr:ncrings} as the degree of the Hilbert polynomial (only
in the commutative case it always coincides with the Krull dimension).  Some
examples for such computations in the context of quantum groups (however,
using Gr\"obner instead of involutive bases) may be found in
\cite{bglc:quantum}.

By contrast, our results on the depth and on the Castelnuovo-Mumford
regularity rely on the fact that for commutative polynomials $f\in\lspan{x_j}$
implies that any term in $f$ is divisible by $x_j$.  In a non-commutative
algebra of solvable type we have the relations $x_i\star
x_j=c_{ij}x_ix_j+h_{ij}$ and in general the polynomial $h_{ij}$ is not
divisible by $x_j$.

For syzygies the situation is complicated, too.  The proof of Theorem
\ref{thm:invschreyer} is independent of the precise form of the multiplication
and thus we may conclude that we can always construct at least a Gr\"obner
basis of the syzygy module.  Our proof of Theorem \ref{thm:pomschreyer} relies
mainly on normal form arguments that generalise.  A minimal requirement is
that the term order $\prec_{\H}$ respects the multiplication $\star$, as
otherwise the theorem does not even make sense.  Furthermore, we must be
careful with all arguments involving multiplicative variables.  We need that
if $x_i$ and $x_j$ are both multiplicative for a generator, then $x_i\star
x_j=c_{ij}x_ix_j+h_{ij}$ must also contain only multiplicative variables which
will surely happen, if $h_{ij}$ depends only on variables $x_k$ with
$k\leq\max{\{i,j\}}$.  This is for example the case for linear differential
operators, so that we may conclude that Theorem \ref{thm:pomschreyer} (and its
consequences) remain true for the Weyl algebra and other rings of differential
operators. 

\begin{example}\label{ex:so3}
  Recall from Example~3.9 of Part I that the universal enveloping algebra of
  the Lie algebra $\mathfrak{so}(3)$ is isomorphic to the ring
  $(\kk[x_1,x_2,x_3],\star)$ with the product $\star$ induced by the relations
  \begin{equation}\label{eq:so3mult}
    \begin{aligned}
      x_1\star x_2&=x_1x_2\;,\qquad & x_2\star x_1&=x_1x_2-x_3\;,\\
      x_1\star x_3&=x_1x_3\;, & x_3\star x_1&=x_1x_3+x_2\;,\\
      x_2\star x_3&=x_2x_3\;, & x_3\star x_2&=x_2x_3-x_1\;.
    \end{aligned}
  \end{equation}
  Obviously $x_1x_2-x_3\in\lspan{x_1}$, but the term $x_3$ is not divisible by
  $x_1$.  It follows from the same relation that $x_2\star x_1$ depends on
  $x_3$ and thus the arguments on multiplicative variables required by our
  proof of Theorem \ref{thm:pomschreyer} break down.\bull
\end{example}

\appendix
\section{Rees Decompositions \`a la Sturmfels-White}\label{sec:rees}

Sturmfels and White \cite{sw:comb} presented an algorithm for the effective
construction of Rees decompositions (based on earlier works by Baclawski and
Garsia \cite{bg:comb}).  We show now that generically it yields a Pommaret
basis.  However, we believe that the involutive approach is much more
efficient.  It does not only allow us to avoid completely computations in
factor algebras, using our results in Section~\ref{sec:pombas} we obtain more
easily and constructively the right coordinates whereas Sturmfels and White
must rely on a probabilistic approach.

We introduce some additional notations.  Let again $\M$ be a finitely
generated module over the ring $\P=\kk[x_1,\dots,x_n]$.  The
\emph{annihilator} of an element $\fv\in\M$ is $\ann(\fv)=\{g\in\P\mid
g\fv=0\}$.  The $\kk$-vector space $Z_{\M}\subseteq\M$ is defined as the set
$Z_{\M}=\{\fv\in\M\mid\ann(\fv)=\lspan{x_1,\dots,x_n}\}$.  The approach of
Sturmfels and White is based on the following fact which may be interpreted as
their version of the concept of $\delta$-regularity.

\begin{lemma}\label{lem:nzd}
  If\/ $Z_{\M}=0$, there exists a non zero divisor\/~$y\in\P_1$, i.\,e.\ 
  $y\fv=0$ implies\/ $\fv=0$ for all\/ $\fv\in\M$.  Identifying\/ $\P_1$
  with\/ $\kk^n$, the set of all non zero divisors contains a Zariski open
  subset.
\end{lemma}

The Sturmfels-White Algorithm~\ref{alg:swa} computes a basis
$\Y=\{y_1,\dots,y_n\}$ of $\P_1$ and a set $\H\subset\M$ of generators such
that (as graded $\kk$-vector spaces)
\begin{equation}\label{eq:rd}
  M\cong\bigoplus_{\hv\in\H}\kk[y_1,\dots,y_{\cls{\hv}}]\cdot\hv\,.
\end{equation}
Here the class $\cls{\hv}$ is automatically assigned in the course of the
algorithm and not necessarily equal to the notion of class we introduced in
the definition of the Pommaret division.  Within this appendix, the latter one
will be referred to as \emph{coordinate class} $\ccls{\hv}{\xv}$, since its
definition depends on the chosen coordinates $\xv$.

\begin{algorithm}
  \caption{Construction of a Rees Decomposition \`a la 
           Sturmfels-White\label{alg:swa}}
  \begin{algorithmic}[1]
     \REQUIRE polynomial module $\M$ over $\P=\kk[x_1,\dots,x_n]$
     \ENSURE basis $\Y$ of $\P_1$, set $\H$ of generators defining
             Rees decomposition (\ref{eq:rd})
     \STATE $k\leftarrow0$;\quad $p\leftarrow0$;\quad $\M'\leftarrow\M$
     \WHILE{$\M'\neq0$}
         \STATE compute $Z_{\M'}$
         \IF{$Z_{\M'}=0$}
             \STATE $k\leftarrow{k+1}$
             \STATE choose a non zero divisor $y_k\in\P_1$ linearly
                    independent of $\{y_1,\dots,y_{k-1}\}$ 
             \STATE $\M'\leftarrow\M'/y_k\M'$
         \ELSE
             \STATE compute $\hv_{p+1},\dots,\hv_{p+\ell}\in\M$ such that
                    $\bigl\{[\hv_i]\mid p<i\leq p+\ell\bigr\}$\\
                    is a basis of $Z_{\M'}$
             \FOR{$i$ \algorithmicfont{from} $p+1$ \algorithmicfont{to}
                  $p+\ell$}
                  \STATE $\cls{\hv_i}\leftarrow k$
             \ENDFOR
             \STATE $p\leftarrow{p+\ell}$
             \STATE $\M'\leftarrow\M'/Z_{\M'}$
         \ENDIF
     \ENDWHILE
     \IF{$k<n$}
         \STATE complete $\{y_1,\dots,y_k\}$ to a basis $\Y$ of $\P_1$
     \ENDIF
     \STATE \algorithmicreturn{$\bigl(\Y,
               \H=\bigl\{(\hv_1,\cls{\hv_1}),\dots,
                         (\hv_p,\cls{\hv_p})\bigr\}\bigr)$} 
  \end{algorithmic}
\end{algorithm}

Whether the individual steps of Algorithm~\ref{alg:swa} can be made effective
depends on how $\M$ is given.  If it is presented by generators and relations,
as we always assume, one may use Gr\"obner bases; Sturmfels and White
formulated their algorithm directly for this case.  Note that they need
repeated Gr\"obner bases calculations in order to perform algorithmically all
the computations in factor modules.  A further problem is to find the non zero
divisors, as Lemma~\ref{lem:nzd} only guarantees their existence but says
nothing about their determination.  Sturmfels and White proposed a
probabilistic approach.  As the non zero divisors contain a Zariski open
subset of~$\P_1$, random choice of an element of $\P_1$ yields one with
probability~$1$.

\begin{theorem}\label{thm:swa}
  The Algorithm~\ref{alg:swa} terminates for any finitely generated polynomial
  module\/~$\M$ with a Rees decomposition.
\end{theorem}

For a proof we refer to~\cite{bg:comb,sw:comb} where also Lemma~\ref{lem:nzd}
is proved.  We will show that generically the Sturmfels-White
Algorithm~\ref{alg:swa} returns a Pommaret basis when applied to a submodule
of a free module.  We begin by studying the relation between the classes and
the coordinate classes of the generators of a Rees decomposition determined
with this algorithm.

\begin{proposition}\label{prop:cls1}
  Let\/ $\H$ define a Rees decomposition of the form (\ref{eq:rd}) for the
  submodule\/ $\M\subseteq\P^r$ with respect to the basis\/ $\Y$ of\/ $\P_1$
  and let\/ $\H$ and\/ $\Y$ be determined with the Sturmfels-White Algorithm
  \ref{alg:swa}.  With respect to the basis\/ $\Y$ we have the inequalities\/
  $\cls{\hv}\geq\ccls{\hv}{\yv}$ for all\/ $\hv\in\H$.
\end{proposition}

\begin{proof}
  $Z_{\M}=0$ for a submodule $\M\subset\P^r$.  Thus Algorithm~\ref{alg:swa}
  produces no generators of class~$0$.  The coordinate class is always greater
  than~$0$.
  
  We follow step by step Algorithm~\ref{alg:swa}.  In the first iteration some
  non zero divisor $y_1\in\P_1$ is chosen and in the second iteration we must
  treat the factor module $\M^{(1)}=\M/y_1\M$.  Now $\fv\in\M$ represents an
  element of $Z_{\M^{(1)}}$, if and only if $y_k\fv\in{y_1\M}$ for all $k>1$.
  Thus $\ccls{(y_k\fv)}{\yv}=1$ for all $k>1$ which is only possible if
  $\ccls{\fv}{\yv}=1$.
  
  If $Z_{\M^{(1)}}\neq0$, Algorithm~\ref{alg:swa} proceeds with
  $\M^{(2)}=\M^{(1)}/Z_{\M^{(1)}}$.  Now $\fv\in\M$ represents an element of
  $Z_{\M^{(2)}}$, if and only if for all $k>1$ the product $y_k\fv$ either is
  an element of $y_1\M$ or represents an element of~$Z_{\M^{(1)}}$.  In both
  cases this is only possible, if $\ccls{\fv}{\yv}=1$.  The same argument
  holds until $Z_{\M^{(\ell)}}=0$ for some~$\ell$.  Thus all generators~$\hv$
  to which Algorithm~\ref{alg:swa} assigns the class~$1$ are divisible by
  $y_1$ and hence all their terms possess the coordinate class~$1$.
  
  The next module to consider is $\M^{(\ell+1)}=\M^{(\ell)}/y_2\M^{(\ell)}$.
  Proceeding as above we see that $\fv\in\M$ represents an element of
  $Z_{\M^{(\ell+1)}}$, if and only if $y_k\fv\in{y_2\M^{(\ell)}}$ for all
  $k>2$ implying that $\ccls{\fv}{\yv}\leq2$.  Using the same argument as
  above, we conclude that all generators of class~$2$ according to the
  Sturmfels-White algorithm consist of terms with a coordinate class less than
  or equal to~$2$.  Following Algorithm~\ref{alg:swa} until the end we obtain
  the assertion, namely that $\cls{\hv}\geq\ccls{\hv}{\yv}$ for all
  generators~$\hv\in\H$.\qed
\end{proof}

As it may happen that $\cls{\hv}>\ccls{\hv}{\yv}$ for some generator $\hv$,
the set $\H$ is not necessarily a Pommaret basis.  More precisely, the
coordinates $\yv$ are not necessarily $\delta$-regular for $\H$.  We show now
similarly to the proof of Theorem~\ref{thm:dreg} that we may always transform
$\yv$ into a $\delta$-regular coordinate system $\zv$.

For simplicity, let us assume that only one generator $\hv$ with
$\ccls{\hv}{\yv}<\cls{\hv}$ exists and that $\cls{\hv}=2$.  Consider a
coordinate transformation $z_k=y_k$ for $k>1$ and $z_1=y_1+cy_2$ where
$c\in\kk$ is chosen such that with respect to the new coordinates
$\ccls{\hv}{\zv}=2$.  The possible values of $c$ form a Zariski open set
in~$\kk$.  By Lemma~\ref{lem:nzd}, the non zero divisors among which $y_1$ was
chosen in Algorithm~\ref{alg:swa} contain a Zariski open subset of~$\P_1$.
Thus there exist values of $c$ such that both $\ccls{\hv}{\zv}=2$ \emph{and}
$z_1$ is a non zero divisor.

As in the proof of Theorem~\ref{thm:dreg}, it is not difficult to show that
this transformation increases the involutive size of the set $\H$.  Applying a
finite number of similar changes of coordinates leads to a new basis $\Z$ of
$\P_1$ in which $\ccls{\hv}{\zv}=\cls{\hv}$.  As we still have a Rees
decomposition, $\H$ is a Pommaret basis of $\M$ and the coordinates $\zv$ are
$\delta$-regular.  Obviously, the one-forms $z_i$ would have been valid
choices for the non zero divisors in Algorithm~\ref{alg:swa}.  Thus we
conclude that this algorithm may be used for the construction of Pommaret
bases.  The following proposition shows that in fact any Pommaret basis may be
constructed this way.
  
\begin{proposition}\label{prop:cls2}
  Let\/ $\H$ be a Pommaret basis of the submodule\/ $\M\subseteq\P^r$ with
  respect to the $\delta$-regular coordinates\/ $\yv$ and a class respecting
  term order~$\prec$.  The one-forms\/ $y_i$ may be used as non zero divisors
  in Algorithm~\ref{alg:swa} and the there obtained generators~$\bar\hv$
  satisfy $\cls{\bar\hv}=\ccls{\bar\hv}{\yv}$.  They are\/ $\kk$-linear
  combinations of the elements of\/ $\H$; one may even use the elements of\/
  $\H$ as generators.
\end{proposition}

\begin{proof}
  The Pommaret basis  $\H$ defines a Rees decomposition 
  \begin{equation}\label{eq:prd}
    \M=\bigoplus_{\hv\in\H}\kk[y_1,\dots,y_{\ccls{\hv}{\yv}}]\cdot\hv\,.
  \end{equation}
  As in the proof of the Proposition~\ref{prop:cls1}, we follow step by step
  the Sturmfels-White Algorithm~\ref{alg:swa}.  Let $\M^{(1)}=\M/y_1\M$ and
  $\H_1=\{\hv\in\H\mid\ccls{\hv}{\yv}=1\}$.  The vector space $Z_{\M^{(1)}}$
  is isomorphic to a subspace of the vector space freely generated by $\H_1$,
  as by Proposition~\ref{prop:cls1} $Z_{\M^{(1)}}$ contains only elements with
  coordinate class~$1$ and the only elements of $\M$ of coordinate class~$1$
  which are not in $y_1\M$ are $\kk$-linear combinations of the elements
  of~$\H_1$.
  
  Let $\hv\in\H_1$ and $k>1$.  We determine the involutive normal form of
  $y_k\hv$ induced by (\ref{eq:prd}).  Every term in $y_k\hv$ has coordinate
  class~$1$, thus $t=\lt{\prec}{(y_k\hv)}$ satisfies $\ccls{t}{\yv}=1$.  Since
  $\H$ is a strong basis, there exists precisely one generator $\hv'\in\H$
  such that $\lt{\prec}{\hv'}\idiv{P}t$.  If $\lt{\prec}{\hv'}=t$, then
  $\hv'\in\H_1$, as $\prec$ is class respecting.  After the corresponding
  reduction step, the initial term of the result is still of coordinate
  class~$1$.  So the normal form of $y_k\hv$ has the following structure
  \begin{equation}\label{eq:nf}
    y_k\hv=\sum_{\tilde\hv\in\H_1}c_{\tilde\hv}\tilde\hv + y_1\fv
  \end{equation}
  for some coefficients $c_{\tilde\hv}\in\kk$ and an element $\fv\in\M$.  The
  vector space $Z_{\M^{(1)}}$ is generated by those $\tilde\hv\in\H_1$ where
  the first summand in (\ref{eq:nf}) is zero; this includes all elements of
  $\H_1$ of maximal degree.
  
  Thus if $\H_1\neq\emptyset$, then $Z_{\M^{(1)}}\neq\emptyset$.
  Algorithm~\ref{alg:swa} proceeds in this case with
  $\M^{(2)}=\M^{(1)}/Z_{\M^{(1)}}$.  If $\dim Z_{\M^{(1)}}<|\H_1|$, then
  $Z_{\M^{(2)}}\neq\emptyset$.  Algorithm~\ref{alg:swa} will iterate
  Line~/14/, until all elements of $\H_1$ have been used up.  When this stage
  is reached, $Z_{\M^{(\ell)}}=0$.  It follows from the direct sum in
  (\ref{eq:prd}) that $y_2$ is a non zero divisor for $\M^{(\ell)}$ and we may
  proceed with $\M^{(\ell+1)}=\M^{(\ell)}/y_2\M^{(\ell)}$.
  
  Let $\H_2=\{\hv\in\H\mid\ccls{\hv}{\yv}=2\}$.  As above, $Z_{\M^{(\ell+1)}}$
  is isomorphic to a subspace of the vector space freely generated by~$\H_2$.
  There are some minor modifications in (\ref{eq:nf}): the first sum is over
  all $\tilde\hv\in\H_2$ and there are additional summands which vanish either
  modulo $y_1\M$ or modulo $y_2\M^{(\ell)}$ or modulo some $Z_{\M^{(i)}}$ for
  $1\leq i\leq\ell$.  Again Algorithm~\ref{alg:swa} will iterate line~/14/,
  until all elements of $\H_2$ have been used up.  The same argument may be
  repeated for $y_3,\dots,y_n$.
  
  Thus, Algorithm~\ref{alg:swa} terminates with a Rees decomposition (with
  respect to the basis $\Y\subset\P_1$) generated by $\bar\H$ where
  $|\bar\H|=|\H|$ and where the elements $\bar\hv\in\bar\H$ with
  $\cls{\bar\hv}=k$ freely generate the same vector space as the elements
  $\hv\in\H$ with $\ccls{\hv}{\yv}=k$.  We may even choose $\bar\H=\H$.  In
  any case, $\cls{\bar\hv}=\ccls{\bar\hv}{\yv}$ and $\bar\H$ is a Pommaret
  basis of $\M$.\qed
\end{proof}

\begin{remark}
  Note the strong similarity between this proof and the proof of
  Proposition~\ref{prop:Iregular}.  This is not surprising, as the minimal
  class assigned by the algorithm is equal to $\depth{\M}$ \cite{sw:comb} and
  the basis $\Y$ determined by Algorithm~\ref{alg:swa} is quasi-regular for
  the module $\M$ in the sense of Serre (see the letter of Serre appended to
  \cite{gs:alg} or \cite{wms:spencer2}).  In fact, Lemma \ref{lem:nzd} follows
  immediately from the results of Serre.  They imply furthermore that
  $\Z_{\M}$ is always finite-dimensional and thus it is not really necessary
  to factor by $\Z_{\M}$.  In \cite{wms:spencer2} it is shown that coordinates
  are $\delta$-regular for the a submodule $\M\subseteq\P^r$, if and only if
  they are quasi-regular for the factor module $\P^r/\M$.  Thus in principle,
  one should always compare Algorithm~\ref{alg:swa} applied to $\P^r/\M$ with
  the Pommaret basis of $\M$ (recall that the latter also leads immediately to
  a Rees decomposition of $\P^r/\M$ via Corollary \ref{cor:stacon}).\bull
\end{remark}

\begin{acknowledgement}
  The author would like to thank V.P.~Gerdt for a number of interesting
  discussions on involutive bases.  M.~Hausdorf and R.~Steinwandt participated
  in an informal seminar at Karlsruhe University where many ideas of this
  article were presented and gave many valuable comments.  The constructive
  remarks of the anonymous referees were also very helpful.  This work
  received partial financial support by Deutsche Forschungsgemeinschaft, INTAS
  grant 99-1222 and NEST-Adventure contract 5006 (\emph{GIFT}).
\end{acknowledgement}

\bibliography{../../BIB/DiffEq,../../BIB/Algebra,../../BIB/Seiler,%
../../BIB/Groebner,../../BIB/DiffAlg,../../BIB/Misc,../../BIB/Lie,%
../../BIB/CA,../../BIB/DiffGeo,../../BIB/ZZProc}
\bibliographystyle{plain}

\end{document}